\documentclass[smallcondensed]{svjour3}    
\smartqed  
\usepackage{graphicx}\graphicspath {{figs/}}
\usepackage{amsmath,amssymb}
\usepackage[pdfborder={0 0 0}]{hyperref}
\usepackage{siunitx}
\usepackage{xcolor}

\renewcommand{\Re}{\operatorname{Re}}

\renewcommand{\L}{\mathcal{L}}
\newcommand{\N}{\mathcal{N}}

\newcommand{\scaledp}{\rho}

\newcommand{\abs}[1]{\left|#1\right|}
\newcommand{\norm}[1]{\left\lVert#1\right\rVert}

\newcommand{\qand}{\quad\text{and}\quad}

\begin{document}
\title{Linearly Stabilized Schemes for the Time Integration of Stiff Nonlinear PDEs
}

\author{Kevin Chow   \and    Steven J. Ruuth 
}

\institute{K. Chow \at
Department of Mathematics,
              Simon Fraser University,
             8888 University Drive,
	  Burnaby, BC V5A 1S6 Canada. \\
	                \email{kjc19@sfu.ca}           
           \and
           S.J. Ruuth \at
              Department of Mathematics,
              Simon Fraser University,
             8888 University Drive,
	  Burnaby, BC V5A 1S6 Canada. \\
	 \email{sruuth@sfu.ca}
}

\date{}

\maketitle

\begin{abstract}
	In many applications, the governing PDE to be solved numerically contains a stiff component. When this component is linear, an implicit time stepping method that is unencumbered by stability restrictions is often preferred. On the other hand, if the stiff component is nonlinear, the complexity and cost per step of using an implicit method is heightened, and explicit methods may be preferred for their simplicity and ease of implementation. In this article, we analyze new and existing linearly stabilized schemes for the purpose of integrating stiff nonlinear PDEs in time. These schemes compute the nonlinear term explicitly and, at the cost of solving a linear system with a matrix that is fixed throughout, are unconditionally stable, thus combining the advantages of explicit and implicit methods. Applications are presented to illustrate the use of these methods.
\end{abstract}

\section{Introduction}
\label{sect: intro}
In this paper, we propose and analyze some new linearly stabilized schemes for the time integration of stiff nonlinear PDEs. The linearly stabilized semi-implicit Euler scheme (see Sect.~\ref{sssect: first order scheme}) is a first order scheme of this type that has been used to approximate the solutions to a variety of PDE problems.   Its first known use appears in a paper by Douglas and Dupont \cite{douglas1971alternating} where it was applied to a variable coefficient heat equation on rectangular domains. In subsequent years, the idea has been rediscovered by Eyre \cite{eyre1998bunconditionally}, who first used the name ``linearly stabilized'', and Smereka \cite{smereka2003semi}. Others have gone on to apply these schemes to Hele-Shaw flows, interface motion, image processing, and solving PDEs on surfaces \cite{eyre1998bunconditionally,salac2008local,glasner2002diffuse,schonlieb2011unconditionally,macdonald2009implicit}.

In each of the references mentioned above, the authors have implemented only a first order time stepping method. More recently, Duchemin and Eggers \cite{duchemin2014explicit} consolidated the approach and produced a second order linearly stabilized scheme they refer to as the explicit-implicit-null (EIN) method. Their method attains second order accuracy by extrapolating the first order results. Moreover, they identified that a key principle for the success of any linearly stabilized scheme is unconditional stability. They show that their method is unconditionally stable under only a mild condition on a parameter that is introduced.

Our derivations for new linearly stabilized schemes will also begin by ensuring that the newly derived schemes are unconditionally stable. The techniques we employ in our stability analysis are those of a standard linear stability analysis, but are applied to a modified test equation. In Sect.~\ref{sect: linear mod}, we formally introduce the notion of linear stabilization. Motivation for this technique is supplied by the need to handle a stiff nonlinear PDE describing axisymmetric mean curvature flow and leads us to the well-studied first order linearly stabilized scheme and the EIN method of Duchemin and Eggers. Following that, the framework in which we analyze the stability of linearly stabilized schemes is introduced.   
\textcolor{black}{
A notable property of this approach to analyzing stability is its assumption that the two operators appearing in the formulation are simultaneously diagonalizable.   Interesting recent work \cite{rosales2017unconditional,seibold2019unconditional} develops an analogous concept that does not require this assumption.   See \cite{rosales2017unconditional,seibold2019unconditional} for details on the approach as well as a corresponding new class of unconditionally stable linear multistep IMEX schemes.
}

In Sect.~\ref{sect: lmm}, we investigate implicit-explicit (IMEX) linear multistep methods within the linear stabilization framework.  A detailed comparison of the schemes based on IMEX methods and the EIN method is conducted in Sect.~\ref{sect: 3 key properties}. Our experiments suggest that three criteria, in addition to unconditional stability, are desired for practical linearly stabilized schemes. Notably, one of these criteria eliminates third and higher order multistep-based linearly stabilized schemes from use and another suggests EIN is far from practical either. 

In Sect.~\ref{sect: exp rk}, we explore the use of exponential Runge-Kutta methods to mend this deficiency. A second order and a fourth order exponential Runge-Kutta method are verified to exhibit unconditional stability over an unbounded parameter range. However, the error constant of both these schemes scales unfavorably in $p$, and this narrows their range of applicability.

In Sect.~\ref{sect: numerical experiments}, application of our linearly stabilized schemes to a number of 2D and 3D problems is presented. Not surprisingly, our second order schemes offer improvements over the commonly used first order linearly stabilized scheme. The experiments show that our schemes provide a substantial efficiency improvement yet the complexity of their implementation is no greater than solving a heat equation with standard implicit methods.

Finally, some concluding remarks are presented in Sect.~\ref{sect: conc}.  

\section{Linear Modification and Unconditional Stability}
\label{sect: linear mod}
To construct time integration schemes for stiff, nonlinear PDEs, we set out two key design principles. Firstly, we want to handle the nonlinearity simply and inexpensively. Secondly, we must be free to select time step-sizes reflecting the accuracy requirement, rather than choosing step-sizes that are severely constrained by stability. Linearly stabilized schemes, as we will see, adhere to both principles and are remarkably easy to implement.

\subsection{Prototype 1D Problem}
\label{subsect: proto 1d}
\begin{figure}[htb!]
	\centering
	\includegraphics[width=84mm]{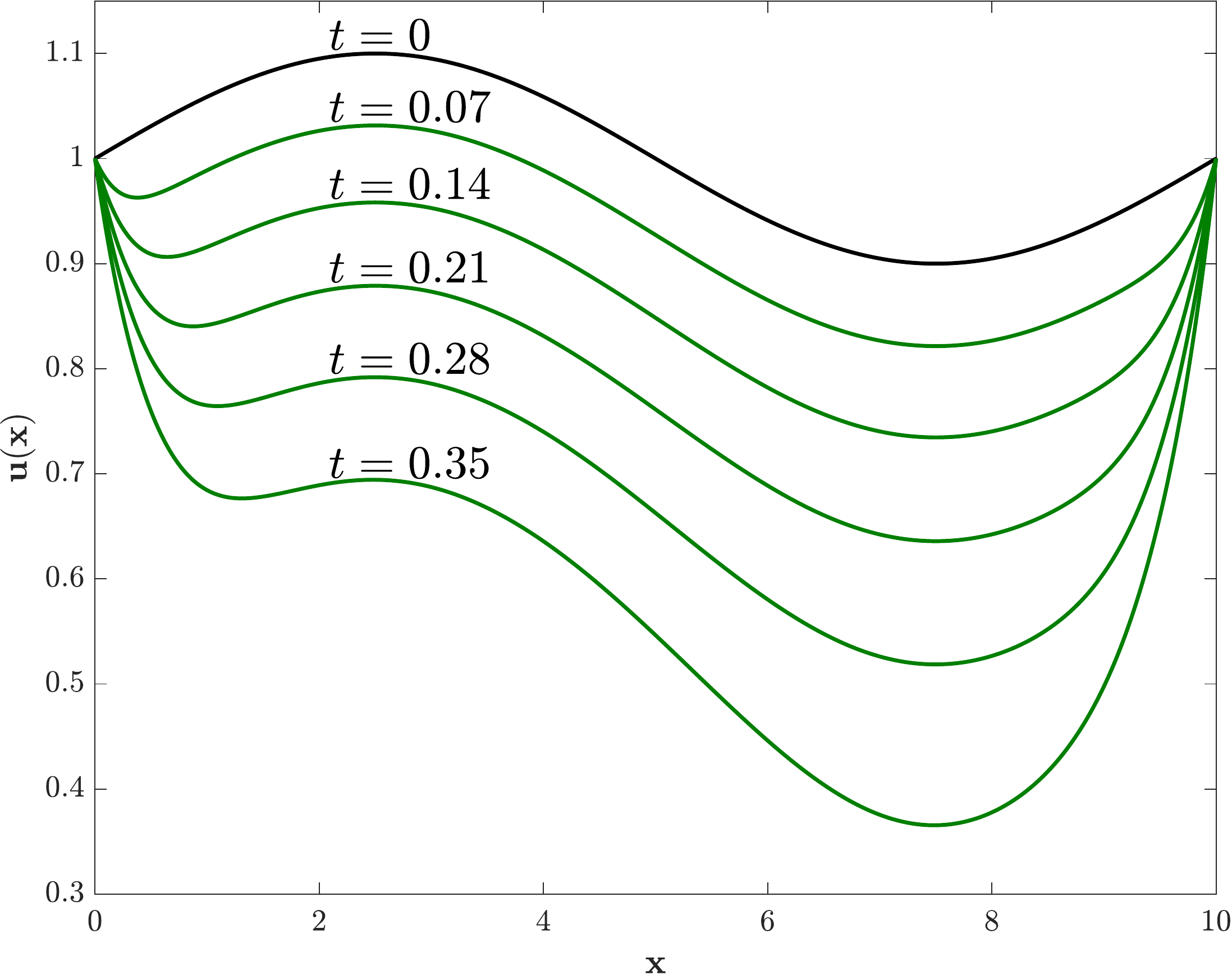}
	\caption[Numerical solution to a 1D axisymmetric mean curvature motion problem]{Time evolution for the solution to \eqref{ammc}. The initial state is a sinusoidal perturbation from an unstable equilibrium \cite{bernoff1998axisymmetric}. The bottommost curve is the final state ($T=0.35$).}
	\label{fig:ammc sol}
\end{figure}

As a prototype, let us consider the following 1D axisymmetric mean curvature motion problem \cite{duchemin2014explicit}: 
\begin{subequations} 
	\begin{align}
		u_t &= \frac{u_{xx}}{1 + u_x^2} - \frac{1}{u}, 
		\qquad 0 < x < 10, \quad t > 0,
		\label{ammc 1}
	\end{align}
	with initial and boundary conditions
	\begin{gather}
		u(x,0) = 1 + 0.10\sin\left( \frac{\pi}{5}x \right),
		\\
		u(0,t)=u(10,t)=1.
	\end{gather}
	\label{ammc}
\end{subequations}
A time evolution of this problem is plotted in Fig.~\ref{fig:ammc sol}.

The presence of the $u_{xx}$ guarantees that \eqref{ammc} is stiff, suggesting that an implicit time stepping scheme may prove to be more efficient. However, instead of a static linear system which can be preprocessed and solved efficiently, this would require solution to a nonlinear system at each time step due to the factor of $(1 + u_x^2)^{-1}$. Thus we are presented with a scenario where neither an implicit nor an explicit approach proves particularly palatable.

\subsubsection{A first order linearly stabilized scheme}
\label{sssect: first order scheme}
As demonstrated in Duchemin and Eggers \cite{duchemin2014explicit} as well as in an earlier paper by Smereka \cite{smereka2003semi}, an efficient method for handling \eqref{ammc} is to add and subtract a linear Laplacian term to the right-hand side, 
\begin{align}
	u_t = \underbrace{\frac{u_{xx}}{1 + u_x^2} 
		- \frac{1}{u} 
		- u_{xx}}_{\N(u)} 
	+ \underbrace{\phantom{\frac{1_1}{1}}u_{xx}\phantom{\frac{1}{1_1}}}_{\L u}, 
	\label{ammc +- uxx}
\end{align} 
and then time step according to
\begin{align}
	\frac{u^{n+1} - u^n}{\Delta t} 
	= \N(u^n) + \L u^{n+1}.
	\label{ammc sbdf1}
\end{align}
Since this is our first instance of a linearly stabilized scheme, we remark on some of the key properties. We first note that in the continuous case, the modified equation \eqref{ammc +- uxx} is unchanged from \eqref{ammc 1}. Next, note that once discretized \eqref{ammc sbdf1}, the nonlinear term is evaluated explicitly; ignoring the $\L u^{n+1}$ term, it corresponds to a forward Euler step. 
On the other hand, it is a step of backward Euler if only the linear terms are considered. This combination of time stepping methods is known as implicit-explicit (IMEX) or semi-implicit Euler \cite{ascher1995implicit,smereka2003semi}. As it is a combination of explicit and implicit Euler steps, the accuracy is first order. We also note that discretization of the Laplacian term 
typically leads to a sparse, symmetric, positive definite matrix.  Inversion of the corresponding system is efficient and easy to implement. Lastly, as a result of the implicit discretization of the $\L u$ term, we may expect this scheme to have improved stability compared to a purely explicit scheme, and indeed this is the case. Discretizing with second order centered differences in space, it can be shown by von Neumann analysis that this scheme is unconditionally stable \cite{duchemin2014explicit}.

\subsubsection{Second order by Richardson extrapolation}
\label{sssect: second order richardson}
As stated at the outset, the time stepping procedure in \eqref{ammc sbdf1} is only first order. The work of Duchemin and Eggers \cite{duchemin2014explicit} extends the method to second order by Richardson extrapolation (see also \cite{smereka2003semi}, where Richardson extrapolation was suggested but not implemented). They generalized the approach with a free parameter, $p$, i.e.,
\begin{align}
	u_t = \frac{u_{xx}}{1 + u_x^2} - \frac{1}{u} - pu_{xx} + pu_{xx}, 
	\label{ammc p}
\end{align}
and derived restrictions on $p$ subject to the condition that the resulting scheme be unconditionally stable. \textcolor{black}{With the semi-implicit Euler approach~(\ref{ammc sbdf1}) and $N$ spatial grid nodes, they found $p \geq \max_{1\leq j\leq N} 0.5/(1 + (D_1 u^n_j)^2)$, where $D_1 u^n_j$ is the second order centered difference approximation to $u_x$ at $x_j$, to be sufficient when stabilized as in \eqref{ammc p}. With an additional Richardson extrapolation step, the restriction becomes $p \geq \max_{1\leq j \leq N} (2/3)/(1 + (D_1 u^n_j)^2)$.
}
\subsection{A Modified Test Equation}
\label{subsect: mte}
\textcolor{black}{
Section~\ref{subsect: proto 1d}  gives an example where a specific problem is discretized and then analyzed for stability via a von Neumann analysis. However, it is frequently the case that we would like to know the stability properties of a scheme in a more standalone fashion. Analogous to the standard linear stability analysis where a numerical scheme is applied to the test equation $w' = \lambda w$, we wish to establish the stability of a linear stabilization scheme with respect to a suitable test equation. }

\textcolor{black}{
With linear stabilization, the nonlinear system $u' = F(u)$ is modified according to
\begin{align}
u' = \underbrace{ (F(u) - pL u)}_{\N(u)} + \underbrace{pL u}_{\L u},
\label{generic} 
\end{align}
where $p > 0$.   To analyze stability, we first linearize, $F(u) \approx F(u^*) +  J (u-u^*)$, where  $J \equiv J_F(u^*)$ denotes the Jacobian of $F$ at the expansion point $u^*$.  Neglecting higher order terms  leads us to the linear equation
\begin{align}
v' = Jv - pLv + pLv.
\end{align}
Assuming that $J, L$ are simultaneously diagonalizable simplifies the problem to one involving scalar equations
\begin{align}
w' = \lambda_F w - p\lambda_L w + p\lambda_L w,
\end{align}
which we reformulate as 
\begin{align}
w' = (1 - \bar p) \lambda w + \bar p \lambda w,
\label{mte}
\end{align}
where $\lambda = \lambda_F$ and $\bar p = p\lambda_L / \lambda_F$.} We will refer to \eqref{mte} as the modified test equation. Note that when $p = 0$, the modified test equation reduces to the standard test equation.  The real $\lambda$ case arises frequently and will be the relevant case for the applications we consider in this paper.
As a consequence, we are primarily interested in the case where $\lambda<0$ is real.
However, complex $\lambda_F$, $\lambda_L$ with
\begin{align*}
\Re(\lambda_F), \Re(\lambda_L) \leq 0, \quad
\arg(\lambda_F) = \arg(\lambda_L), \quad 
p > 0, 
\end{align*}
can also  arise (cf. \cite{duchemin2014explicit}) and is amenable to our analysis.  We therefore present results for that case as well.

\textcolor{black}{
The assumptions made here regarding $J,L$ lead to a relatively straightforward analysis by way of the modified test equation. However, one could ask whether these assumptions may be relaxed. To that end, the paper by Rosales et al. \cite{rosales2017unconditional} and its companion by Seibold et al. \cite{seibold2019unconditional} develop an analogous concept and attain unconditional stability for a new class of linear multistep IMEX schemes, without the assumption that the two operators $J, L$ commute. See \cite{rosales2017unconditional,seibold2019unconditional} for details on the stability criteria, new schemes up to fifth order, as well as a variety of illustrative and illuminating numerical experiments.
}

We discuss next the stability properties of three time stepping methods as applied to the modified test equation \eqref{mte}.

\subsubsection{Forward Euler}
\label{sssect: fe}
Forward Euler is a first order time stepping method that treats the right-hand side explicitly. Application to \eqref{mte} is therefore no different than to the standard test equation. 
As a consequence, we cannot obtain unconditional stability.

\subsubsection{Linearly stabilized semi-implicit Euler}
\label{sssect: lin stab sbdf1}
Semi-implicit Euler time stepping~\eqref{ammc sbdf1} was applied to the 1D axisymmetric mean curvature motion problem \eqref{ammc}, and its stability analyzed in \cite{smereka2003semi,duchemin2014explicit}. For the modified test equation \eqref{mte}, we identify $\N(w^n) = (1-\bar p)\lambda w^n$ and $\L w^{n+1} = \bar p\lambda w^{n+1}$, to get 
\begin{align}
	\frac{w^{n+1} - w^n}{\Delta t} 
	= (1 - \bar p)\lambda w^n + \bar p\lambda w^{n+1} 
	\iff 
	w^{n+1} 
	= \xi_E w^n,
\end{align}
where 
\[
\xi_E = \left(1 + \frac{\lambda \Delta t}{1 - \bar p\lambda\Delta t} \right).
\]
Enforcing unconditional stability, i.e. $\abs{\xi_E} \leq 1$, for all $\lambda \Delta t \leq 0$, we find 
\begin{align}
	\abs{\xi_E} \leq 1 \iff 
	-2 \leq \frac{\lambda\Delta t}{1 - \bar{p}\lambda\Delta t} \le 0
	\iff \bar{p} \geq 0 \quad\text{and}\quad (2\bar{p}-1)\lambda\Delta t \leq 2.
\end{align}
Thus unconditional stability is guaranteed if $\bar{p} \geq 1/2$.  The same bound arises for the complex $\lambda$ case.

Going forward, we shall refer to the linearly stabilized semi-implicit Euler method as SBDF1.

\subsubsection{Explicit-implicit-null}
\label{sssect: ein}
In \cite{duchemin2014explicit}, the SBDF1 approach is extended to second order by using Richardson extrapolation, and their methodology is referred to as explicit-implicit-null (EIN). For EIN, the amplification factor, $\xi_{\mathrm{EIN}}$, can be expressed in terms of $\xi_E$, 
\begin{align}
	\xi_\mathrm{EIN}
	= 2\xi^2_{E}(\lambda\Delta t/2) - \xi_E(\lambda\Delta t) 
	= 1 + \frac{z\left(\bar{p}( 3\bar{p}-2 )z^2 + 2( 1-4\bar{p} )z + 4 \right)}{( 1 - \bar{p}z )( 2-\bar{p}z )^2}
	\label{ein amp fac}
\end{align} 
where $z = \lambda \Delta t$.
They show that unconditional stability is guaranteed if $\bar{p} \geq 2/3$.  Once again,  the same bound arises for the complex $\lambda$ case.

\begin{example}
Suppose we wish to determine $p$ for the EIN method applied to the 1D axisymmetric mean curvature motion problem \eqref{ammc}.  
For a centered difference spatial discretization, the relevant eigenvalues are \cite{duchemin2014explicit}
	\begin{align}
		\lambda_F = \frac{2}{\Delta x^2}\frac{\cos(k\Delta x)-1}{1 + (D_1 u_j^n)^2}
		\qand 
		\lambda_L = \frac{2}{\Delta x^2}(\cos(k\Delta x) - 1).
		\label{ammc eig}
	\end{align}
	\textcolor{black}{
Since $\bar{p} \geq 2/3$, and $\bar p = p\lambda_L / \lambda_F$, we require
	\begin{align}
		p \geq \max_{1\leq j \leq N} \frac{2}{3}(1 + (D_1 u^n_j)^2)^{-1}.
		\label{ammc ein p}
	\end{align} }
Notice that the restriction on $p$ varies in time. Although a time-adaptive and/or space-adaptive approach may be possible, we do not explore that here. A constant value of $p$ is set to satisfy the time stepping scheme's restriction throughout the evolution of the system (for example, in \cite{duchemin2014explicit} the numerical experiments were reported with $p = 0.7$).

\end{example}

\section{IMEX Linear Multistep Methods}
\label{sect: lmm}
For equations whose right-hand side is comprised of a stiff linear component and a nonstiff nonlinear part, a popular class of methods to apply are the implicit-explicit linear multistep methods\footnote{We will refer to these simply as IMEX methods.}. The simplest of these is the SBDF1 scheme that we reviewed in Sect.~\ref{sssect: lin stab sbdf1}.

In this section, we investigate the use of selected second, third and fourth order IMEX methods within the context of linearly  stabilized schemes. 
In our approach, the added linear term
will be discretized implicitly, while the remaining terms, including the stiff nonlinear term, will be treated explicitly.

\subsection{IMEX Formulas}
\label{subsect: imex formulas}
In \cite{ascher1995implicit}, IMEX schemes up to order four are investigated and a select number are singled out for their extensive use in the literature or for desired properties such as strong high frequency damping. We present these schemes relative to
the ODE system
\begin{align*}
	u' = f + g, 
\end{align*}
where standard usage has $g$ representing a stiff linear term, and $f$ representing the remaining nonlinear/nonstiff terms.   
Recognizing that the standard, first order IMEX scheme is simply SBDF1, we proceed immediately to second order methods, with higher order methods following. 
\\ \ \\ \noindent
\textbf{Second order methods}
\\
Second order IMEX schemes that have appeared in the literature include the following three schemes: 
\\ \ \\
CNAB:
\begin{align}
	\frac{u^{n+1}-u^n}{\Delta t} 
	= \frac{3}{2} f^n - \frac{1}{2}f^{n-1} 
	+ \frac{1}{2}(g^{n+1} + g^n), 
	\label{cnab}
\end{align}
\\
CNLF:
\begin{align}
	\frac{u^{n+1}-u^{n-1}}{2\Delta t}
	= f^n + \frac{1}{2}(g^{n+1} + g^{n-1}),
	\label{cnlf}
\end{align} 
\\
SBDF2:
\begin{align}
	\frac{3u^{n+1}-4u^n+u^{n-1}}{2\Delta t} 
	= 2f^n - f^{n-1} + g^{n+1}.
	\label{sbdf2}
\end{align}

CNAB combines the second order Crank-Nicholson and Adams-Bashforth schemes. It has a small error constant, but gives a slow decay of high frequency
error modes in dissipative problems.  
\footnote{To address this deficiency, \cite{ascher1995implicit} recommended mCNAB, 
a scheme closely related to CNAB but with stronger damping of high frequencies. As it turns out, the two are equivalent within this linear stabilization framework.}
CNLF also combines two well-known schemes, Crank-Nicholson and Leap Frog.  In our context, this scheme does not appear particularly promising. 
Relative to CNAB or SBDF2, it has a small stability region and offers a very weak decay of high frequency error modes in  dissipative problems.  
Finally, SBDF2 has as its implicit part the second order BDF scheme.  It has the attractive feature of strongly damping high frequency error components. 
\\

\noindent
\textbf{Higher order methods}
\\
Third and fourth order IMEX schemes based on the corresponding BDF schemes have also been used: 
\\ \ \\
SBDF3:
\begin{align}
	\frac{1}{\Delta t}\left(\frac{11}{6}u^{n+1} - 3u^n + \frac{3}{2}u^{n-1} - \frac{1}{3}u^{n-2} \right) 
	= 3f^n - 3f^{n-1} + f^{n-2} + g^{n+1},
	\label{sbdf3}
\end{align}
\\
SBDF4:
\begin{align}
	\frac{1}{\Delta t}\left(\frac{25}{12}u^{n+1} - 4u^n + 3u^{n-1} - \frac{4}{3}u^{n-2} + \frac{1}{4}u^{n-3} \right) 
	= 4f^n - 6f^{n-1} + 4f^{n-2} -f^{n-3} + g^{n+1}.
	\label{sbdf4}
\end{align}
Similar to their lower order counterparts, SBDF3 and SBDF4  produce a strong decay of high frequency error in dissipative problems.
\\

In the next subsections, these IMEX schemes will be applied to the modified test equation to determine for each scheme the range of $\bar{p}$ suitable for linear stabilization.

\subsection{Analysis of the Amplification Polynomials}
\label{subsect: analysis amp poly}
Before we apply the above IMEX schemes to the modified test equation \eqref{mte}, 
let us observe that the polynomial arising from the application of an $n$th order IMEX method will be a degree $n$ polynomial in the amplification factor, $\xi$. 
The goal is to identify for each IMEX scheme any restriction on the parameter $\bar{p}$ which when satisfied will allow a user to freely choose the time step-size without being subject to a stability constraint. 
For the analysis of these amplification polynomials, we turn to the theory of von Neumann polynomials \cite[Chapter 4]{strikwerda2004finite}. 
In particular, we require the resulting amplification polynomials to be simple von Neumann.  

Two-step, second order IMEX methods form a two parameter family \cite{ascher1995implicit}.  
Interestingly,  when applied to \eqref{generic}, the number of free parameters reduces by one.  
To see this, we apply the general two-parameter family to (\ref{generic}), yielding  
\begin{align}\begin{split} 
\frac{1}{\Delta t}&\left[ 
	\left(\gamma + \frac{1}{2}\right) u^{n+1} - 2\gamma u^n + \left( \gamma - \frac{1}{2}\right) u^{n-1}
	\right] 
\\&=
	(\gamma + 1)F(u^n) - \gamma F(u^{n-1}) +  \scaledp L u^{n+1} -2\scaledp Lu^n + \scaledp L u^{n-1},
	\end{split}
	\label{2order1param}
\end{align}
where $\gamma \in [0,1]$ is a free parameter,
and $\scaledp$ is a scaled version of $p$ that has absorbed the second parameter. 
Simplifying \eqref{2order1param} to the case of the modified test equation gives the amplification polynomial 
\begin{align}
	\left( \gamma + \bar{\scaledp} z + \frac{1}{2} \right) \xi^2 
-	\left( 2\gamma + z(\gamma + 1 - 2\bar{\scaledp}) \right) \xi 
+	\gamma - \frac{1}{2} + z(\gamma - \bar{\scaledp})
= 0,
\end{align}
where  $\bar{\scaledp} \equiv  (\lambda_F / \lambda_L)\cdot \scaledp \geq (1 + 2\gamma)/4$ is necessary and sufficient for unconditionally stability for both the real and complex-valued $\lambda$ cases.
We present the corresponding range of stable $\bar{p}$ values for CNAB, CNLF, and SBDF2 in Table~\ref{table:amp poly 2}.  The minimal value required for unconditional
stability (1, 1/2,  3/4 for CNAB, CNLF and SBDF2, respectively) will be used later in our numerical experiments and will be denoted by $\bar{p}_{\text{min}}$.  
The computer algebra system, \textsc{Maple}\texttrademark, was used to facilitate the calculations.

\begin{table}
	\caption{Amplification polynomials for select second order linearly stabilized IMEX methods. The rightmost column is the guide to setting $\bar{p}$.}
	\label{table:amp poly 2}
	\centering
	\begin{tabular}{lll}
		\hline\noalign{\smallskip}
		Method & Amplification Polynomial & $\bar{p} \in$  
		\\
		\noalign{\smallskip}\hline\noalign{\smallskip}
		CNAB 
		& $\left( 1 - \frac{1}{2}z\bar{p} \right)\xi^2
		- \left( 1 + z\left(\frac{3}{2}-\bar{p}  \right) \right) \xi + \frac{1}{2}z( 1-\bar{p} )$
		& $[1,\infty)$
		\\
		CNLF 
		& $\left( 1-\bar{p}z \right) \xi^2 -2z( 1-\bar{p} )\xi -( 1+\bar{p}z )$
		& $[1/2,\infty)$
		\\
		SBDF2 
		& $\left( \frac{3}{2} - z\bar{p} \right) \xi^2
		- 2\left( 1 + z( 1-\bar{p} )\right) \xi 
		+ \frac{1}{2} + z( 1-\bar{p} )
		$
		& $[3/4,\infty)$ 
		\\
		\noalign{\smallskip}\hline
	\end{tabular}
\end{table}

Applying the same analysis to SBDF3 and SBDF4, we observe  two crucial differences.
First, we find that for the real $\lambda$ case that 
the parameter $\bar{p}$ must be restricted to a finite interval to achieve unconditional stability. 
See  Table~\ref{table:amp poly 34} for the corresponding results.
Second, in contrast to SBDF1, EIN, and the second order IMEX methods, the stability region for complex $\lambda$ does not contain the entire left half-plane regardless of the choice of $\bar{p}$. 
As a consequence, we limit ourselves to the real-valued $\lambda$ case.

The significance of having a finite interval is addressed as part of Sect.~\ref{sect: 3 key properties} where it is demonstrated that the finite interval property 
renders the linearly stabilized SBDF3 and SBDF4 methods\footnote{For simplicity, going forward we will refer to linearly stabilized IMEX methods without prefacing by ``linearly stabilized''.
For example, we will refer to the ``the linearly stabilized CNAB method'' as CNAB and the ``the linearly stabilized SBDF3 method'' as SBDF3, etc.} ineffective in many situations.

\begin{table}
	\caption{Amplification polynomials for the linearly stabilized SBDF3 and SBDF4 schemes.   The rightmost column is the guide to setting $\bar{p}$.}
	\label{table:amp poly 34}
	\centering
	\begin{tabular}{lll}
		\hline\noalign{\smallskip}
		Method & Amplification Polynomial & $\bar{p} \in$  
		\\
		\noalign{\smallskip}\hline\noalign{\smallskip}
		SBDF3
		& $\left(\frac{11}{6} - z\bar{p} \right)\xi^3
		- 3\left( 1 + z(1-\bar{p} ) \right) \xi^2 
		+ \frac{3}{2}\left( 1  + 2z( 1-\bar{p} ) \right) \xi 
		$
		& $[7/8, 2]$
		\\
		& \phantom{$\left(\frac{11}{6} - z\bar{p} \right)\xi^3
			- 3\left(1 + z( 1-\bar{p} )  \right) \xi^2$}$- \frac{1}{3}\left(1 +  3z( 1-\bar{p} ) \right)$
		\\
		SBDF4
		& $\left(\frac{25}{12} - z\bar{p} \right) \xi^4
		- 4\left(1 + z( 1-\bar{p} ) \right)\xi^3
		+ 3\left(1 + 2z( 1-\bar{p} ) \right)\xi^2 
		$
		& $[15/16, 5/4]$ 
		\\
		& \phantom{$\left(\frac{25}{12} - z\bar{p} \right) \xi^4$}$- \frac{4}{3}\left(1 + 3z( 1-\bar{p} ) \right) \xi
		+ \frac{1}{4}\left(1 + 4z( 1-\bar{p} )\right)$
		\\
		\noalign{\smallskip}\hline
	\end{tabular}
\end{table}

%
%

\section{Analysis of the Methods: 3 Key Properties}
\label{sect: 3 key properties}
This section presents numerical experiments designed to explore the relative
performance of the proposed methods.   
As a result, three criteria are proposed for selecting effective linearly stabilized schemes (in addition to the usual
requirement of unconditional stability). Two problems will motivate our criteria.

\subsection{Test problem 1 and the unbounded \texorpdfstring{$p$}{p}-parameter restriction}
\label{subsect: test problem 1}
Convergence of the proposed schemes will be tested on the 1D axisymmetric mean curvature motion problem:
\begin{subequations}
	\begin{align}
		u_t = \frac{u_{xx}}{1 + u_x^2} - \frac{1}{u} - pu_{xx} + pu_{xx},
		\quad 0< x< 10,\quad t>0,
	\end{align}
	with initial and boundary conditions
	\begin{gather}
		u(x,0) = 1 + 0.10\sin\left(\frac{\pi}{5}x \right), 
		\\
		u(0,t) = u(10,t) = 1.
	\end{gather}
	\label{ammc pmpm}
\end{subequations}

We solve this problem to time $T=0.35$ using $N=2048$ spatial grid nodes. 
A reference solution is generated using Heun's third order Runge-Kutta method \cite{hundsdorfer2013numerical} with a time step-size $\Delta t = \num{1.46e-5}$.
Starting values for multistep schemes are found using the same third order Runge-Kutta method. 
By comparing with the reference solution, an approximation of the max norm relative error is obtained for various time step-sizes $\Delta t$. 
The values of $p$ used for the second order IMEX schemes are chosen with reference to Table~\ref{table:amp poly 2}.   As we shall see, it is impossible to choose valid $p$-values for third and fourth order SBDF.
Numerical experiments at various values of $p$ were conducted and logged for third and fourth order SBDF to illustrate the issue.

\begin{figure}
	\centering
	\includegraphics[width=84mm]{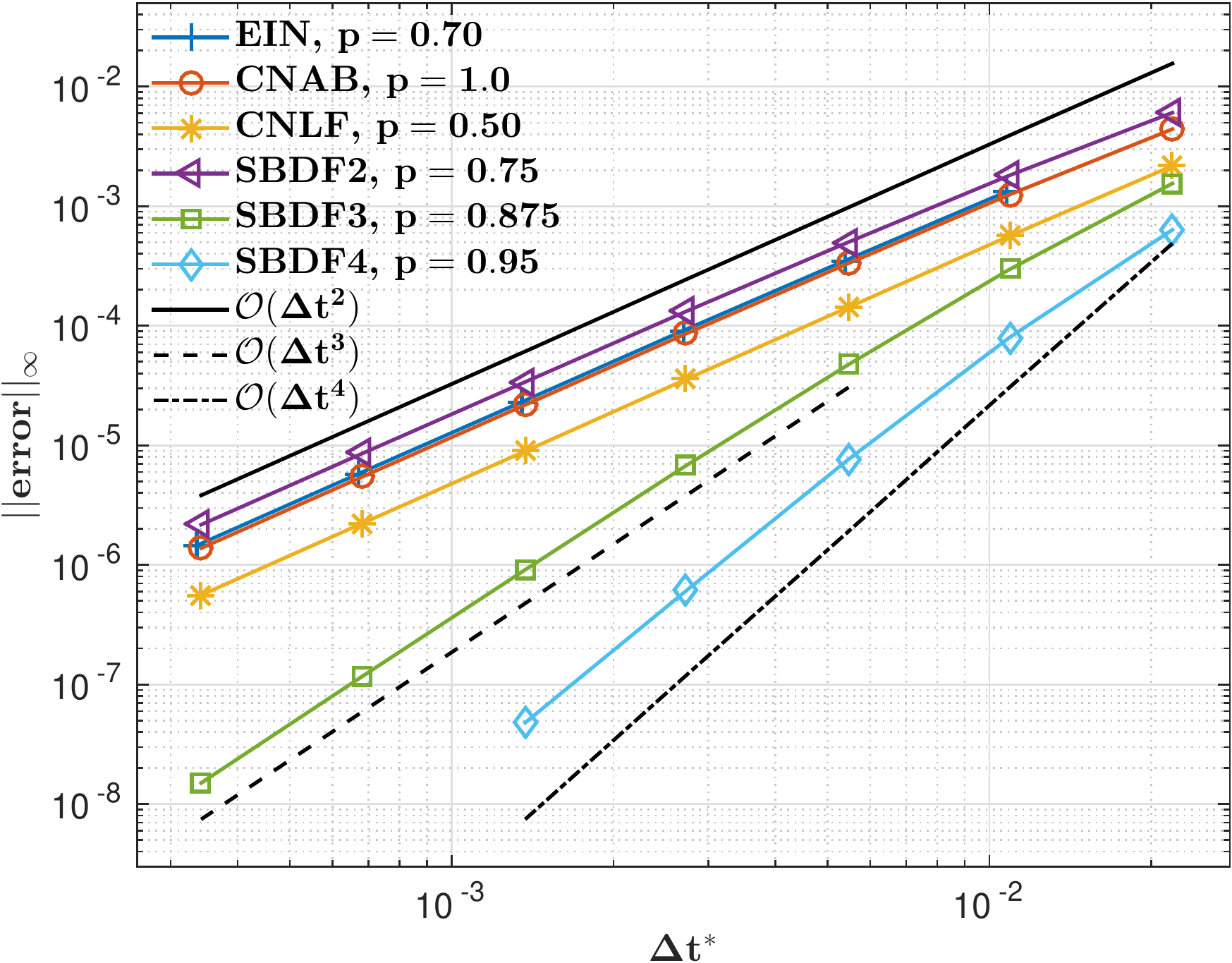}
	\caption{Numerical convergence study to \eqref{ammc pmpm} with IMEX methods. 
Convergence of EIN is also included for comparison; a scaled time step-size is presented to reflect its higher per step cost.}
	\label{fig: lmm conv test}
\end{figure}

Results of a numerical convergence study are shown in Fig.\ \ref{fig: lmm conv test}. Each of the second order methods converge with the expected order of accuracy, with SBDF2 having the largest errors, followed by EIN, CNAB, and CNLF. We note that in the case of EIN, the step-size presented has been scaled down by a factor of 2.5 to account for its nearly 3 times greater cost per step relative to the IMEX methods. 

For the third and fourth order IMEX variants, it appears that SBDF3 converges nicely with $p=0.875$. 
However, SBDF4 does not exhibit fourth order convergence and in fact fails for both $\Delta t =\num{6.84e-04}$ and $\Delta t=\num{3.42e-04}$. We discuss next the cause of SBDF4's failure, and show also that SBDF3 suffers from the same defect.

\subsubsection{Rejecting methods with bounded \texorpdfstring{$p$}{p}-parameter restrictions}
\label{sssect: reject methods with bounded p}
To begin, we tabulate the experimentally observed convergence rates for SBDF3 for various values of $p$. Table \ref{table: sbdf3 nonconvergence} documents three cases. 
The first case ($p=0.875$) is the one already considered in Fig.\ \ref{fig: lmm conv test}. 
The second case ($p=1.475$) exhibits a drastic drop in the observed convergence rate. In the third case ($p=1.675$), 
the method diverges as the time step-size is reduced. We attribute the divergence of SBDF3 and SBDF4 to the fact that their parameter restrictions correspond to bounded intervals,
 in contrast to the unbounded intervals  that we have seen for the second order IMEX schemes.

\begin{table}
	\caption{Observed convergence rates of SBDF3 for various $p$.}
	\label{table: sbdf3 nonconvergence}
	\centering
	\begin{tabular}{lccc}
		\hline\noalign{\smallskip}
		& \multicolumn{3}{c}{Observed convergence rate}
		\\
		$\Delta t$ & $p=0.875$ & $p=1.475$ & $p=1.675$ 
		\\
		\noalign{\smallskip}\hline\noalign{\smallskip}
		\num{2.19e-2} & -- & -- & -- 
		\\ 
		\num{1.09e-2} & \num{2.39} & \num{2.39} & \num{2.39}
		\\
		\num{5.47e-3} & \num{2.64} & \num{2.64} & \num{1.23}
		\\
		\num{2.73e-3} & \num{2.80} & \num{2.80} & \num{-1.51}
		\\ 
		\num{1.37e-3} & \num{2.90} & \num{2.90} & \num{-3.95}
		\\
		\num{6.84e-4} & \num{2.95} & \num{2.95} & diverge
		\\ 
		\num{3.42e-4} & \num{2.97} & \num{2.08} & diverge
		\\ 
		\noalign{\smallskip}\hline
	\end{tabular}
\end{table}

\begin{figure}[htb!]
	\centering
	\begin{minipage}[t]{0.50\textwidth}
		\includegraphics[width=0.96\textwidth]{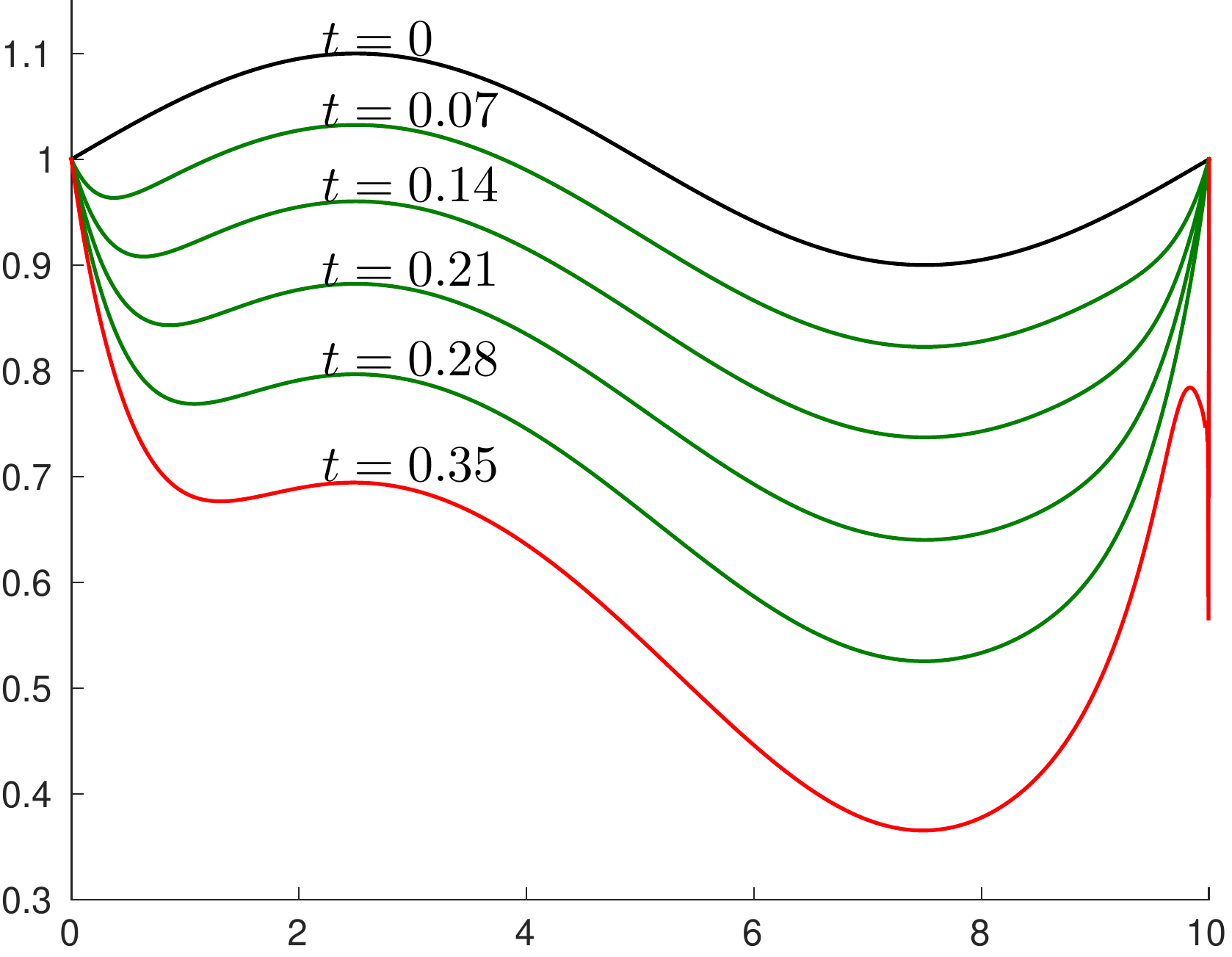}
	\end{minipage}%
	\begin{minipage}[t]{0.50\textwidth}
		\includegraphics[width=0.96\textwidth]{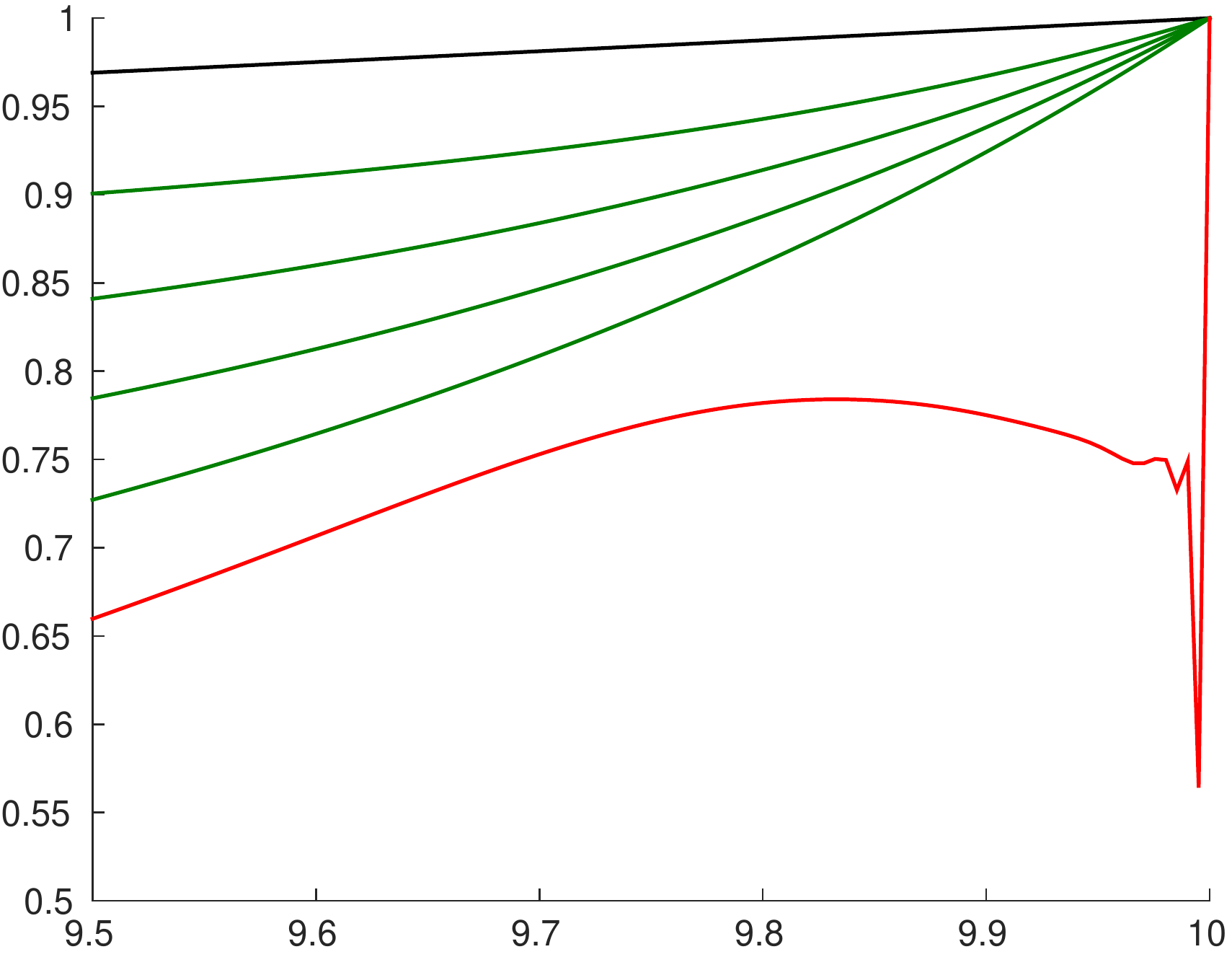}
	\end{minipage}
	\caption{Left: SBDF3 approximations of the solution using  $p=1.675$ and $\Delta t=\num{9.2e-4}$. We observe an instability develop near the right-hand boundary of the bottommost curve. 
Right: A zoom-in to the right-hand boundary.}
	\label{fig:sbdf3 instab}
\end{figure}

\textcolor{black}{To see this, recall the relation \eqref{ammc eig} and the subsequent calculation \eqref{ammc ein p}. If given a parameter restriction $\bar p \in[ \bar{p}_\text{min},\bar{p}_\text{max}]$, it would be necessary to select a $p$ satisfying
\begin{align}
p \geq \max_{1\leq j\leq N} \bar{p}_\text{min} \frac{1}{1 + (D_1 u^n_j)^2}
\quad\text{and}\quad 
p \leq \min_{1\leq k\leq N} \bar{p}_\text{max} \frac{1}{1 + (D_1 u^n_k)^2}. 
\end{align}
For SBDF3, we would need
\begin{align} 
	\max_{1\leq j\leq N}\frac{7}{8(1 + (D_1 u^n_j)^2)} \leq p 
	\leq \min_{1\leq k\leq N}\frac{2}{1 + (D_1 u^n_k)^2}.
\end{align} }
The method fails as no single value of $p$ is able to satisfy the parameter constraint at every grid node simultaneously. 
From Fig.\ \ref{fig:ammc sol}, we see that $\max_j (D_1 u^n_j)^2$ is increasing as the solution evolves and occurs near the boundaries. 
For long enough times, we expect instabilities to develop, and to develop in those regions first. This analysis is corroborated by Fig.\ \ref{fig:sbdf3 instab}, where we see that instabilities develop near the right-hand boundary.

With SBDF4, the instability is more pronounced because  the restriction is tighter. 
While the result in Fig.~\ref{fig: lmm conv test} appeared acceptable at coarse step-sizes, this was a consequence of using a small  number of time steps, as there were too few steps to allow instabilities to grow to an extent that they dominate the solution. 
We conclude that linear stabilization with SBDF3 or SBDF4 is not recommended. 

A natural follow-up question is whether {\it all} third and fourth order IMEX schemes are unsuited for combination with linear stabilization. To this we provide a partial answer. Third order, three step schemes form a three parameter family, and fourth order, four step schemes form a four parameter family \cite{ascher1995implicit}. 
An extensive search through this parameter space was conducted, but we were unable to find any schemes with an unbounded $p$-parameter restriction.

This leaves us a number of competing second order methods to consider. 
Next, we introduce a second test problem and compare the performance of our IMEX based schemes and the EIN method of \cite{duchemin2014explicit}.

\subsection{Test problem 2: Error constants and amplification factors}
\label{subsect: test problem 2}
Of the schemes that we have proposed, only the second order variants are worth further investigation. 
Including the EIN method, we have a total of four second order linearly stabilized schemes to consider. 
We now proceed with a comparison of these methods.

Let us consider as a test problem the following nonlinear PDE from \cite{vdHouwen1982on}:
\begin{subequations}
	\begin{align}
		u_t = \Delta (u^5),
		\quad 0\leq x,y\leq 1,
		\quad t > 0,
		\label{nl5a}
	\end{align}
	with initial and boundary conditions set so that the exact solution is 
	\begin{align}
		u(x,y,t) = \left(\frac{4}{5}(2t+x+y)\right)^{1/4}.
		\label{nl5b}
	\end{align}
	\label{nl5}%
\end{subequations}
Discretizing with a uniform grid and second order centered differences in space, 
the eigenvalues of the linearization of $\Delta (u^5)$  are estimated to lie in the interval
\begin{align}
	\left[-\frac{64}{h^2}(1+t), -16\pi^2(t+h) \right].
	\label{nl5 eigs}
\end{align}
To solve \eqref{nl5}, we propose stabilization with $p\Delta u$, i.e., replace \eqref{nl5a} with 
\begin{align}
	u_t = \Delta (u^5) - p\Delta u + p\Delta u, \quad 0\leq x,y\leq 1,
	\quad t > 0.
	\label{nl5 p}
\end{align}
The parameter $p$ will then be chosen according to the ratio
\begin{align}
	\frac{p\lambda_L}{\lambda_F}
	\approx \frac{-8p/h^2}{-64(1+t)/h^2}
	= \frac{p}{8(1+t)}.
\end{align}

For equation~\eqref{ammc p}, a von Neumann analysis provided tight eigenvalue estimates. On the other hand, for test problem 2 and others, the estimates may be rough and even grow with $t$.	Nonetheless, $p$ need not be updated; in practice it is initialized and fixed at that initial value for all time steps. Consequently, $p$ may at times be substantially greater than necessary. Moreover, $p$ must also compensate for the fact that the Laplacian is less stiff than the original nonlinear term. Both factors force us to select a relatively large value of $p$.
 
In practice, large $p$-values arise frequently, making test problem 2 particularly interesting for understanding how $p$ affects the discretization errors.

\subsubsection{Loss of accuracy with EIN}
\label{sssect: ein loss of accuracy}
We test our second order methods on \eqref{nl5 p} with initial and boundary conditions set by \eqref{nl5b}. We solve to time $T=0.40$ with a uniform spatial grid size $h =0.015$. 
To generate a reference solution, we use Heun's third order Runge-Kutta method with time step-size $\Delta t=\num{6.25e-6}$. 
(Other standard explicit time-stepping methods require similarly strict step-sizes for stability.) 
Using the linearly stabilized schemes, we solve \eqref{nl5 p} with a variety of time step-sizes up to $\Delta t = \num{1.25e-2}$ and compute the max norm relative error. 

\begin{figure}[htb!]
	\centering
	\includegraphics[width=84mm]{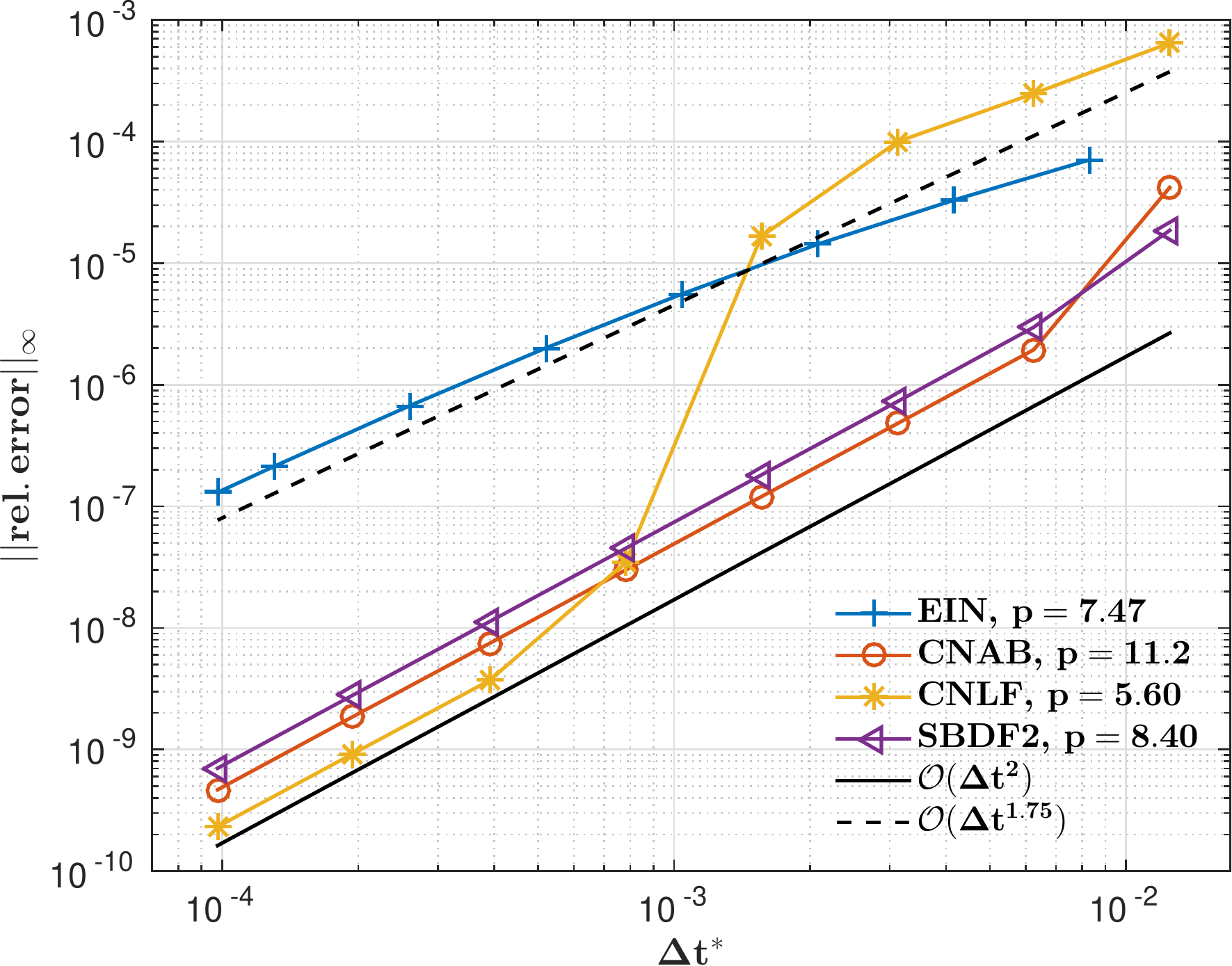}
	\caption{Numerical convergence study for second order methods. 
The time step-size for EIN is scaled to reflect its higher per step cost.}
	\label{fig:nl5 conv a}
\end{figure}

Results of the numerical convergence test are plotted in Fig.~\ref{fig:nl5 conv a}. 
We note that the step-size presented for EIN has been scaled down by a factor of three to reflect 
its approximately 3 times greater cost per step relative to the IMEX methods. 
The results paint an unfavorable picture for the EIN method and for CNLF. 
We first discuss the mediocre performance of the EIN method, after which we comment on the relative performance of the IMEX-based schemes. 

Comparing the performance of EIN in Fig.~\ref{fig: lmm conv test} and Fig.~\ref{fig:nl5 conv a}, 
we observe a significant reduction in the order of accuracy. 
In the former figure, EIN converged with second order accuracy and is competitive with the IMEX-based schemes. 
However, in the latter, we do not (yet) observe second order convergence. 
With further refinement, we find that the EIN method only begins 
to exhibit the full second order rate of convergence for time steps $\Delta t$ below $\num{1e-5}$.

In fact, we argue that this same issue may be observed in the original paper by Duchemin and Eggers \cite{duchemin2014explicit}. 
In their experiments with Hele-Shaw interface flows and with the Kuramoto-Sivashinsky equation, 
the EIN method does not accurately reproduce the reference figures taken from prior publications \cite{hou1994removing,kassam2005fourth}. 
In both cases, a large value of $p$ was needed to obtain unconditional stability. 

We offer an explanation. For each method, consider the local error when applied to \eqref{generic},
\begin{align}
	u(t^{n+1}) - u^{n+1}_*
= C_{k+1}(p) \Delta t^{k+1} + C_{k+2}(p) \Delta t^{k+2} + \ldots ,
\label{local error}
\end{align}
where  $u^{n+1}_*$ is the numerical approximation that is obtained if the past values $u^n$ and $u^{n-1}$ are taken equal to $u(t^n)$ and $u(t^{n-1})$. Note that for first order schemes, $k = 1$, and for second order schemes, $k = 2$. Turning to the coefficients $C_{k+1}(p), C_{k+2}(p)$, etc., we observe that for SBDF1 and the second order IMEX methods, $C_{k+1}(p)$ is linear in $p$, $C_{k+2}(p)$ is quadratic, etc., whereas for EIN, $C_{k+1}(p)$ is quadratic, $C_{k+2}(p)$ is cubic, etc. 
Thus if $p$ is large, EIN requires $\Delta t$ to be set much smaller than the value for SBDF1, etc., before its non-leading order error terms are insignificant. In other words, the observed convergence of the EIN method may suffer in a way similar to Fig.\ \ref{fig:nl5 conv a} whenever large
$p$-values arise.    We provide the leading order, local error constants $C_{k+1}(p)$ in Table~\ref{tab:loc_error_const}.

\begin{table}
	\caption{Local error constants, $C_{k+1}(p)$}
	\label{tab:loc_error_const}
	\centering
	\begin{tabular}{lc} 
		\hline\noalign{\smallskip}
		Method & $C_{k+1}(p) $
		\\
		\noalign{\smallskip}\hline\noalign{\smallskip}
		SBDF1 & $\frac{1}{2}u''(t^n) - pLu'(t^n)$
		\\
			\noalign{\smallskip}\noalign{\smallskip}
		CNAB & $\frac{5}{12} u'''(t^n) - \frac{1}{8} pLu''(t^n)$ 
		\\ 
			\noalign{\smallskip}\noalign{\smallskip}
		CNLF & $\frac{1}{3} u'''(t^n) - pL u''(t^n)$ 
		\\ 
			\noalign{\smallskip}\noalign{\smallskip}
		SBDF2 & $\frac{4}{9}u'''(t^n) - \frac{2}{3}pL u''(t^n)$ 
		\\  
			\noalign{\smallskip}\noalign{\smallskip}
		EIN	& $\frac{1}{2} p^2 L^2 u'(t^n) - \frac{1}{8} p(Lu''(t^n) + 2F'(u)Lu'(t^n)) + \frac{1}{24}u'''(t^n) + \frac{1}{8}F'(u) u''(t^n)$
		\\
	\noalign{\smallskip}\hline\noalign{\smallskip}
	\end{tabular}
\end{table}

\subsubsection{Amplification factors at infinity}
\label{sssect: amp at inf}
In the previous section, we uncovered a deficiency of the EIN method: a large value of $p$ may significantly degrade the observed order of accuracy. Thus an effective linearly stabilized time stepping scheme should have a leading order error term that is linear with respect to $p$. However, this does not explain the miserable performance of CNLF or the sharp dip in the observed convergence of CNAB  near $\Delta t = \num{1e-2}$ (see Fig.\ \ref{fig:nl5 conv a}). To posit an explanation, we think back to our discussion on stability and amplification factors. 
Although we have found stable schemes, we have not yet considered whether the accumulation of slow decaying high frequency error modes 
can drive up the error and force us to use smaller time steps to adequately damp and get the expected convergence order.

To explore this aspect, we consider each method's amplification factor as $z=\lambda\Delta t \to -\infty$. For example, with the EIN method, we have the amplification factor \eqref{ein amp fac}. As $z\to -\infty$, 
\begin{align}
	\lim_{z\to-\infty} \abs{\xi_\mathrm{EIN}} = \abs{\frac{ \bar{p}^2-3\bar{p}+2 }{ \bar{p}^2 }}.
\end{align}

For the multistep schemes, we first find the limiting expression of the amplification polynomial, and then take the max of the magnitude of the two roots. 
Consider CNLF. Starting from $\Phi_\mathrm{CNLF}(z)$ (in Table \ref{table:amp poly 2}) and recognizing $\bar{p} \geq 1/2$, we find 
\begin{align}
\lim_{z\to-\infty} \abs{ \xi_{\text{CNLF}} } = 
	\max\left\{ \frac{\abs{ \bar{p}-1+\sqrt{-2\bar{p}+1} }}{ \bar{p} }, \frac{\abs{1-\bar{p} + \sqrt{ -2\bar{p} +1 }}}{ \bar{p} }\right\} = 1. \label{CNLF_HF}
\end{align}
It is known that the standard CNLF scheme is weakly damping at high frequencies and should not be used for strongly diffusive problems \cite{ascher1995implicit}. 
Equation \eqref{CNLF_HF} shows that this is equally true in the linear stabilization framework: Unless very small time steps are taken with this scheme, the method gives very poor damping of high frequency modes.

\begin{figure}[htb!]
	\centering 
	\includegraphics[width=84mm]{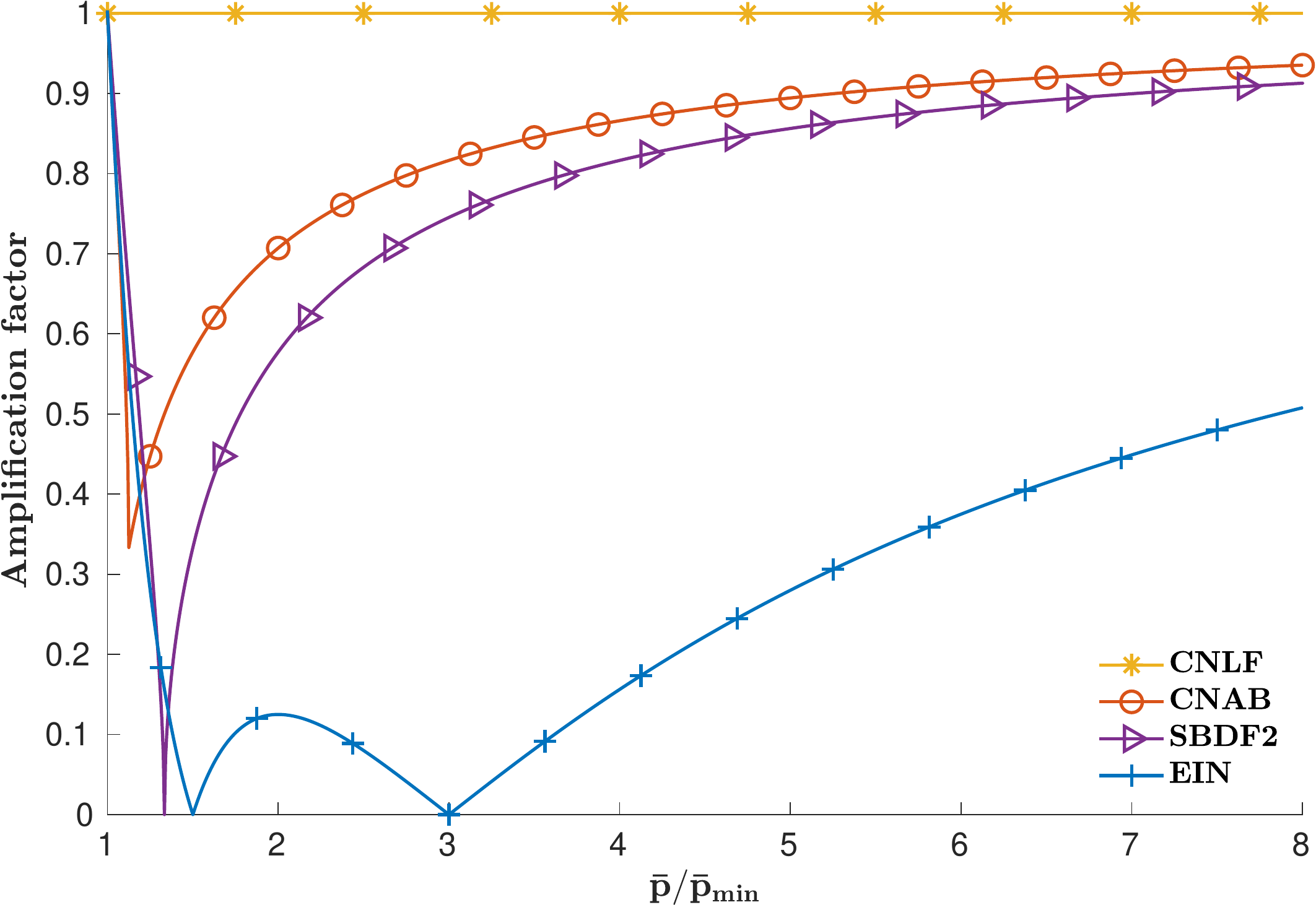}
	\caption[Amplification factors as $z\to-\infty$]{Amplification factors as $z\to-\infty$. The normalization along the horizontal axis is with respect to the lower limit of the parameter restriction of each scheme.}
	\label{fig:damp fac at inf}
\end{figure}

Plots of the amplification factors as $z\to-\infty$ are provided in Fig.~\ref{fig:damp fac at inf}  for all our second order schemes. 
Recalling that $\bar{p}$ must belong to an interval with some lower bound $\bar{p}_{\text{min}}$, we plotted along a normalized parameter
range starting at the abscissa value $\bar{p}/\bar{p}_{\text{min}} = 1$. 
Out of the second order methods, we find that the EIN method gives a strong decay of high frequency error modes
over a wide range of $\bar p$-values (we are interested in a wide range since $\bar{p}$ is in practice frequently larger than its theoretical minimum).  Out of the second order IMEX schemes, SBDF2 provides the strongest damping.  
The decay of high frequency error is slower for CNAB when large time steps are chosen.   
Nonetheless, this scheme remains a useful alternative due to its relatively small error constant.

\section{Higher Order with Exponential Integrators}
\label{sect: exp rk}
The investigation with IMEX methods left us with a major question: Since the linearly stabilized schemes based on SBDF3 and SBDF4 were shown to be unsuitable for practical use, is it possible to construct practical high order linearly stabilized time stepping methods? 
In this section, we consider a second and a fourth order exponential Runge-Kutta method from Cox and Matthews \cite{cox2002exponential} and assess whether they are suited to linear stabilization with respect to the three key properties discussed in the previous section.

\subsection{Exponential Runge-Kutta}
\label{subsect: exp rk}
As we did in Section~\ref{subsect: imex formulas}, consider the ODE
	\begin{align}
	u' = f+ g,
	\label{eq:ODE}
	\end{align}
where $g$ represents a stiff linear term and $f$ represents the remaining nonlinear/nonstiff terms.
Exponential time differencing methods (or exponential integrators)  treat the linear part of (\ref{eq:ODE}) exactly, and approximate the nonlinear part by some suitable quadrature formula.

Our investigation covers explicit exponential Runge-Kutta methods only. This family of one-step methods has the form 
\begin{subequations}
	\begin{align}
		u^{n+1} &= e^{\Delta t g} u^n
		+ \Delta t \sum^s_{i=1} b_i(\Delta t g)f(U^{n,i}), 
		\\
		U^{n,i} &= e^{c_i\Delta t g} u^n 
		+ \Delta t\sum^{i-1}_{j=1} a_{ij}(\Delta t g) f(U^{n,j}),
	\end{align}
	\label{exp rk general}
\end{subequations}
and can be presented in the familiar Butcher tableau:
\newcommand\raisepunct[1]{\,\mathpunct{\raisebox{-4.20ex}{#1}}}
\begin{align}
	\begin{tabular}{c|cccc}
		$c_1$
		&  
		\\
		$c_2$ & $a_{21}$ & 
		\\
		$\vdots$ & $\vdots$ & $\ddots$ 
		\\
		$c_s$ & $a_{s1}$ & $\cdots$ & $a_{s,s-1}$ 
		\\ \hline 
		& $b_1$ & $\cdots$ & $b_{s-1}$ & $b_s$
	\end{tabular}\raisepunct{.}
\end{align}
In particular, we focus on the second and fourth order exponential Runge-Kutta formulas of Cox and Matthews \cite{cox2002exponential}:
\renewcommand\raisepunct[1]{\,\mathpunct{\raisebox{-0.9ex}{#1}}}
\begin{align} 
	\begin{tabular}{c|cc}
		$0$ 
		& 
		\\
		\num{1} 
		& $\varphi_{1,2}$ 
		&
		\\ \hline
		& $\varphi_1 - \varphi_2$ 
		& $\varphi_2$
	\end{tabular}\raisepunct{,}
	\label{etdrk2 butcher}
\end{align}
\renewcommand\raisepunct[1]{\,\mathpunct{\raisebox{-3.80ex}{#1}}}
\begin{align} 
	\begin{tabular}{c|cccc}
		$0$ 
		&  
		\\
		\num{1/2} 
		& $\frac{1}{2}\varphi_{1,2}$ 
		&
		\\ 
		\num{1/2} 
		& 0 
		& $\frac{1}{2}\varphi_{1,3}$ 
		& 
		\\
		\num{1} 
		& $\frac{1}{2}\varphi_{1,3}(\varphi_{0,3}-1)$ 
		& 0 
		& $\varphi_{1,3}$ 
		& 
		\\ \hline 
		& $\varphi_1 - 3\varphi_2 + 4\varphi_3$ 
		& $2\varphi_2 - 4\varphi_3$ 
		& $2\varphi_2 - 4\varphi_3$ 
		& $4\varphi_3 - \varphi_2$ 
	\end{tabular}\raisepunct{,}
	\label{etdrk4 butcher}
\end{align}
where 
\begin{align}
	\varphi_{k+1}(z) = \frac{\varphi_k(z) - 1/k!}{z}, 
	\quad \varphi_0(z) = \exp(z), 
	\qand 
	\varphi_{i,j}(z) = \varphi_{i}(c_j z).
	\label{phi functions}
\end{align}
We refer to this pair of exponential Runge-Kutta methods as ETDRK2 and ETDRK4, respectively.

\subsection{Linearly stabilized ETDRK2 and ETDRK4}
\label{subsect: etdrk2 and etdrk4}
In Section~\ref{sect: 3 key properties}, we identified criteria for assessing the practicality of linearly stabilized methods. We provide in this section a first assessment of ETDRK2 and ETDRK4 relative to the criteria. Then, numerical experiments in Section~\ref{sect: numerical experiments} will provide further insight.

First, we apply the schemes \eqref{etdrk2 butcher} and \eqref{etdrk4 butcher} to the modified test equation \eqref{mte} and impose unconditional stability. We are only interested in schemes with an unbounded parameter restriction. For ETDRK2 and ETDRK4, with the help of the computer algebra system, \textsc{Maple}\textsuperscript{\texttrademark}, we determined the parameter restriction to be $\bar{p} \geq 1/2$ for the real $\lambda$ case, $\lambda_F, \lambda_L < 0$.
For complex $\lambda$, we were unable to find a parameter restriction with guaranteed stability.   However,  numerical evidence suggests that the stability region of ETDRK2 and ETDRK4 
will contain the left half-plane for $\bar p \geq 3/4$ and $\bar p \geq 1$, respectively.
%

\begin{figure}[htb!]
	\centering
	\includegraphics[width=84mm]{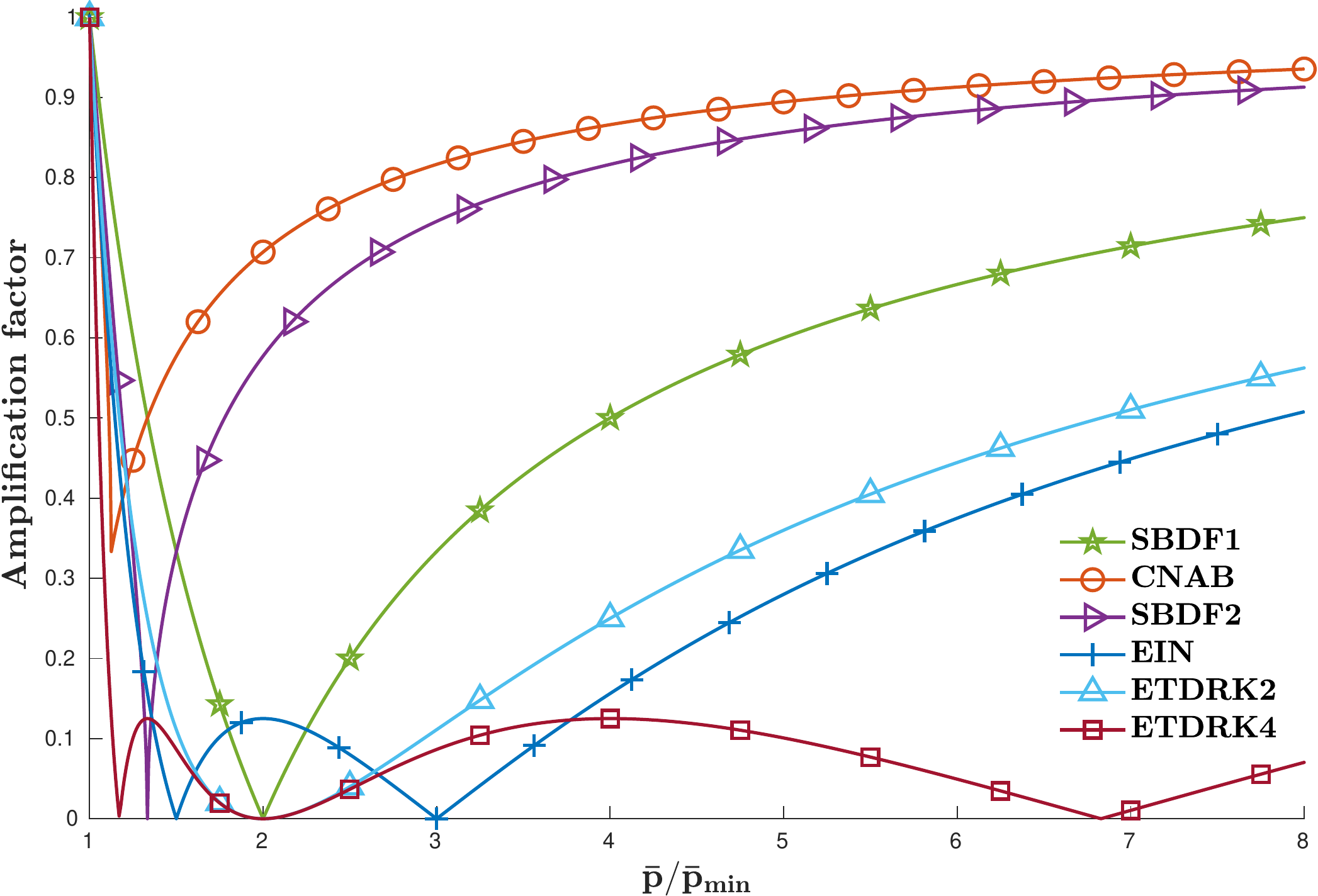}
	\caption{Amplification factors as $z\to-\infty$ with ETDRK schemes. Included also for the purpose of comparison are some first and second order schemes. The normalization along the horizontal axis is with respect to the lower limit of the parameter restriction of each scheme}
	\label{fig:damp fac at inf etd}
\end{figure}

We next consider our remaining two criteria, specifically,  the amplification factor as $z\to-\infty$, and the local error of the numerical scheme. The former is plotted in Fig.\ \ref{fig:damp fac at inf etd} for both ETDRK2 and ETDRK4. It shows that the ETDRK schemes provide strong damping as $z\to-\infty$ for a wide range of $p$.
The latter property, the local error  expanded as \eqref{local error}, is less favorable to the two methods. Recall that for EIN the observed convergence was poor for large $p$ because the coefficients $C_{k+1}(p), C_{k+2}(p
),$ etc., are degree $2,3$, etc., polynomials in $p$, respectively. Likewise for ETDRK2, we find $C_{k+1}(p), C_{k+2}(p), \ldots$ are quadratic, cubic, etc., and for ETDRK4, $C_{k+1}(p), C_{k+2}(p), \ldots$  are quartic, quintic, etc.    
Consequently, these schemes may fare poorly when $p$ is large.  
We emphasize, however, that exponential time differencing schemes should not be entirely ruled out.    Indeed, experiments presented in Section~\ref{subsect: mcm} show that ETDRK2 and ETDRK4 can outperform SBDF2 and CNAB in applications involving small or moderate $p$.

We end this section with a note on the implementation of ETDRK schemes. Any implementation of ETDRK schemes must contend with the stable evaluation of the matrix exponential. In our examples using ETDRK2 and ETDRK4, we follow the direction of Kassam and Trefethen \cite{kassam2005fourth} where a contour integral approach coupled with the trapezoidal rule is used to evaluate functions in the form of \eqref{phi functions}, and, for simplicity, we have only problems on a periodic grid.
Other methods based on scaling and square, Pad\'{e} approximants \cite{moler2003nineteen,higham2008functions}, and Krylov subspace methods \cite{hochbruck1997krylov,sidje1998expokit,simoncini2007recent} could be considered but were not used here.



\section{Numerical Experiments}
\label{sect: numerical experiments}
In this section, we solve a number of stiff PDEs with applications to image inpainting and capturing interface motion. 
For both types of problems, we will give the PDE models and then discuss how to stabilize and select the parameters.
Our experiments will show the practicality of linearly stabilized schemes in 2D and 3D. 

Before proceeding further, we would like to make a few notes on our implementation. 
As stated from the outset, our goal is to provide simple, accurate, and efficient time stepping methods for nonlinear PDEs. 
Consistent with these objectives, the choice of $p$ is fixed throughout the time evolution. 
Alternatively, one could adapt $p$ as the solution evolves to avoid overestimates of $p$ that could lead to larger errors. 
However, we do not pursue that here. 
So while our theory speaks of approximating the eigenvalues of the linearized system, we do not incur this cost in our computations. 
We further note that a static value of $p$ offers the advantage that the linear system to be solved is the same at each time step, i.e.\ the matrix to be inverted is static. Any expensive preprocessing/factorizing of this matrix needs only to be done once.

\subsection{Image Inpainting}
\label{subsect: inpainting}
Image inpainting is the task of repairing corrupted images and damaged artwork \cite{bertalmio2000image}. In the inpainting examples to follow, the user identifies the region to be inpainted in the image.  
From there, a PDE model is evolved to fill-in the inpainting region using the neighboring information.

Two PDE models are selected. The first is a second order model from Shen and Chan \cite{shen2002mathematical}, 
\begin{align}
	u_t  = \nabla \cdot \left(\frac{\nabla u}{\sqrt{\abs{\nabla u}^2 + \epsilon^2}} \right) 
	+ \lambda_D(u_0 - u), 
	\label{bv inpaint}
\end{align}
and the second is a fourth order model from Sch{\"o}nlieb and Bertozzi \cite{schonlieb2011unconditionally},
\begin{align}
	u_t  = -\Delta \nabla \cdot \left(\frac{\nabla u}{\sqrt{\abs{\nabla u}^2 + \epsilon^2}} \right) 
	+ \lambda_D(u_0 - u).
	\label{tvhneg inpaint}
\end{align}
We refer to these as TV inpainting and TV-H$^{-1}$ inpainting, respectively.
In both inpainting models, $u$ is the solution and the restored image, $u_0$ is the initial corrupted image, and $\epsilon > 0$ is a regularization parameter. 
Denoting the image domain $\Omega$ and the inpainting region $D$, $\lambda_D$ is then defined as 
\begin{align}
	\lambda_D(x)
	= \begin{cases}
		\lambda_0, & x\in \Omega\setminus D
		\\
		0, &\text{otherwise},
	\end{cases}
\end{align}
for some $\lambda_0 > 0$. In our experiments, we set $\lambda_0 = 100$.

\begin{figure}[htb!]
	\centering
	\begin{minipage}{0.48\linewidth}
		\includegraphics[width=\textwidth]{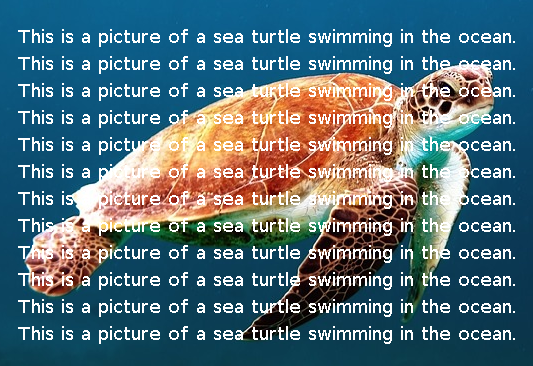}
	\end{minipage} %
	\begin{minipage}{0.48\linewidth}
		\includegraphics[width=\textwidth]{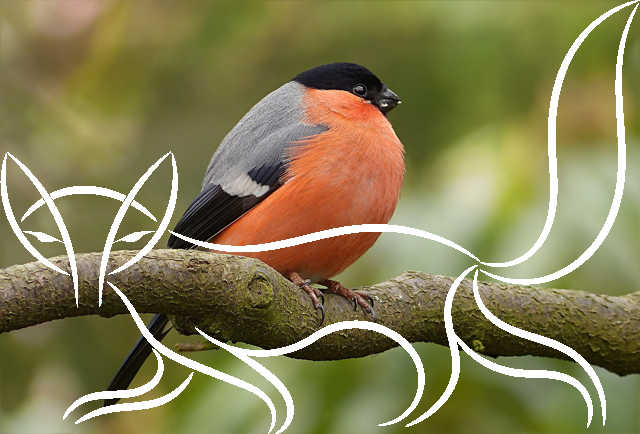}
	\end{minipage} %
	\caption[Vandalized!]{Left: Photograph of a sea turtle overwritten with text. Right: Photograph of a bullfinch vandalized with a cartoon fox.}
	\label{fig:vandalized}
\end{figure}

For initial conditions, we have vandalized two images as shown in Fig.\ \ref{fig:vandalized}. One is a photograph of a sea turtle covered with text that we would like to restore by removing the text. The second is a photo of a bullfinch where the fox-shaped figure requires removal. 
Although the latter may look simpler, it is in fact a more challenging example because the thickness of the inpainting region requires 
an effective extension of level lines over longer distances \cite{schonlieb2011unconditionally}.

The images are restored by evolving each color channel according to the PDE models. 
Spatial discretization is by second order centered differences with uniform spacing, $h$.
As a stopping criterion, we iterate until the relative error of each of the three colour channels (RGB) falls below a  prescribed threshold, $\delta$. That is, 
\begin{align}
	\max\left\{
	\frac{\norm{U^{n+1}_R - U^n_R}}{\norm{U^{n+1}_R}},
	\frac{\norm{U^{n+1}_G - U^n_G}}{\norm{U^{n+1}_G}},
	\frac{\norm{U^{n+1}_B - U^n_B}}{\norm{U^{n+1}_B}}
	\right\} 
	< \delta.
\end{align}
During testing, we found that the choice of $\delta$ depends on the method used and the time step-size. Our choices for $\delta$ and $\Delta t$ reflect only a minimal amount of trial and error testing predicated on attaining satisfactory subjective image quality and (close to) minimal iteration count. 

Next, we report how each method performed on the TV inpainting model.

\subsubsection{TV inpainting}
\label{sssect: tv inpainting}
We first show that the TV inpainting model \eqref{bv inpaint} can easily be approximated by our methods. In the model, there are two terms on the right-hand side, both potentially stiff.  The second term is stabilized by adding and subtracting $-\bar{p}_{\text{min}}\lambda_0 u$, where $\bar{p}_{\text{min}}$ is the minimum value required for unconditional stability when applying the time stepping method to the modified test equation. For the first term, we stabilize by adding and subtracting $p_1\Delta u$. To determined $p_1$, we bound the first term as 
\begin{align}
	\begin{split}
		\nabla \cdot \left(\frac{\nabla u}{\sqrt{\abs{\nabla u}^2 + \epsilon^2}} \right) 
		&= \frac{u_{xx}(u_y^2 + \epsilon^2) + u_{yy}(u_x^2 + \epsilon^2)}{(u_x^2 + u_y^2 + \epsilon^2)^{3/2}} 
		- \frac{2u_xu_y u_{xy}}{(u_x^2 + u_y^2 + \epsilon^2)^{3/2}},
		\\
		&\leq \frac{(u_{xx}+u_{yy})(u_x^2 + u_y^2 + \epsilon^2)}{(u_x^2 + u_y^2 + \epsilon^2)^{3/2}} + \frac{(u_x^2 + u_y^2 + \epsilon^2)u_{xy}}{(u_x^2 + u_y^2 + \epsilon^2)^{3/2}}, 
		\\
		&= \frac{u_{xx} + u_{yy} + u_{xy}}{\sqrt{u_x^2 + u_y^2 + \epsilon^2}}.
	\end{split}
	\label{p1 estimate inpaint 1}
\end{align}
We then consider the auxiliary equation $u_t = u_{xx} + u_{yy} + u_{xy}$ discretized by centered differences in space and forward Euler in time and apply a von Neumann analysis with $u^n_{jk} = \xi^n \exp(i\omega_1jh)\exp(i\omega_2kh)$ to get 
\begin{align}
	\frac{\xi  - 1}{\Delta t} 
	= \frac{1}{h^2}(-4 + 2\cos(\omega_1 h) + 2\cos(\omega_2 h) - \sin(\omega_1 h)\sin(\omega_2h))
	\geq -\frac{8}{h^2}.
	\label{p1 estimate inpaint 2}
\end{align}
Combined with the assumption of the extreme case, $\sqrt{u_x^2 + u_y^2 + \epsilon^2} \geq \epsilon$, we set $p_1$ according to 
\begin{align}
	p_1\frac{8/h^2}{8/(\epsilon h^2)} \geq \bar{p}_{\text{min}} \iff 
	p_1 \geq \frac{1}{\epsilon}\bar{p}_{\text{min}}.
	\label{eq: tv_p}
\end{align}

\begin{figure}[htb!]
	\centering
	\begin{minipage}{0.48\linewidth}
		\includegraphics[width=\textwidth]{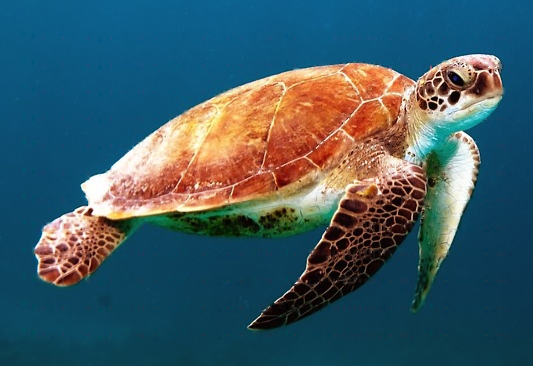}
	\end{minipage} %
	\begin{minipage}{0.48\linewidth}
		\includegraphics[width=\textwidth]{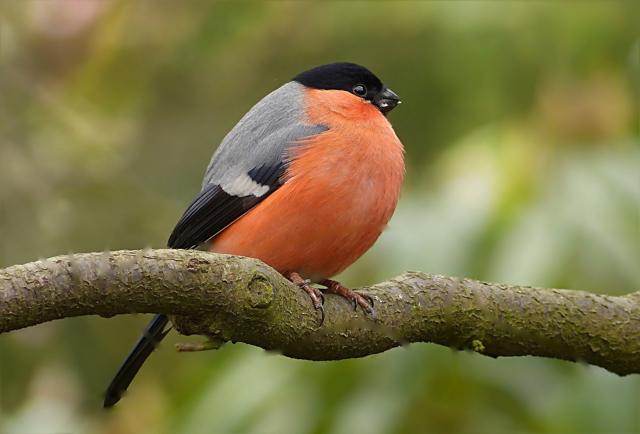}
	\end{minipage} %
	\caption[Image restoration by TV inpainting]{Image restoration by TV inpainting using  CNAB (here we use 41 iterations for the sea turtle and 63 iterations for the bullfinch).}
	\label{fig:bv inpainting}
\end{figure}

\begin{figure}[htb!]
	\centering
	\begin{minipage}{0.48\linewidth}
		\includegraphics[width=\textwidth]{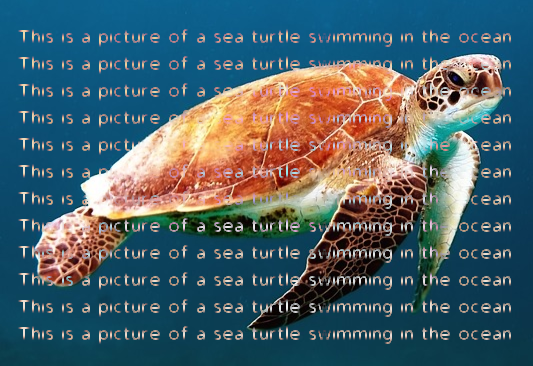}
	\end{minipage}
	\begin{minipage}{0.48\linewidth}
		\includegraphics[width=\textwidth]{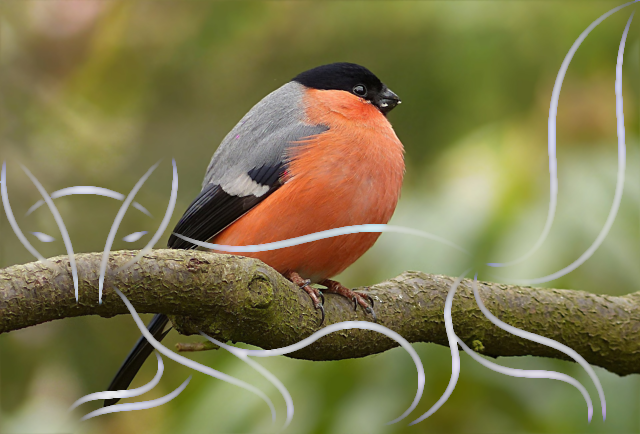}
	\end{minipage}
	\caption[Image restoration by TV inpainting]{Incomplete image restoration by TV inpainting using SBDF1.  The number of iterations is chosen to match the CNAB iteration count (41 iterations for the turtle, and 63 iterations for the bullfinch).   Artifacts are clearly visible.}
	\label{fig:bv inpainting-SBDF1}
\end{figure}

The images in Fig.\ \ref{fig:bv inpainting} are restored by TV inpainting via CNAB. The iteration count, the time step-size $\Delta t$, and the tolerance $\delta$ are listed in Table \ref{tab:bv iter counts} for the restoration of the sea turtle image and for the bullfinch image. The results can be compared to those with SBDF2 and SBDF1. We note that CNAB is slightly faster than SBDF2, and both significantly outperform SBDF1. For the more difficult case of the bullfinch image, SBDF1 needed over four times as many iterations to process to the same image quality as the second order methods. In particular, if we use SBDF1 with the CNAB iteration count  (i.e., 41 iterations for the sea turtle and 63 iterations for the bullfinch) we obtain an incomplete restoration with strong artifacts; see Figure\ \ref{fig:bv inpainting-SBDF1}.   Finally, we note that with all three methods we have set $\epsilon=0.10$ and $\lambda_0=100$.
 

\begin{table}[htb!]
	\caption{Iteration counts for TV image restoration.}
	\label{tab:bv iter counts}
	\begin{minipage}[b]{0.48\linewidth}
		\centering
		\begin{tabular}{cccc}
			\hline\noalign{\smallskip}
			& \multicolumn{3}{c}{Sea Turtle} 
			\\ 
			& Iterations & $\Delta t$ & $\delta (\times 10^{-4})$
			\\ 
			\noalign{\smallskip}\hline\noalign{\smallskip}
			SBDF1 &122 & 0.88 & 20 
			\\
			SBDF2 &48 & 0.10 & 16
			\\
			CNAB &41 & 0.12 & 24
			\\
			\noalign{\smallskip}\hline
		\end{tabular}
	\end{minipage}
	\hspace{0.15cm}
	\begin{minipage}[b]{0.48\linewidth}
		\centering
		\begin{tabular}{cccc}
			\hline\noalign{\smallskip}
			&\multicolumn{3}{c}{Bullfinch} 
			\\ 
			&Iterations & $\Delta t$ & $\delta (\times 10^{-4})$
			\\ 
			\noalign{\smallskip}\hline\noalign{\smallskip}
			SBDF1 & 322 & 0.88 & 20 
			\\
			SBDF2 & 77 & 0.10 & 3.0
			\\
			CNAB & 63 & 0.12 & 3.9
			\\
			\noalign{\smallskip}\hline
		\end{tabular}
	\end{minipage}
\end{table}

\subsubsection{TV-H\texorpdfstring{$^{-1}$}{-1} inpainting}
\label{sssect: tvhneg inpainting}
For TV-H$^{-1}$ inpainting, we stabilize \eqref{tvhneg inpaint} as
\begin{align}
	u_t = -\Delta \nabla \cdot \left(\frac{\nabla u}{\sqrt{\abs{\nabla u}^2 + \epsilon^2}} \right) + \lambda_D(u_0 - u)  + p_1\Delta^2 u + \bar{p}_{\text{min}}\lambda_0 u - p_1\Delta^2 u - \bar{p}_{\text{min}}\lambda_0 u.
	\label{tvhneg stabilized}
\end{align}
As with TV inpainting, we determine a bound for setting $p_1$:
\begin{align}
	p_1\frac{(8/h^2)^2}{(8/h^2)(8/(\epsilon h^2))} \geq \bar{p}_{\text{min}} 
	\iff p_1 \geq \frac{1}{\epsilon}\bar{p}_{\text{min}}.
\end{align}

\begin{figure}[htb!]
	\centering
	\begin{minipage}{0.48\linewidth}
		\includegraphics[width=\textwidth]{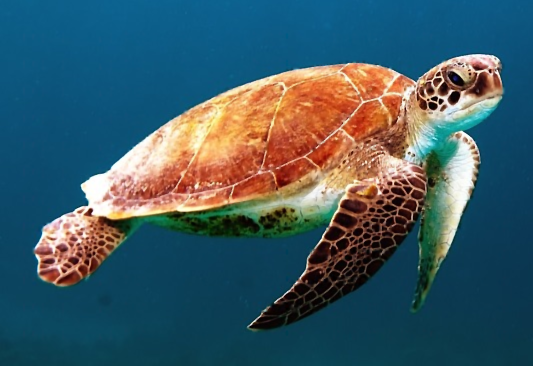}
	\end{minipage} %
	\begin{minipage}{0.48\linewidth}
		\includegraphics[width=\textwidth]{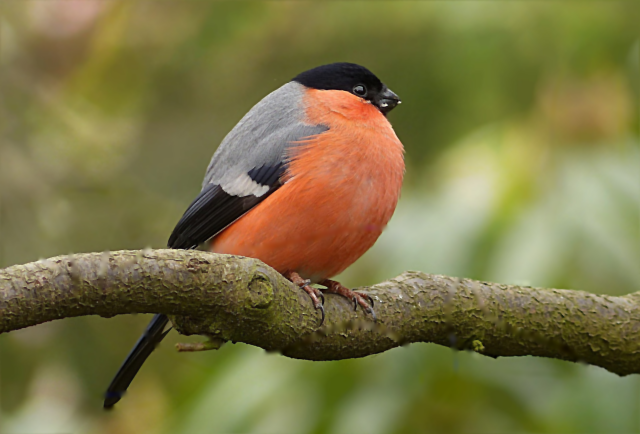}
	\end{minipage} %
	\caption[Image restoration by TV-H$^{-1}$ inpainting]{Image restoration by TV-H$^{-1}$ inpainting using  CNAB (here we use 34 iterations for the sea turtle and 80 iterations for the bullfinch).}
	\label{fig:bvhneg inpainting}
\end{figure}

\begin{figure}[htb!]
	\centering
	\begin{minipage}{0.48\linewidth}
		\includegraphics[width=\textwidth]{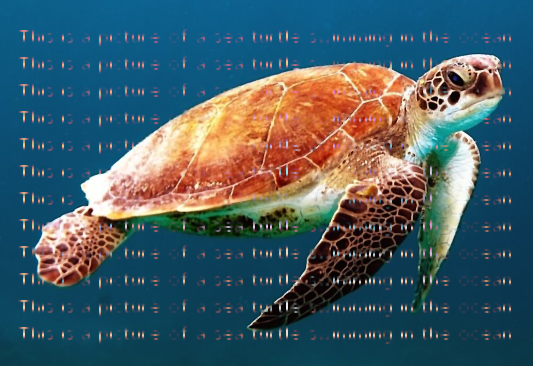}
	\end{minipage}
	\begin{minipage}{0.48\linewidth}
		\includegraphics[width=\textwidth]{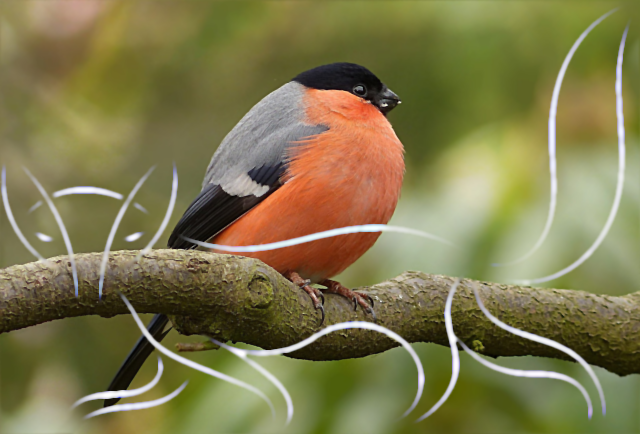}
	\end{minipage}
	\caption[Image restoration by TV-H$^{-1}$ inpainting]{Incomplete image restoration by TV-H$^{-1}$ inpainting using SBDF1.  The number of iterations is chosen to match the CNAB iteration count  (34 iterations for the turtle, and 80 iterations for the bullfinch).   Artifacts are clearly visible.}
	\label{fig:bvhneg inpainting-SBDF1}
\end{figure}

Notably, in the same paper where they propose \eqref{tvhneg inpaint} for image inpainting, the authors offer exactly \eqref{tvhneg stabilized} and time stepping with SBDF1 as the solution algorithm. In Table \ref{tab:bvhneg iter counts}, we list the iteration counts required for each of SBDF1, SBDF2, and CNAB, again with parameters $\epsilon=0.10$ and $\lambda_0=100$.   Once more, the second order methods vastly outperform SBDF1, with SBDF1 needing well over five times the number of iterations as either second order method for the restoration of the bullfinch image.
See Fig.\ \ref{fig:bvhneg inpainting} for the TV-H$^{-1}$  restoration with CNAB.   For comparison purposes, we also give an incomplete SBDF1 restoration using the same number of iterations as CNAB;  see Fig.\ \ref{fig:bvhneg inpainting-SBDF1}.

We should emphasize that more research is needed to decipher the relationship between image size, inpainting region thickness, and how to best choose $\Delta t$ and $\delta$. This will be left to future work. 

Lastly, we mention other relevant developments. In \cite{bredies2010total,papafitsoros2014combined,papafitsoros2013combined}, a number of image restoration models are proposed that involve high order derivatives interacting nonlinearly. Of interest would be to test the effectiveness of our schemes on other inpainting models and run them against the methods that were considered.

\begin{table}[htb!]
	\caption{Iteration counts for TV-H$^{-1}$ image restoration.}
	\label{tab:bvhneg iter counts}
	\begin{minipage}[b]{0.48\linewidth}
		\centering
		\begin{tabular}{cccc}
			\hline\noalign{\smallskip}
			& \multicolumn{3}{c}{Sea Turtle} 
			\\ 
			& Iterations & $\Delta t$ & $\delta (\times 10^{-4})$
			\\ 
			\noalign{\smallskip}\hline\noalign{\smallskip}
			SBDF1 & 122 & 0.88 & 5.3 
			\\
			SBDF2 & 39 & 0.08 & 19
			\\
			CNAB & 34 & 0.08 & 25
			\\
			\noalign{\smallskip}\hline
		\end{tabular}
	\end{minipage}
	\hspace{0.15cm}
	\begin{minipage}[b]{0.48\linewidth}
		\centering
		\begin{tabular}{cccc}
			\hline\noalign{\smallskip}
			&\multicolumn{3}{c}{Bullfinch} 
			\\ 
			&Iterations & $\Delta t$ & $\delta (\times 10^{-4})$
			\\ 
			\noalign{\smallskip}\hline\noalign{\smallskip}
			SBDF1 & 564 & 0.98 & 2.3 
			\\
			SBDF2 & 96 & 0.14 & 2.8
			\\
			CNAB & 80 & 0.16 & 3.1
			\\
			\noalign{\smallskip}\hline
		\end{tabular}
	\end{minipage}
\end{table}


\subsection{Motion by Mean Curvature}
\label{subsect: mcm}
In this section, we study the problem of interface evolution under mean curvature flow. The level set equation for motion by mean curvature is 
\begin{align}
	u_t 
	= \kappa\abs{\nabla u} 
	= \abs{\nabla u} \nabla \cdot \left(\frac{\nabla u}{\abs{\nabla u}} \right).
	\label{mcm}
\end{align}
Our interest is in the time evolution of the interface, $\Gamma=\Gamma(t)$, described by the zero level set of the function $u$,
\begin{equation}
	\Gamma(t) = \{x \in \mathbb{R}^d \mid u(x,t) =0\}.
\end{equation}

We will demonstrate the effectiveness of our schemes on examples similar to those of Smereka \cite{smereka2003semi}. 
In \cite{smereka2003semi}, linearly stabilized SBDF1 was used to take large, stable time steps. In that same paper it was also suggested that Richardson extrapolation may be used to attain second order convergence (although this was not implemented). 
Following \cite{smereka2003semi}, we stabilize \eqref{mcm} with a Laplacian term, $p\Delta u$, to obtain
\begin{align}
	u_t = \kappa \abs{\nabla u} - p\Delta u + p\Delta u.
\end{align}
An analysis similar to \eqref{p1 estimate inpaint 1} and \eqref{p1 estimate inpaint 2} yields $p\geq \bar{p}_{\text{min}}$ to be sufficient for unconditional stability.

Let us point out a key difference between this problem and the inpainting problem of the previous section. In the inpainting problem, the system was to be driven to steady state. As such, we were afforded a range of time step-sizes where the solution method was computationally efficient. Indeed, the step-size did not affect the visual quality. For mean curvature flow, computing time and accuracy are directly related to the choice of step-size. 
Thus we seek large step-sizes subject to maintaining an acceptable level of accuracy.

\subsubsection{Shrinking dumbbell in 2D}
\label{sssect: shrinking dumbbell}
\begin{figure}[htb!]
	\centering
	\begin{minipage}{0.50\textwidth}
		\includegraphics[width=0.96\textwidth]{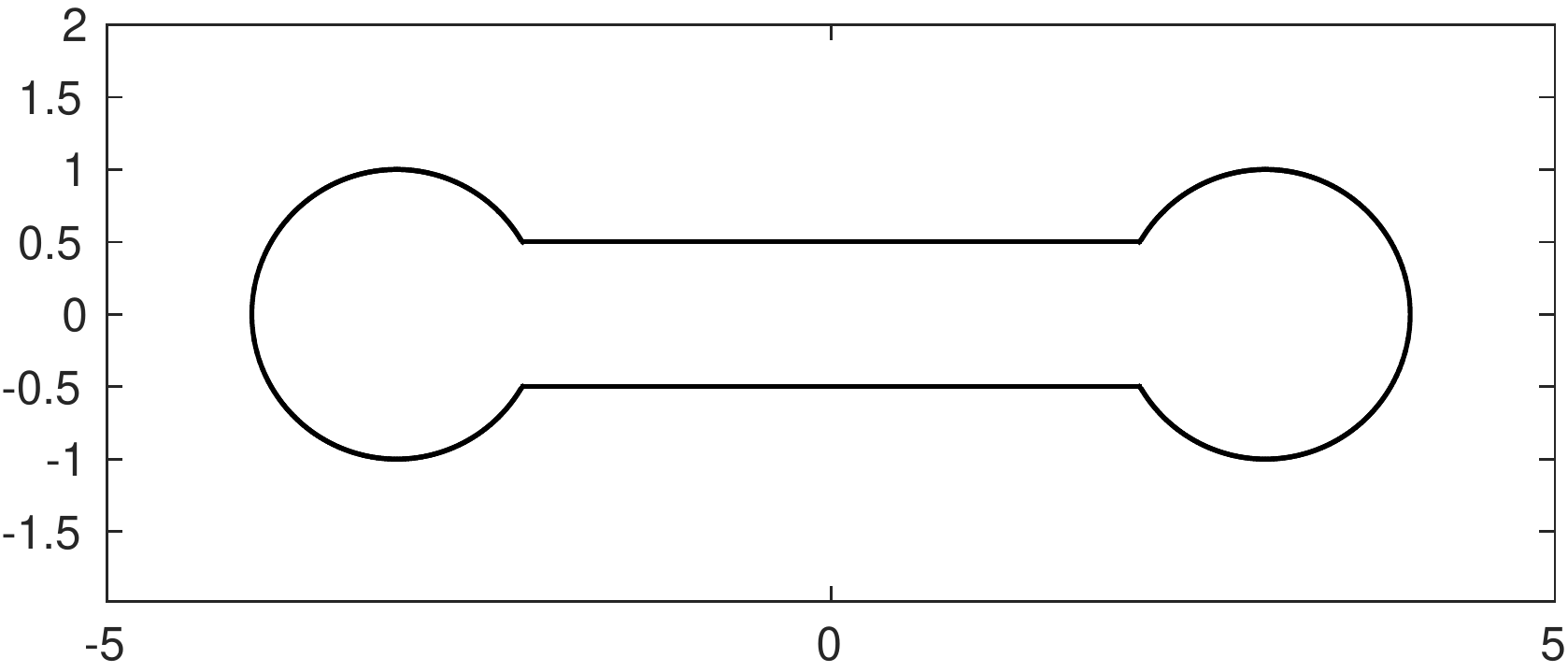}
	\end{minipage}%
	\begin{minipage}{0.50\textwidth}
		\includegraphics[width=0.96\textwidth]{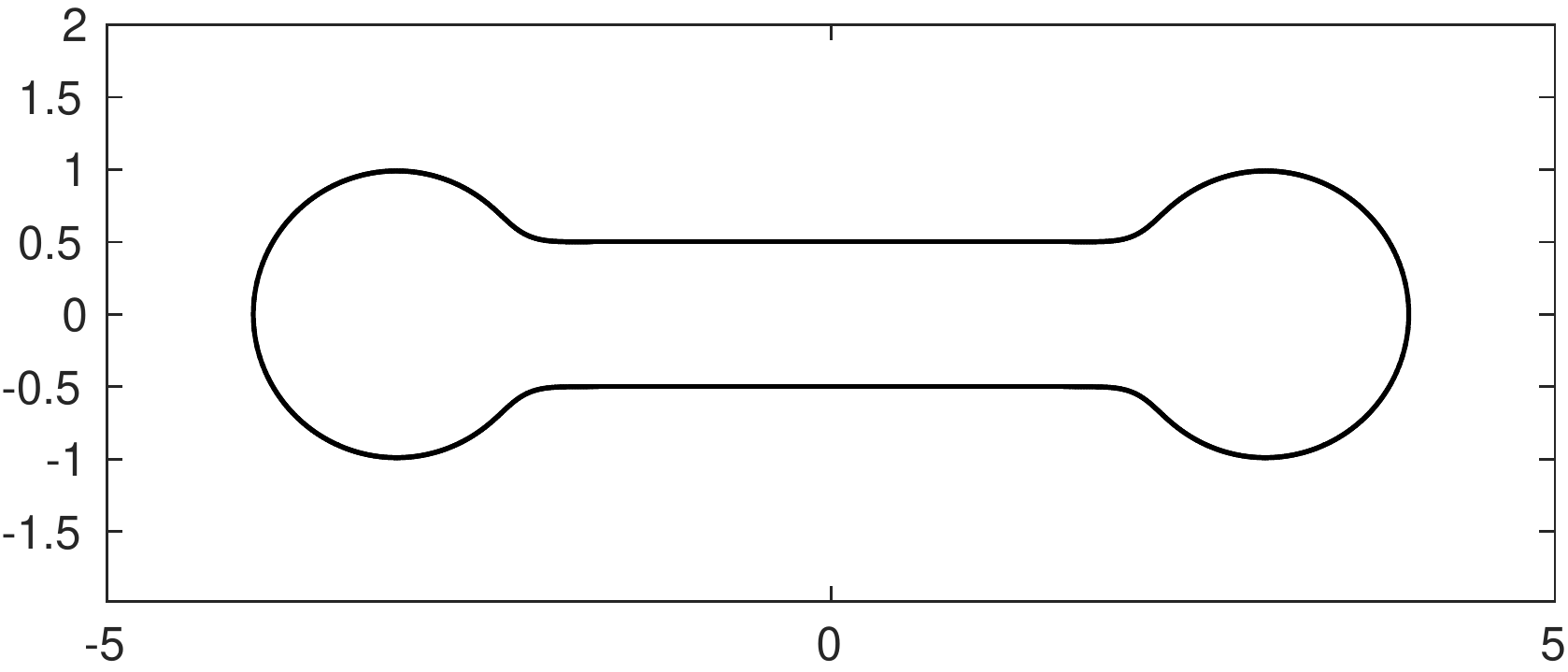}
	\end{minipage}
	\begin{minipage}{0.50\textwidth}
		\includegraphics[width=0.96\textwidth]{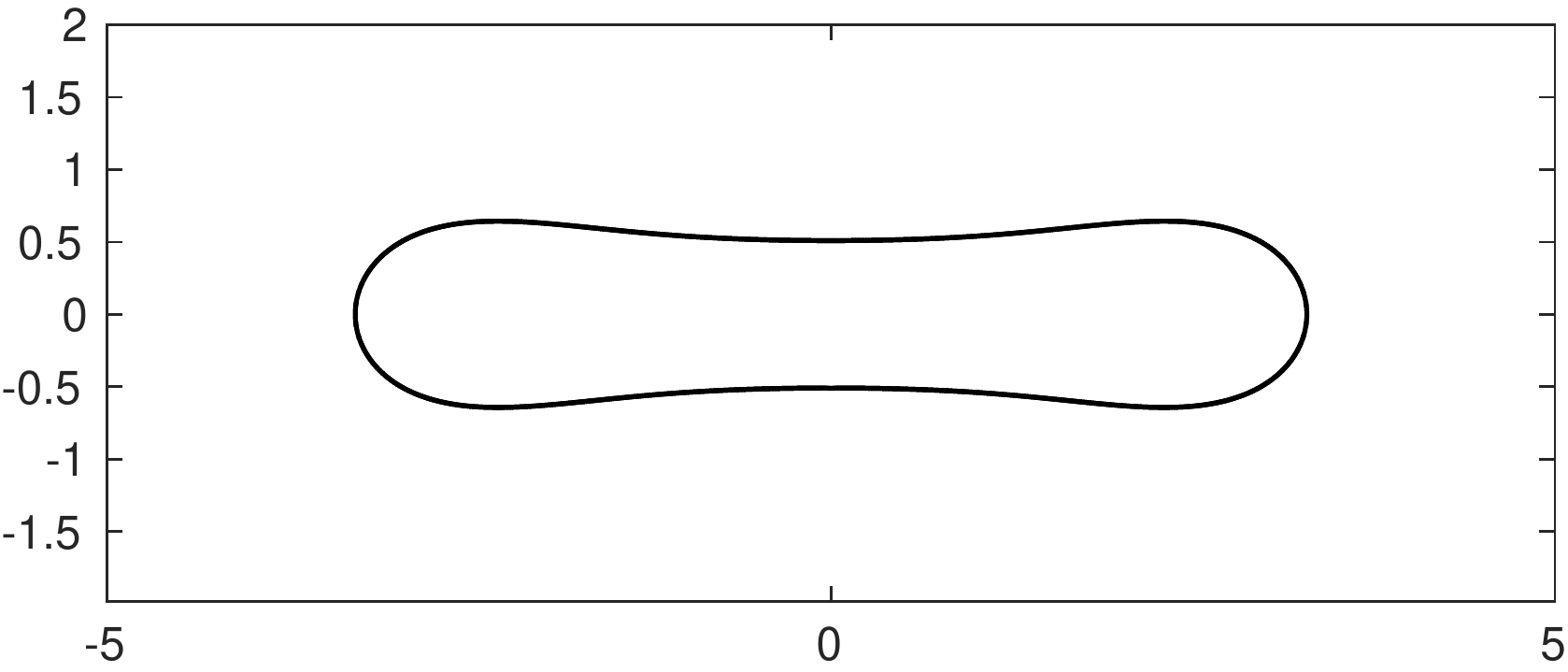}
	\end{minipage}%
	\begin{minipage}{0.50\textwidth}
		\includegraphics[width=0.96\textwidth]{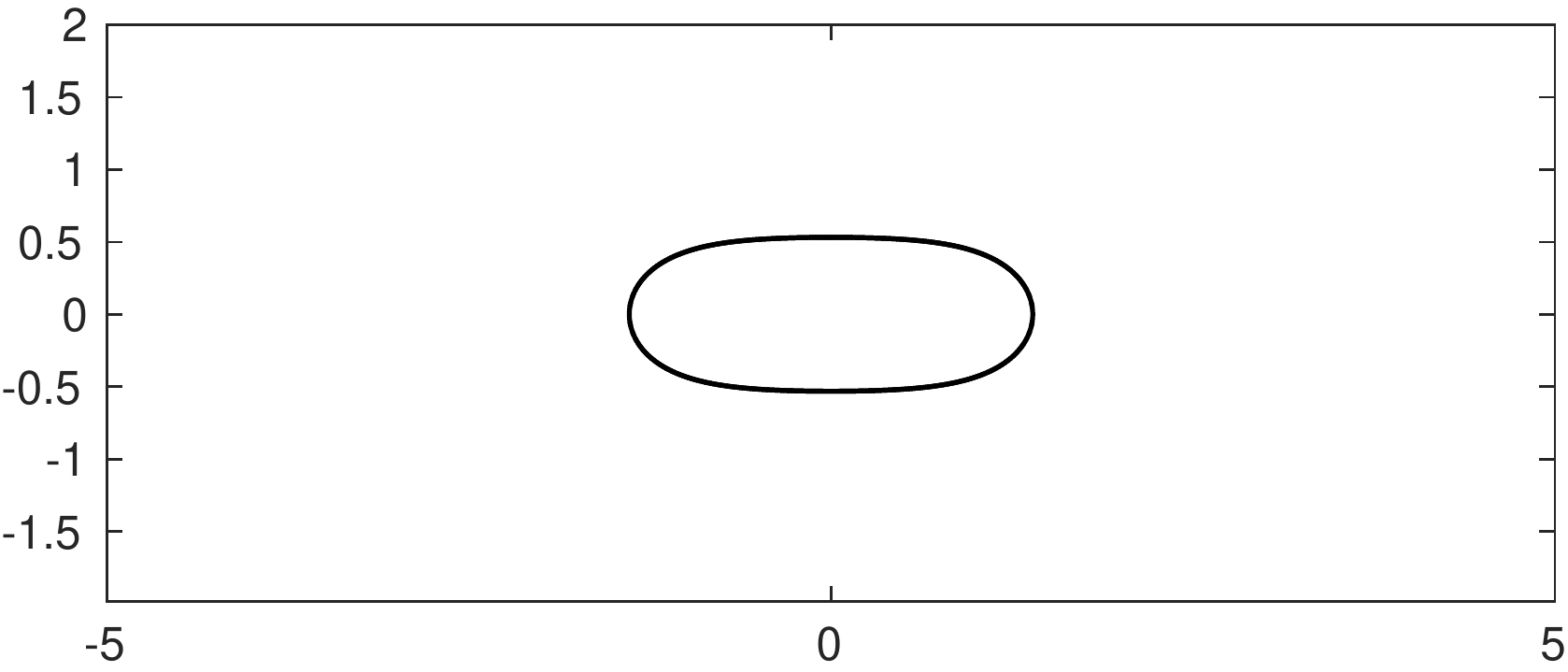}
	\end{minipage}%
	\caption{Mean curvature flow of a dumbbell-shaped curve in 2D. From the top left to the bottom right, the plots show the evolution at times $t=0$, $0.01$, $0.50$, $1.25$.}
	\label{fig:mcm 2d dumbbell}
\end{figure}

Our first example is the motion by mean curvature of a dumbbell-shaped curve in 2D; see Fig.\ \ref{fig:mcm 2d dumbbell}. From the initial dumbbell shape we see that all the corners smooth out rapidly. Then as the evolution continues, the curve shortens, and if given enough time will eventually collapse down to a point.
Our reference solution was computed to time $T=1.25$ using an explicit Runge-Kutta method and a small time step-size
(for standard explicit methods, the number of time steps needed for stability is on the order of $10^4$).
In all computations, spatial derivatives are approximated using second order centered differences on
a periodic grid of size $256\times 512$.

\begin{figure}[htb!]
	\centering
	\begin{minipage}{0.50\textwidth}
		\includegraphics[width=0.96\textwidth]{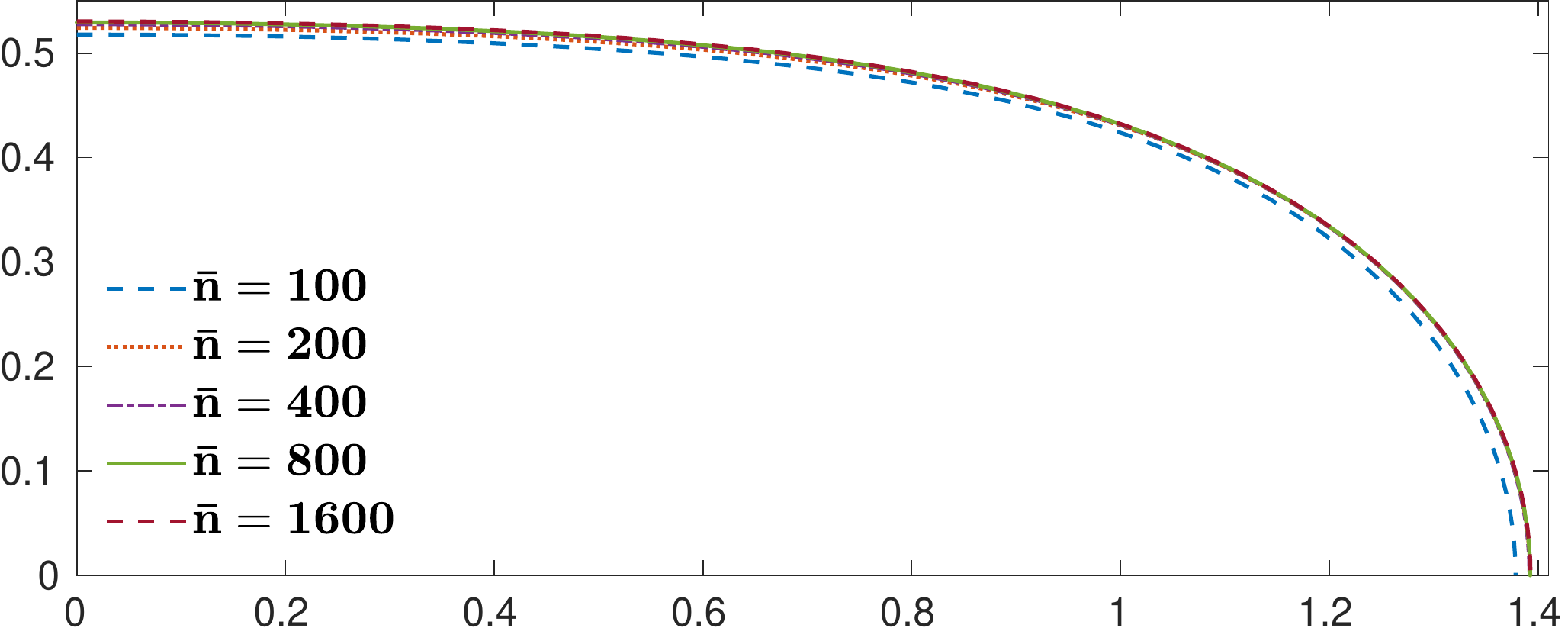}
	\end{minipage}%
	\begin{minipage}{0.50\textwidth}
		\includegraphics[width=0.96\textwidth]{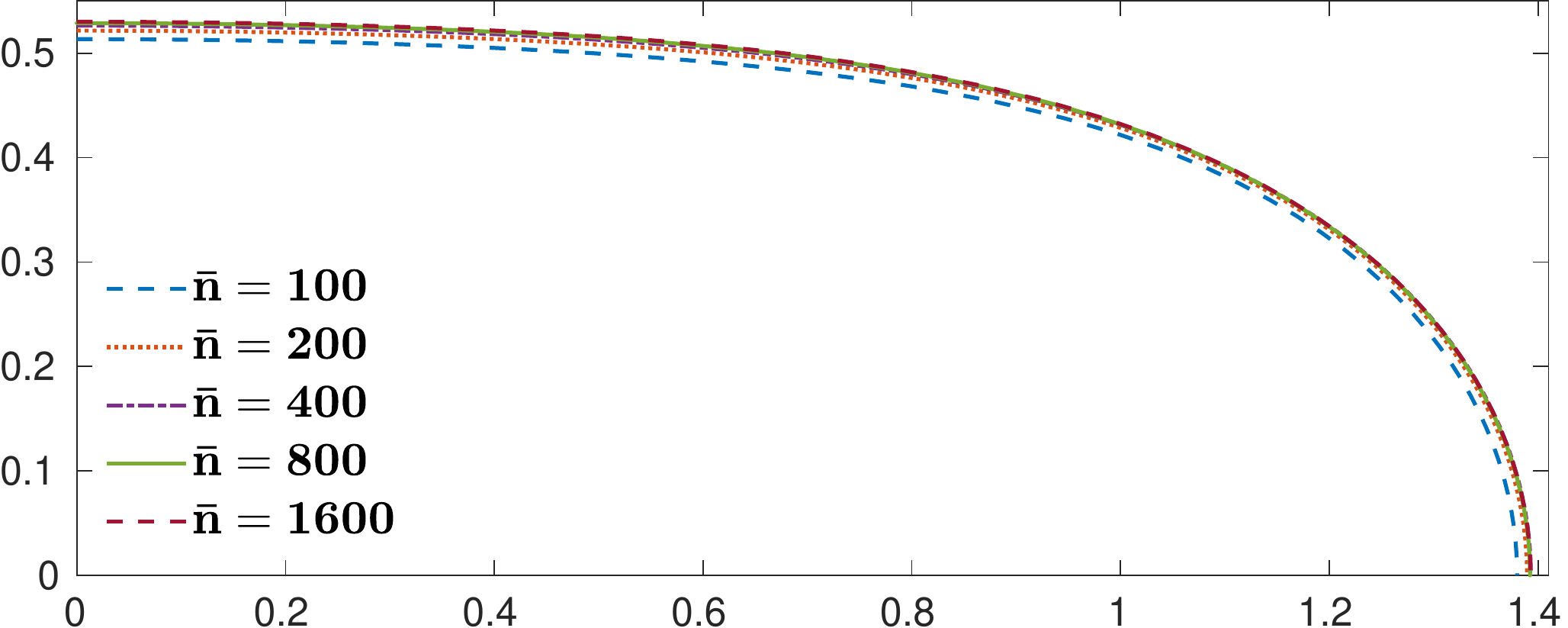}
	\end{minipage}
	
	\begin{minipage}{0.50\textwidth}
		\includegraphics[width=0.96\textwidth]{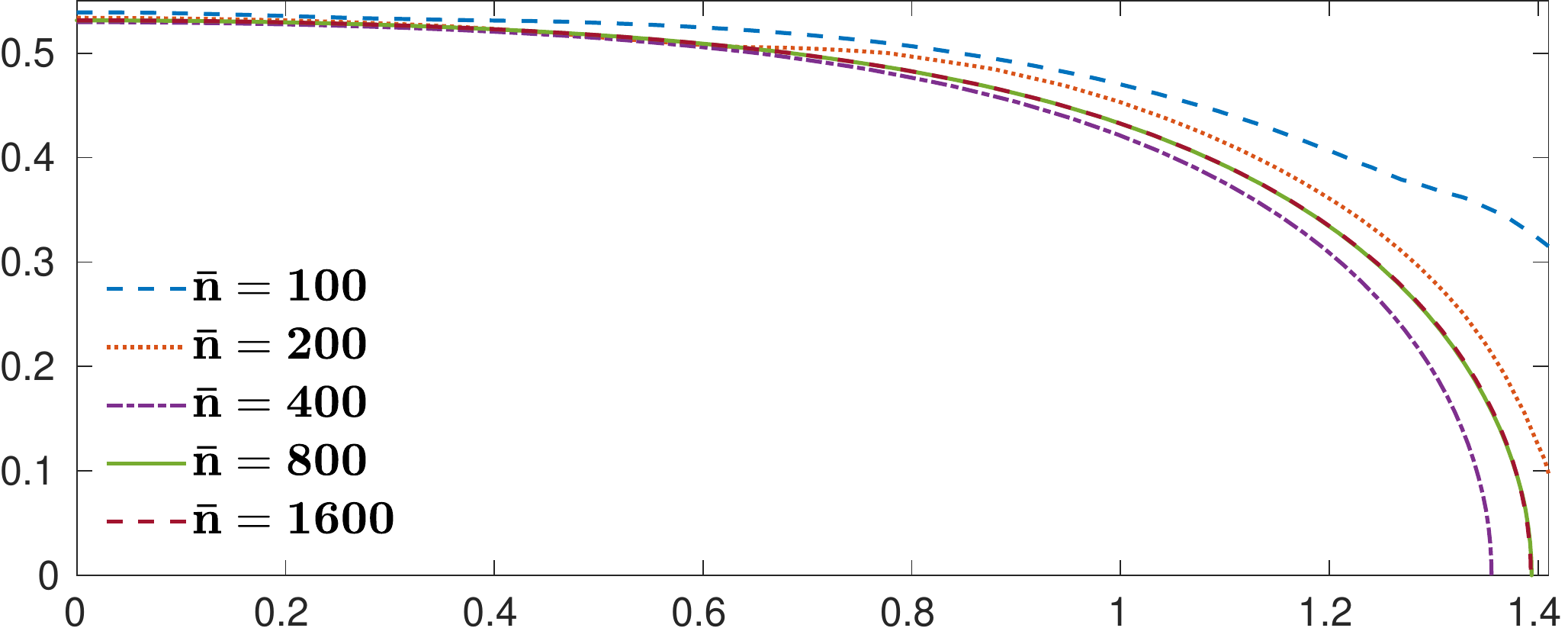}
	\end{minipage}%
	\begin{minipage}{0.50\textwidth}
		\includegraphics[width=0.96\textwidth]{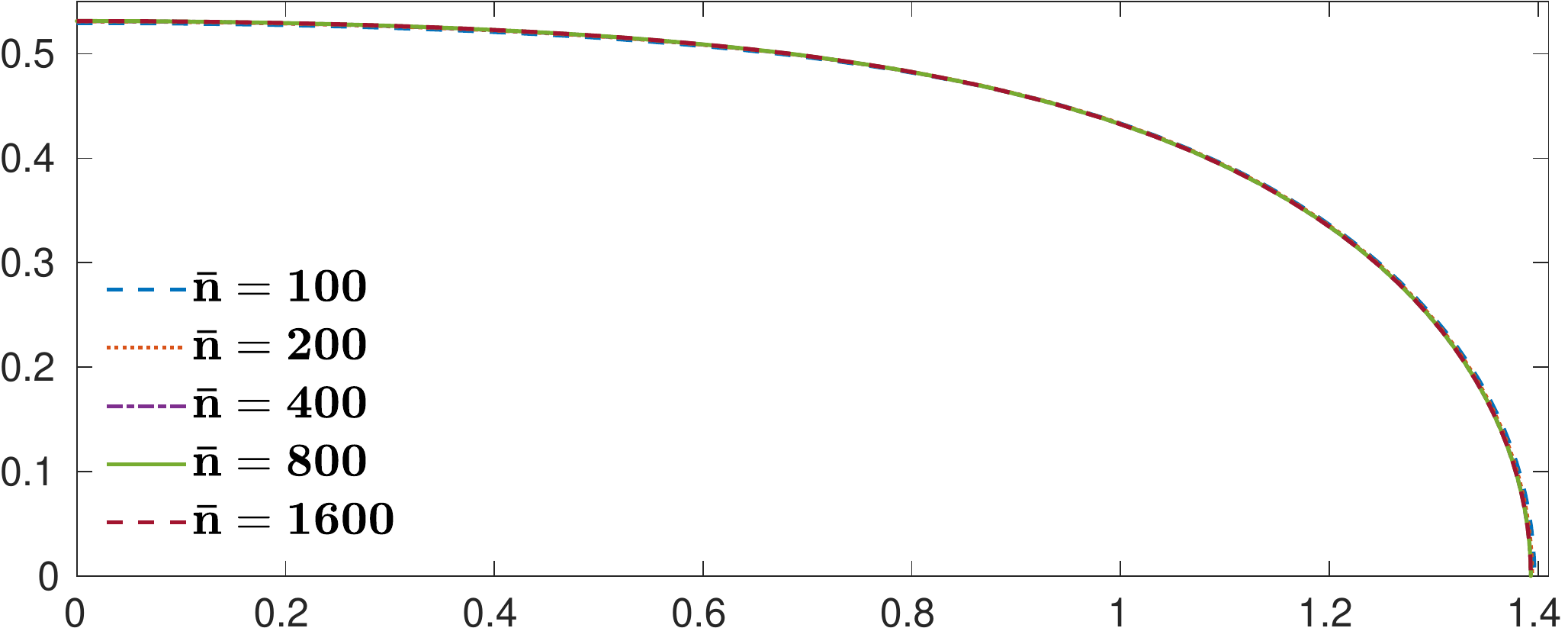}
	\end{minipage}
	\begin{minipage}{0.50\textwidth}
		\includegraphics[width=0.96\textwidth]{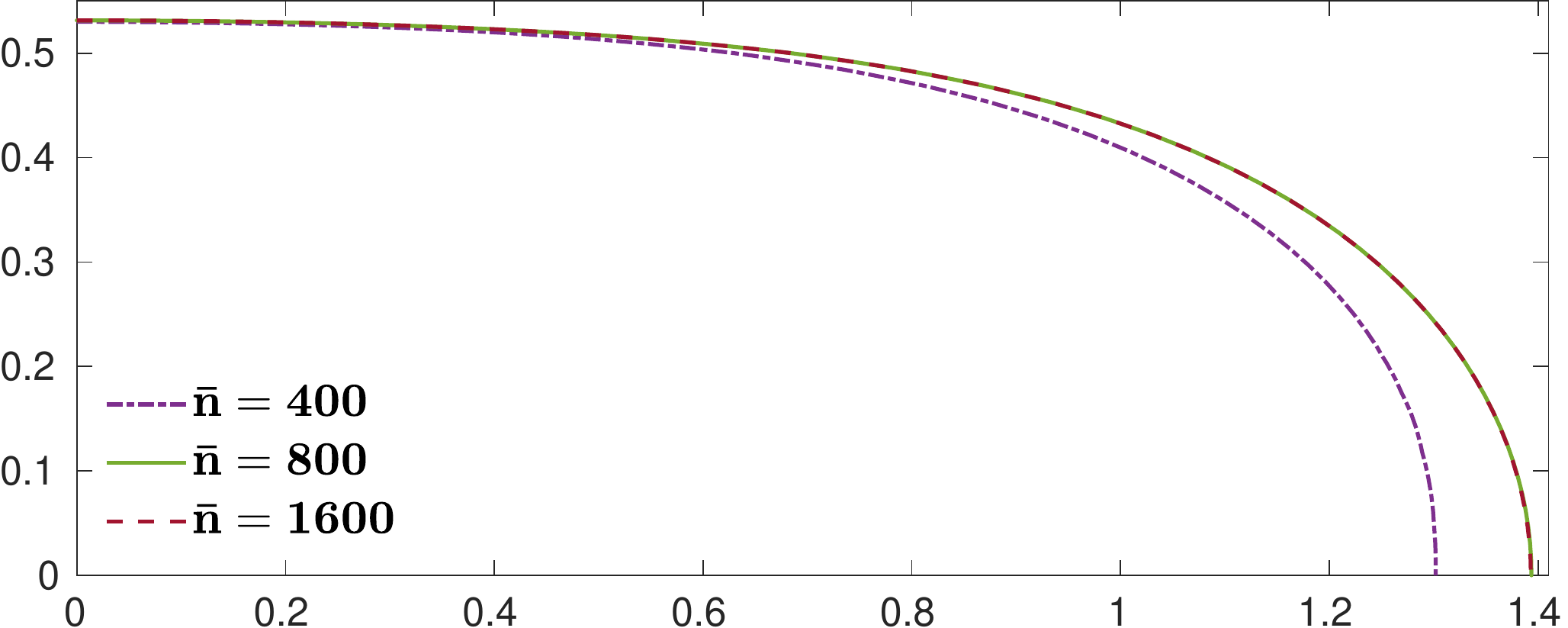}
	\end{minipage}%
	\begin{minipage}{0.50\textwidth}
		\includegraphics[width=0.96\textwidth]{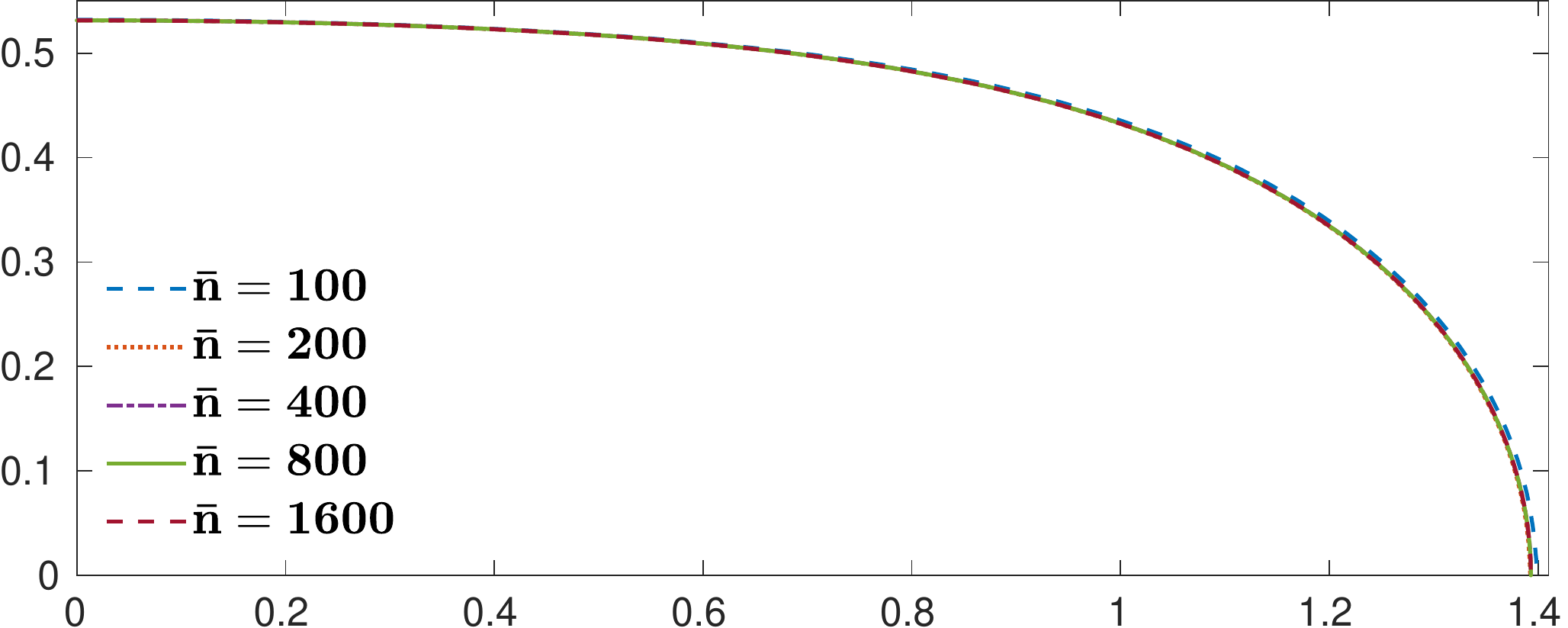}
	\end{minipage}
	\caption{Solution curves of various linearly stabilized time stepping schemes for the problem of curvature motion of a dumbbell-shaped curve in 2D. 
For clarity, we display only the first quadrant.  In the left column from top to bottom, we have SBDF1, SBDF2, and CNAB. 
In the right column from top to bottom, we have EIN, ETDRK2, and ETDRK4. 
Within each plot, each curve has an associated work estimate, $\bar{n}$, as indicated.}
	\label{fig:mcm 2d conv}
\end{figure}

Fig.~\ref{fig:mcm 2d conv} shows the convergence of SBDF1, SBDF2, CNAB, EIN, ETDRK2, and ETDRK4 for the curvature flow problem at 
the final time $T=1.25$. To compare performance, each solution curve is labeled with a work estimate, $\bar{n}$, where 
\begin{align*}
\bar{n} = (\#\text{ of time steps})\times (\#\text{ of RHS evals per time step}).
\end{align*}
The work estimate $\bar{n}$ is based on observations that RHS evaluation was the dominant cost in computing the solution to this problem in 2D and 3D, and that the relative computing times scaled proportionally.

\begin{figure}[htb!]
	\centering
	\includegraphics[width=84mm]{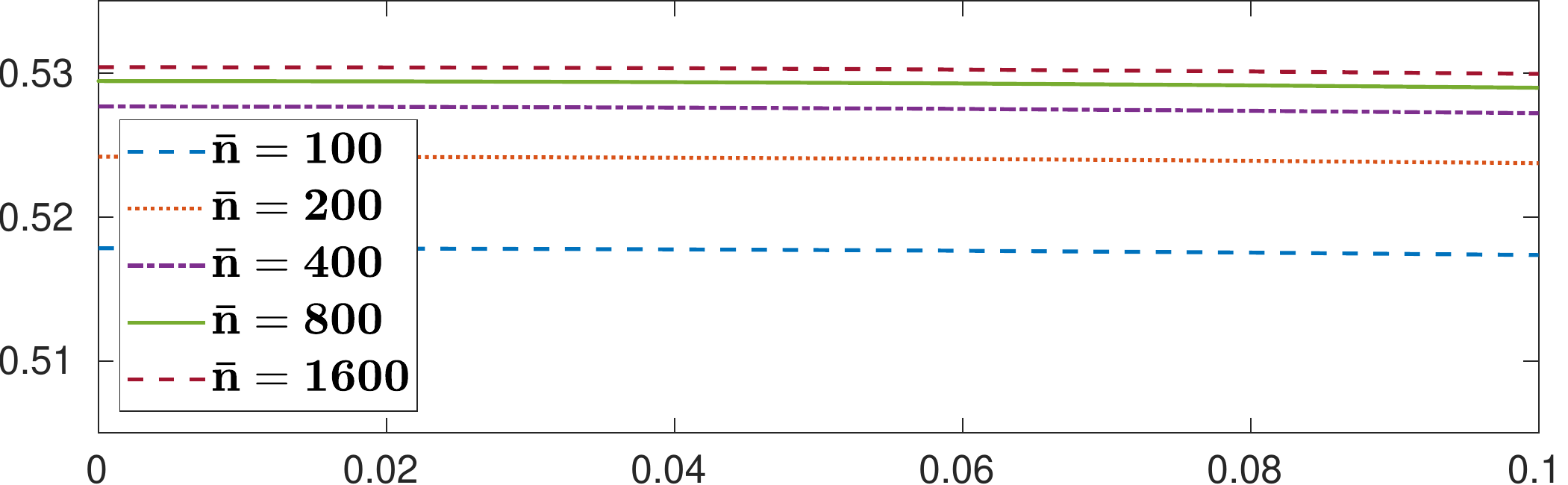}
	\includegraphics[width=84mm]{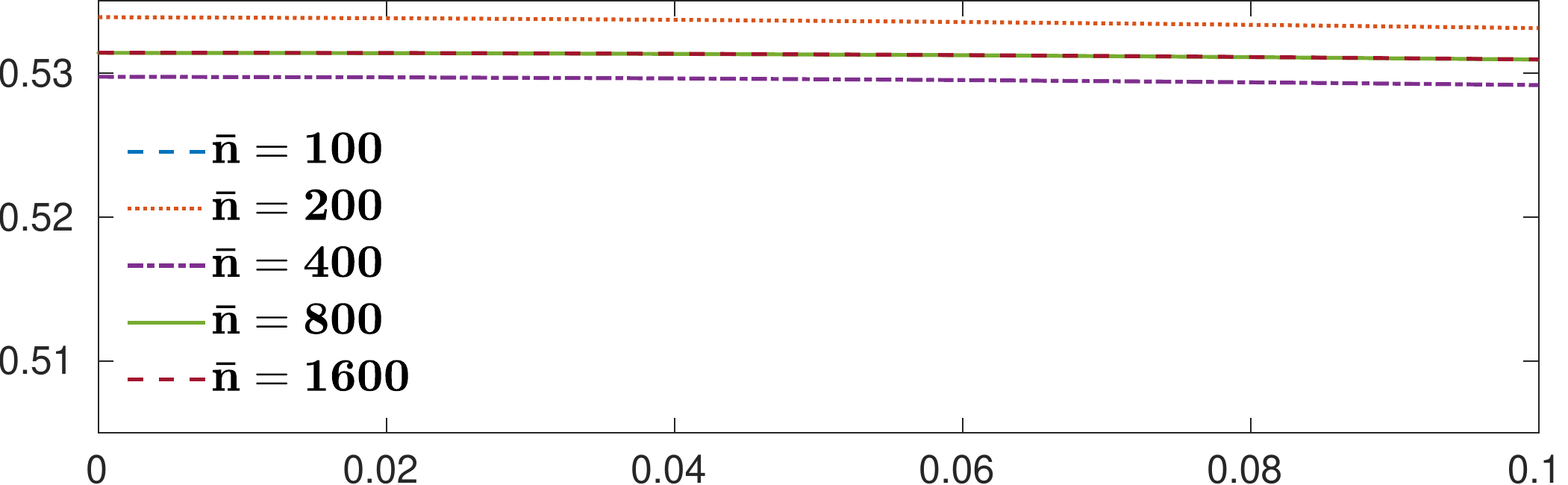}
	\includegraphics[width=84mm]{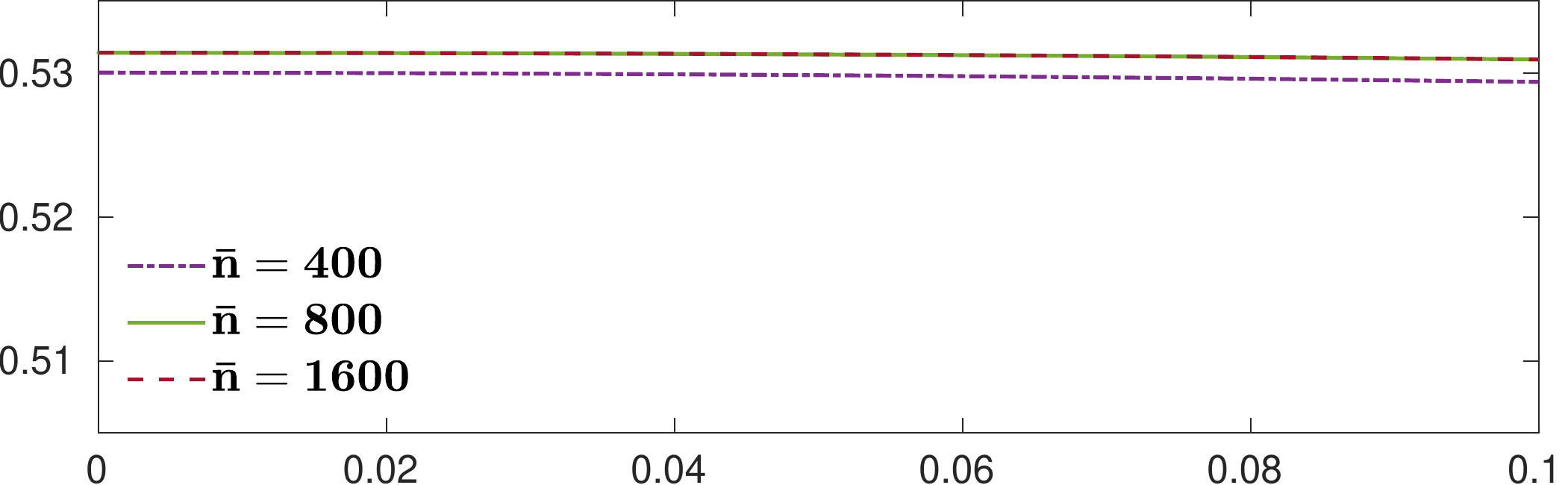}
	\includegraphics[width=84mm]{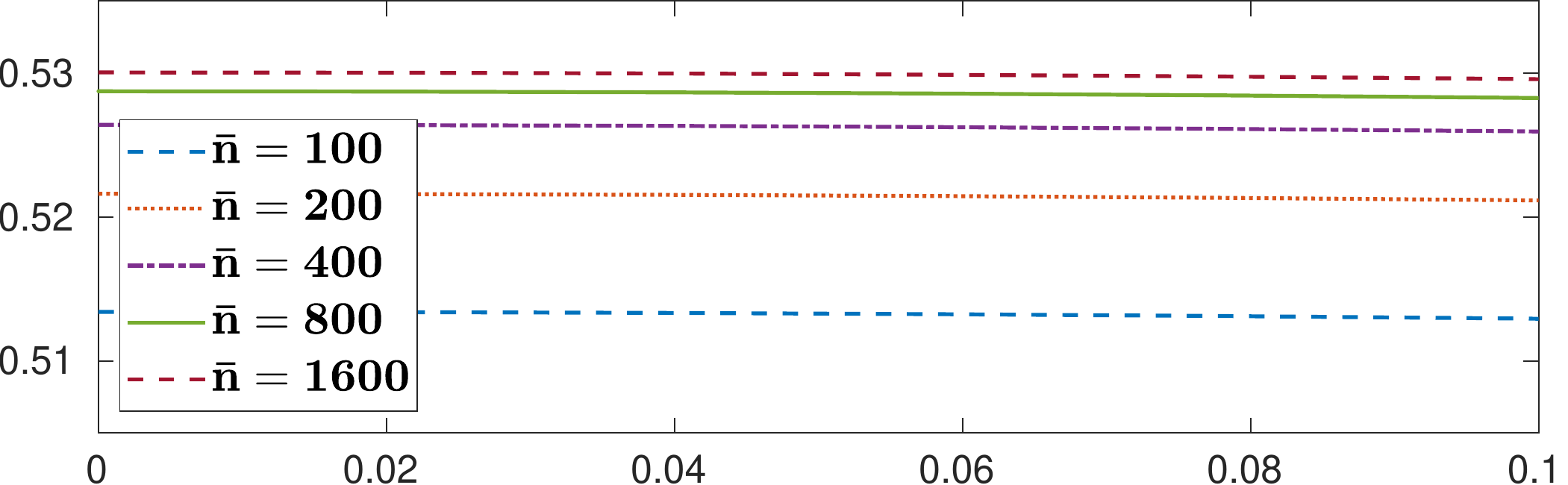}
	\includegraphics[width=84mm]{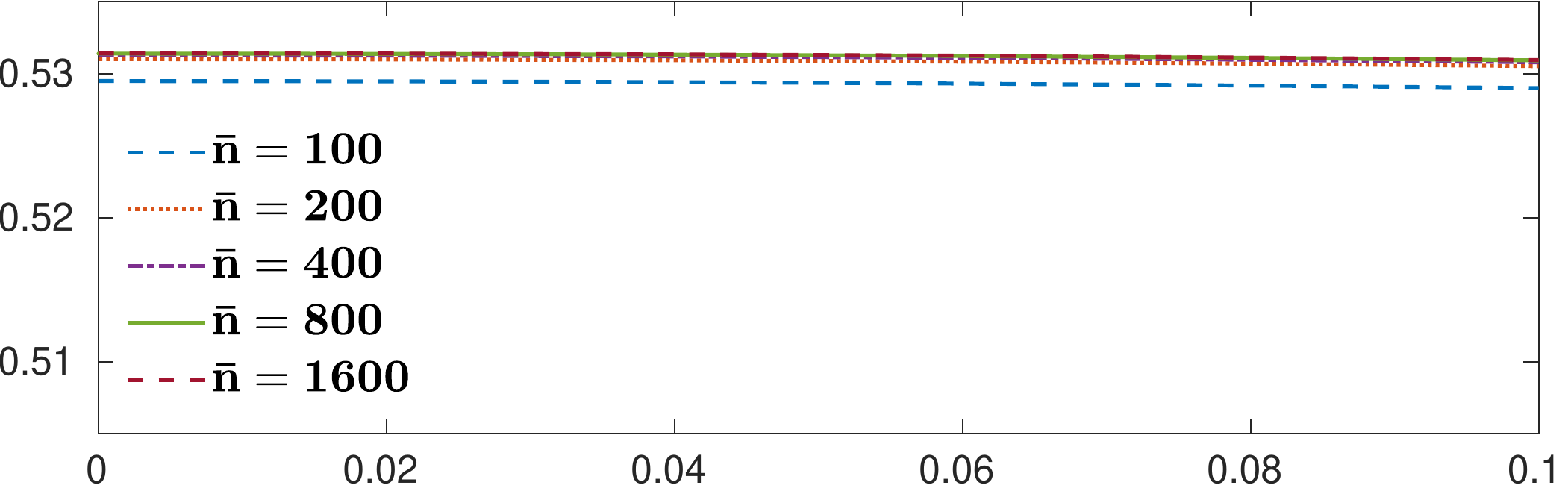}
	\includegraphics[width=84mm]{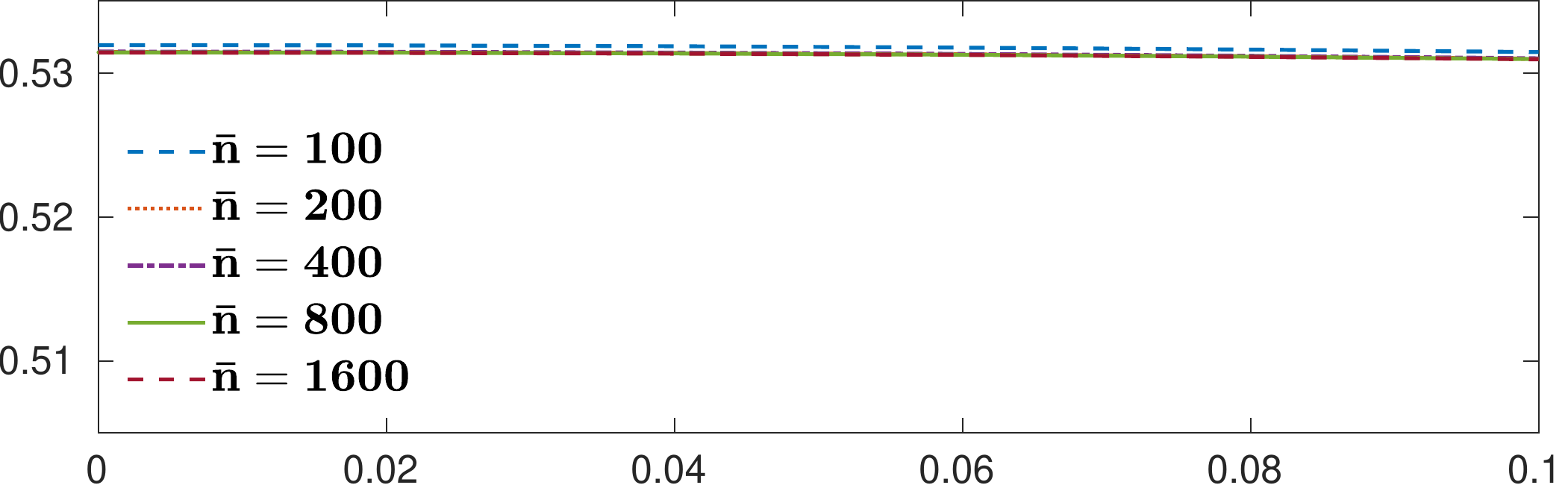}
	\caption{Zoom-in over $[0.0, 0.10]\times [0.505, 0.535]$ for the curvature flow of a dumbbell problem. 
From top to bottom we have  SBDF1, SBDF2, CNAB, EIN, ETDRK2, and ETDRK4. We see slow convergence of SBDF1 and the EIN method, and good convergence of the ETDRK schemes.}
	\label{fig:mcm zoom in 1}
\end{figure}

\begin{figure}[htb!]
	\centering
	\begin{minipage}{0.33\textwidth}
		\includegraphics[width=0.96\textwidth]{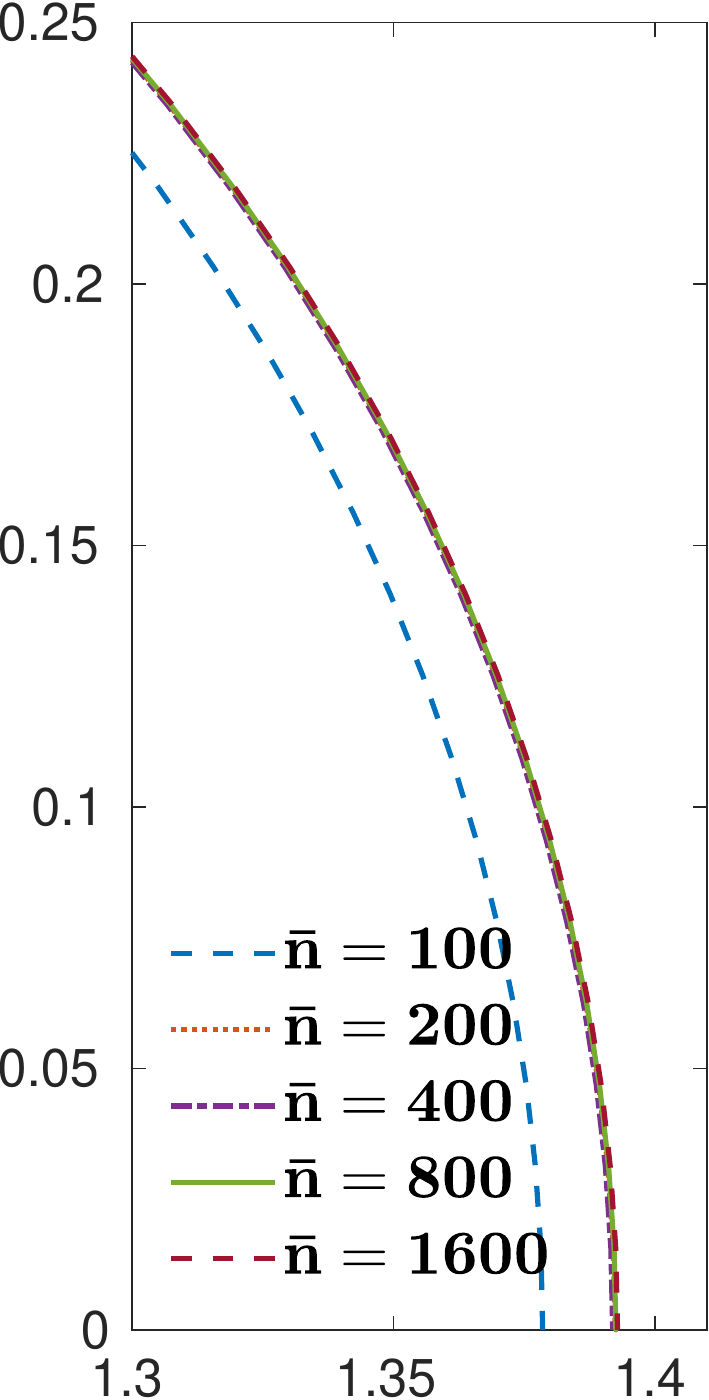}
	\end{minipage}%
	\begin{minipage}{0.33\textwidth}
		\includegraphics[width=0.96\textwidth]{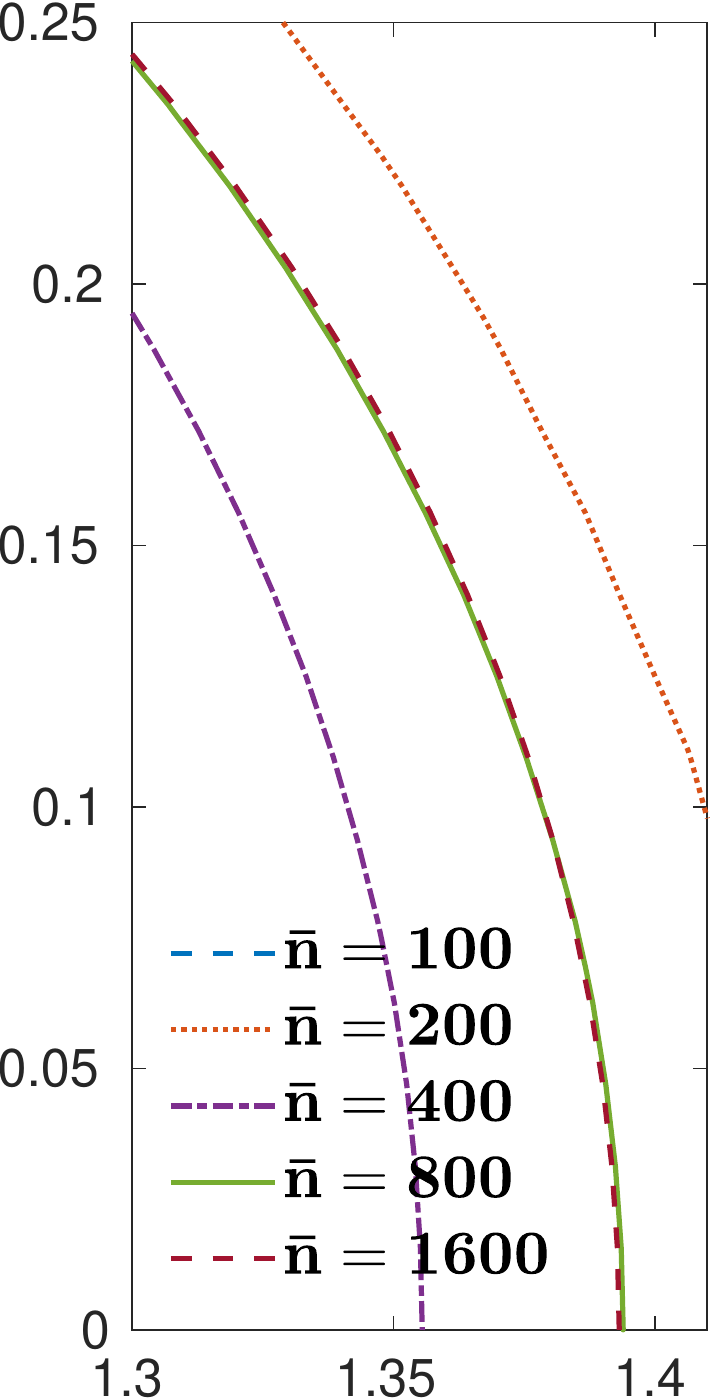}
	\end{minipage}%
	\begin{minipage}{0.33\textwidth}
		\includegraphics[width=0.96\textwidth]{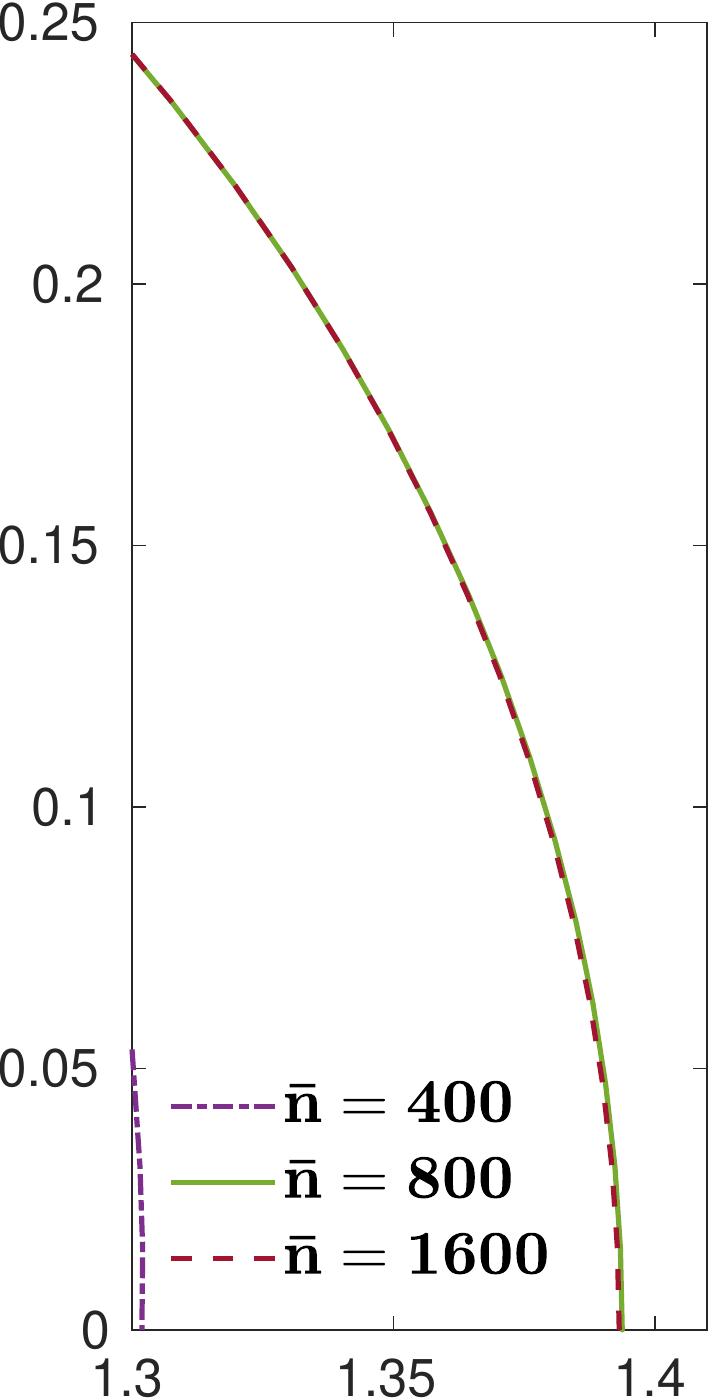}
	\end{minipage}
	\begin{minipage}{0.33\textwidth}
		\includegraphics[width=0.96\textwidth]{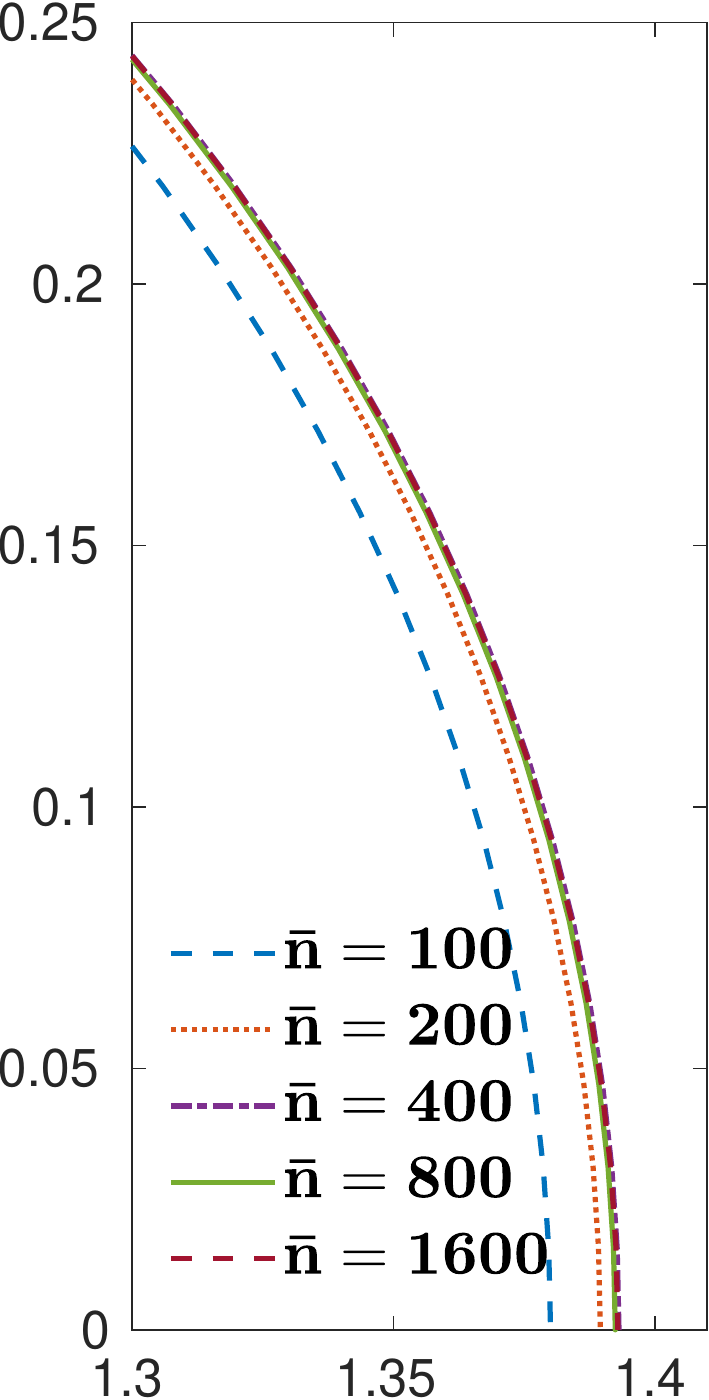}
	\end{minipage}%
	\begin{minipage}{0.33\textwidth}
		\includegraphics[width=0.96\textwidth]{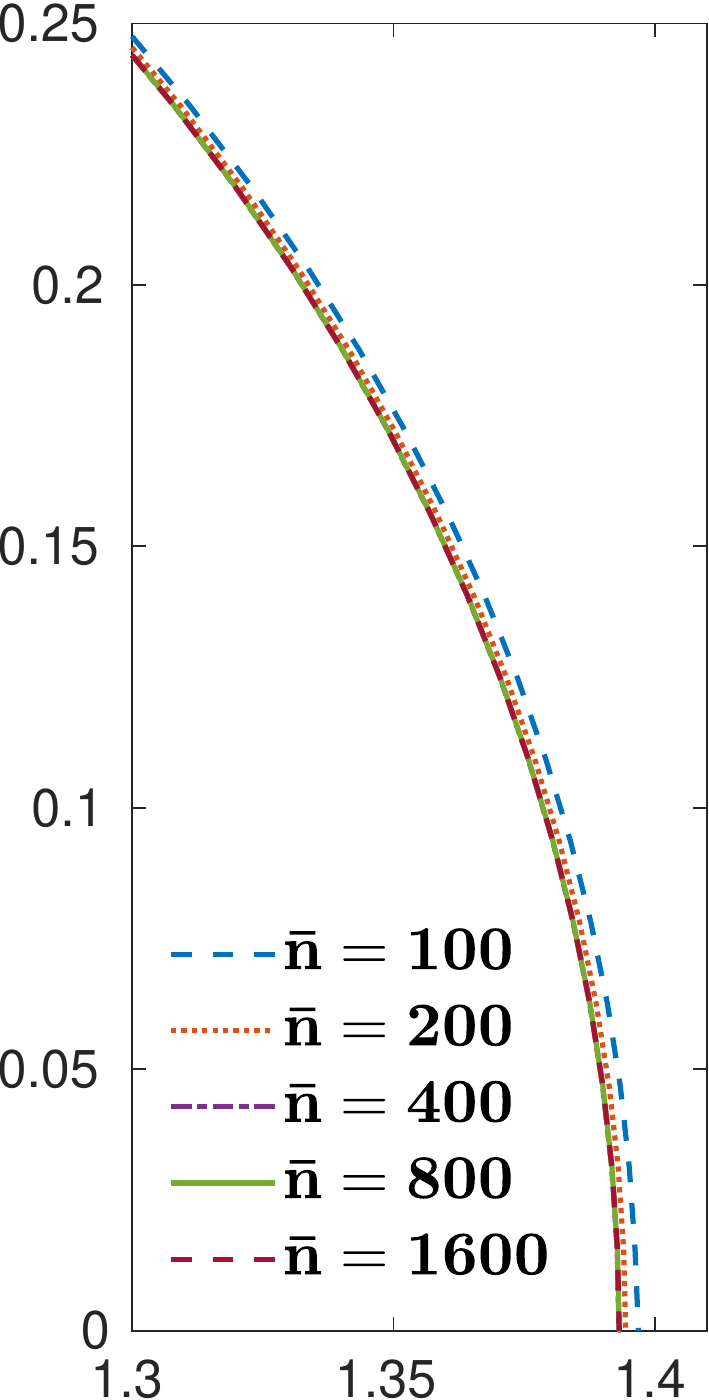}
	\end{minipage}%
	\begin{minipage}{0.33\textwidth}
		\includegraphics[width=0.96\textwidth]{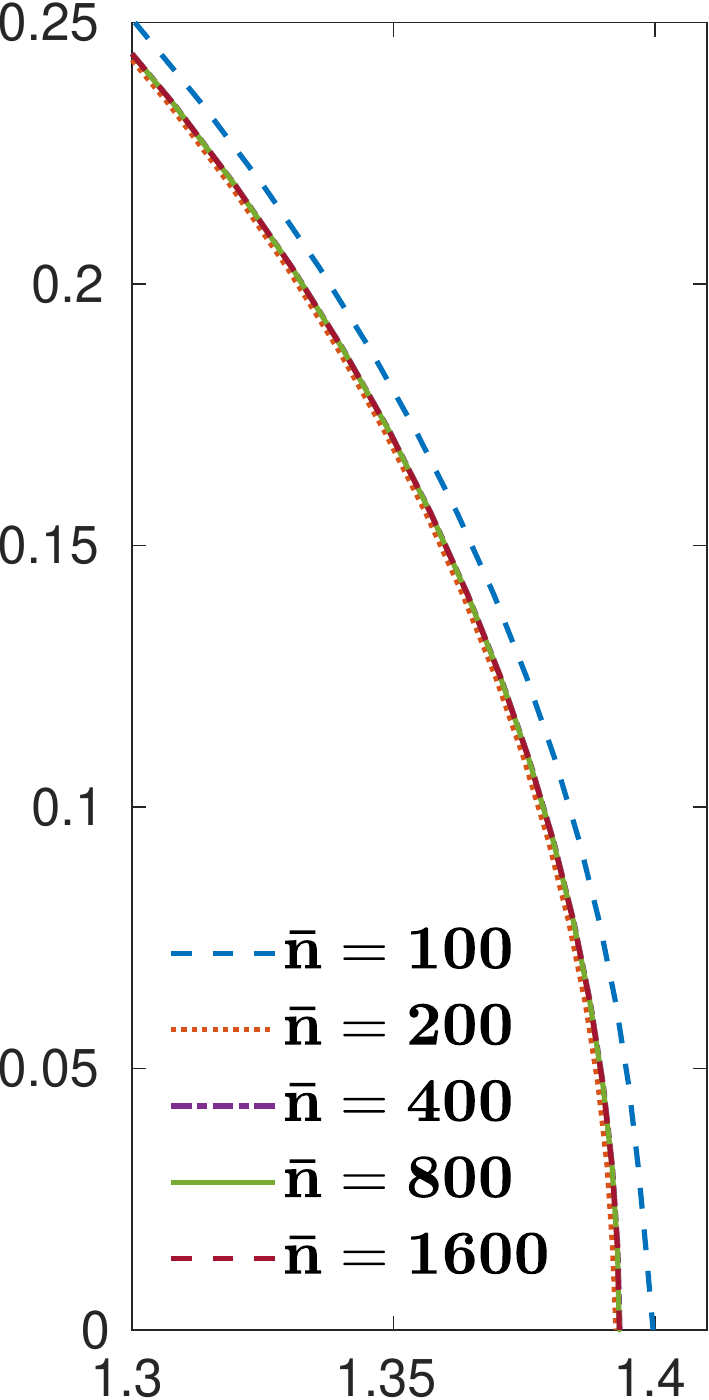}
	\end{minipage}
	\caption{Zoom-in over $[1.30,1.41]\times [0.0, 0.25]$. Top row: SBDF1, SBDF2, CNAB. Bottom row: EIN, ETDRK2, ETDRK4. We see slow convergence of the EIN method, and good convergence of the ETDRK schemes.}
	\label{fig:mcm zoom in 2}
\end{figure}

At a glance, Fig.~\ref{fig:mcm 2d conv} shows SBDF1, EIN, ETDRK2, and ETDRK4, the schemes that have strong damping as $z\to-\infty$ (see Fig.\ \ref{fig:damp fac at inf etd}), performing well at large step-sizes. In Figs.~\ref{fig:mcm zoom in 1} and \ref{fig:mcm zoom in 2}, we take a closer look at the performance of each scheme. The zoom-ins show that the ETDRK schemes offer the best accuracy and convergence for a given amount of work. 
SBDF1 offers surprisingly good results that are intermediate in quality between the EDTRK schemes and the EIN method.

The observed performance of SBDF2 and CNAB requires further comment. 
SBDF2 and CNAB were not competitive with the other schemes for $\bar{n}<800$.   
In particular, SBDF2 and CNAB exhibited poor accuracy for $\bar{n}=100$ and $\bar{n}<400$, respectively.  
Upon refining the step-size so that $\bar{n}\ge800$,  SBDF2 and CNAB both outperform the results from SBDF1 and EIN, with the result of CNAB being somewhat more accurate. 
As a point of comparison, CNLF doesn't generate acceptable solution curves until $\bar n > 2150$.

\subsubsection{Shrinking dumbbell in 3D}
Next, we take this example into 3D to illustrate the speed up of linearly stabilized schemes over the standard choice of explicit time stepping schemes. In 2D, one could argue that the computations can be completed within reasonable computing times using forward Euler or an explicit Runge-Kutta method. In 3D, time step restrictions for explicit time stepping schemes may lead to excessively long computations that necessitate trade offs in the grid size, or computing only over very short times.

\begin{figure}[htb!]
	\centering
	\begin{minipage}{0.50\textwidth}
		\includegraphics[width=0.96\textwidth]{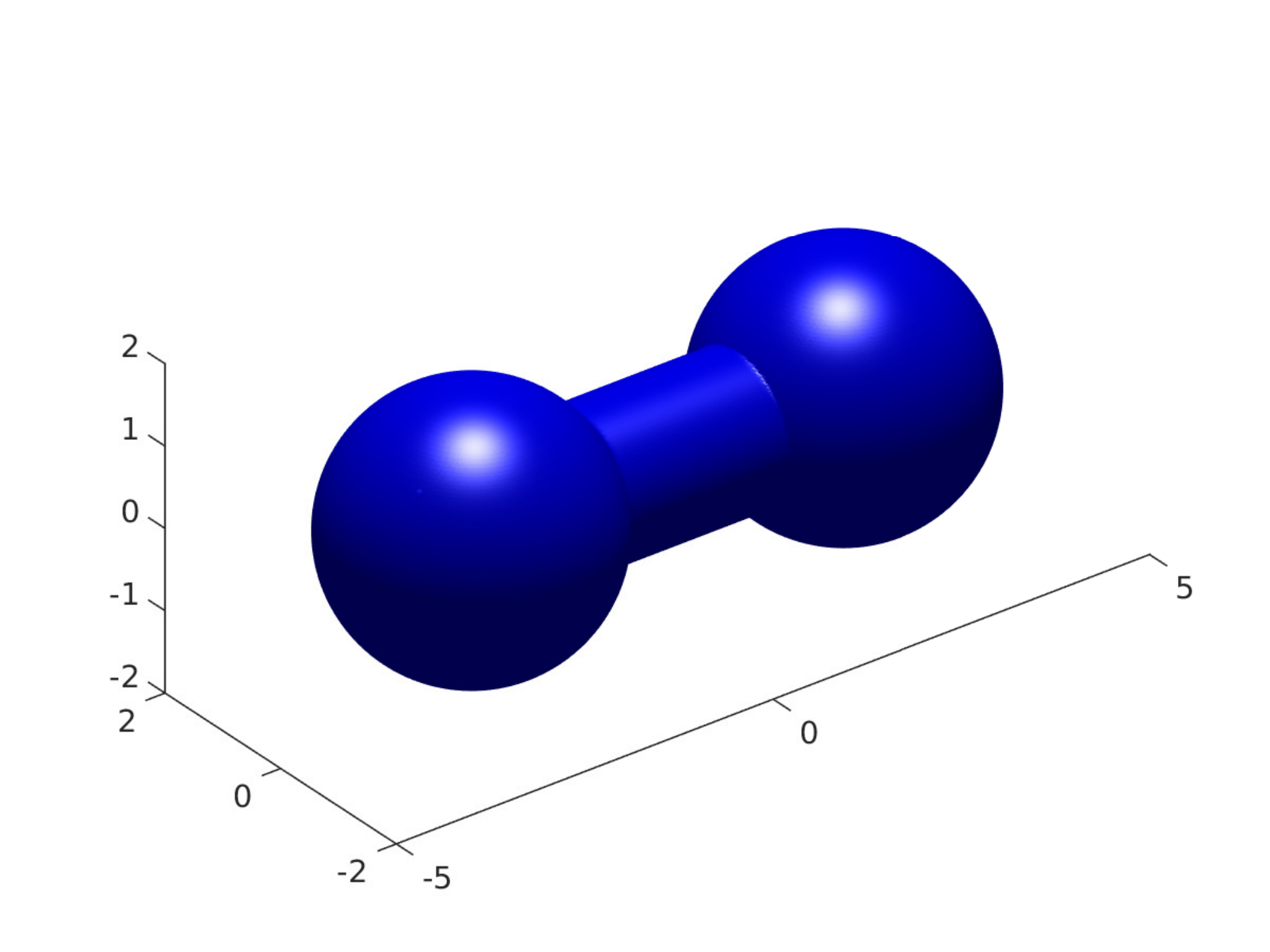}
	\end{minipage}%
	\begin{minipage}{0.50\textwidth}
		\includegraphics[width=0.96\textwidth]{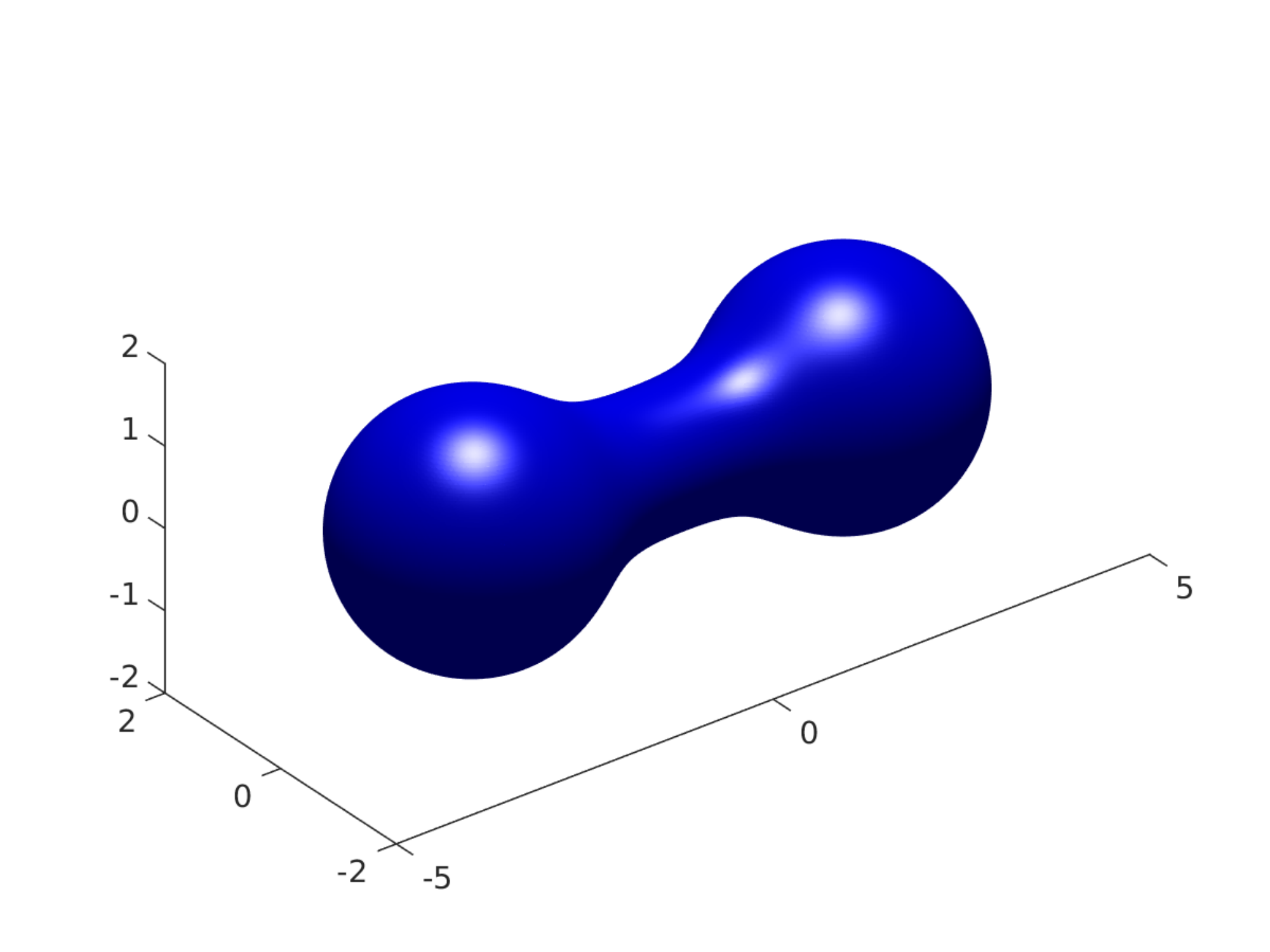}
	\end{minipage}
	\begin{minipage}{0.50\textwidth}
		\includegraphics[width=0.96\textwidth]{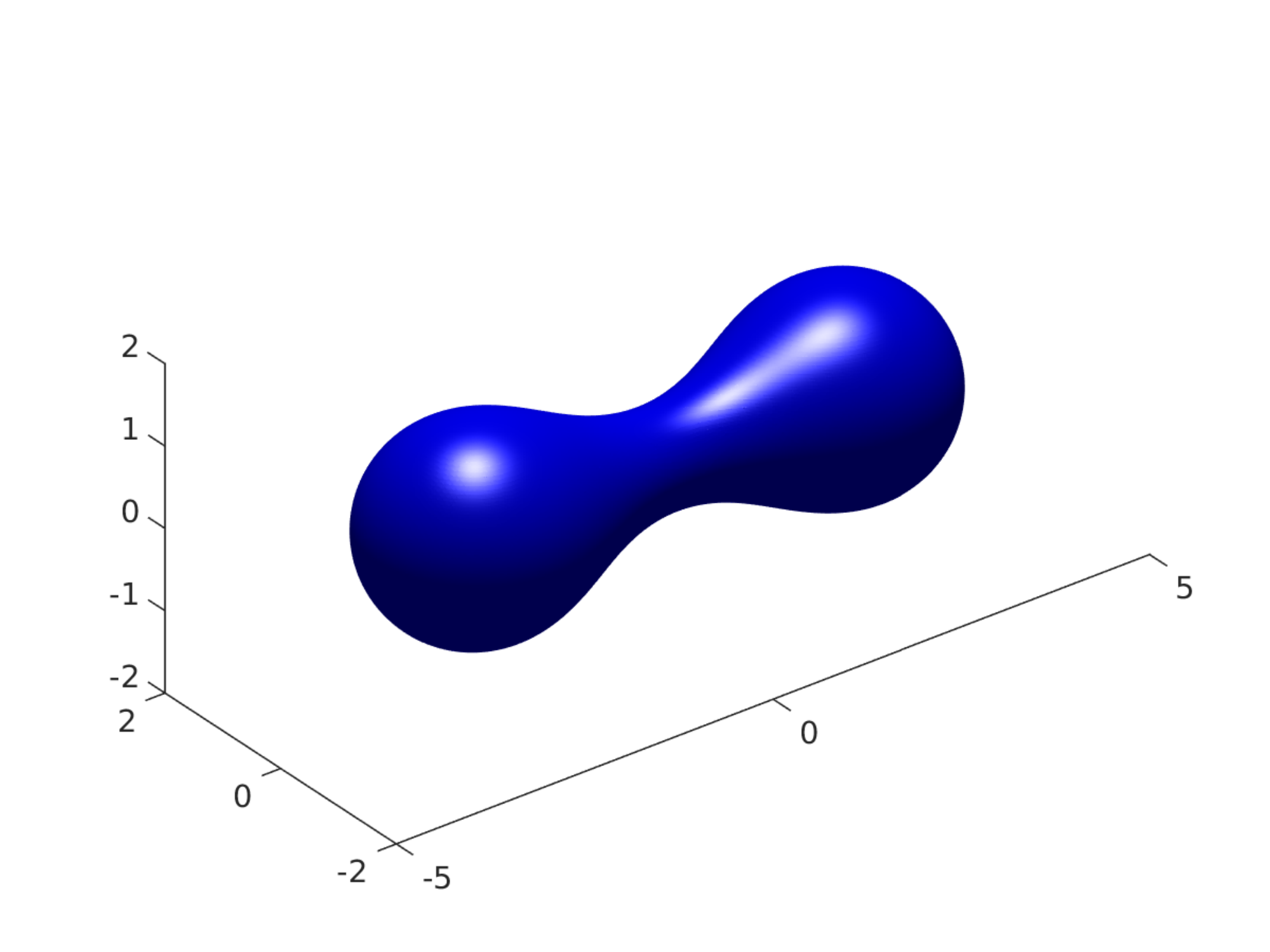}
	\end{minipage}%
	\begin{minipage}{0.50\textwidth}
		\includegraphics[width=0.96\textwidth]{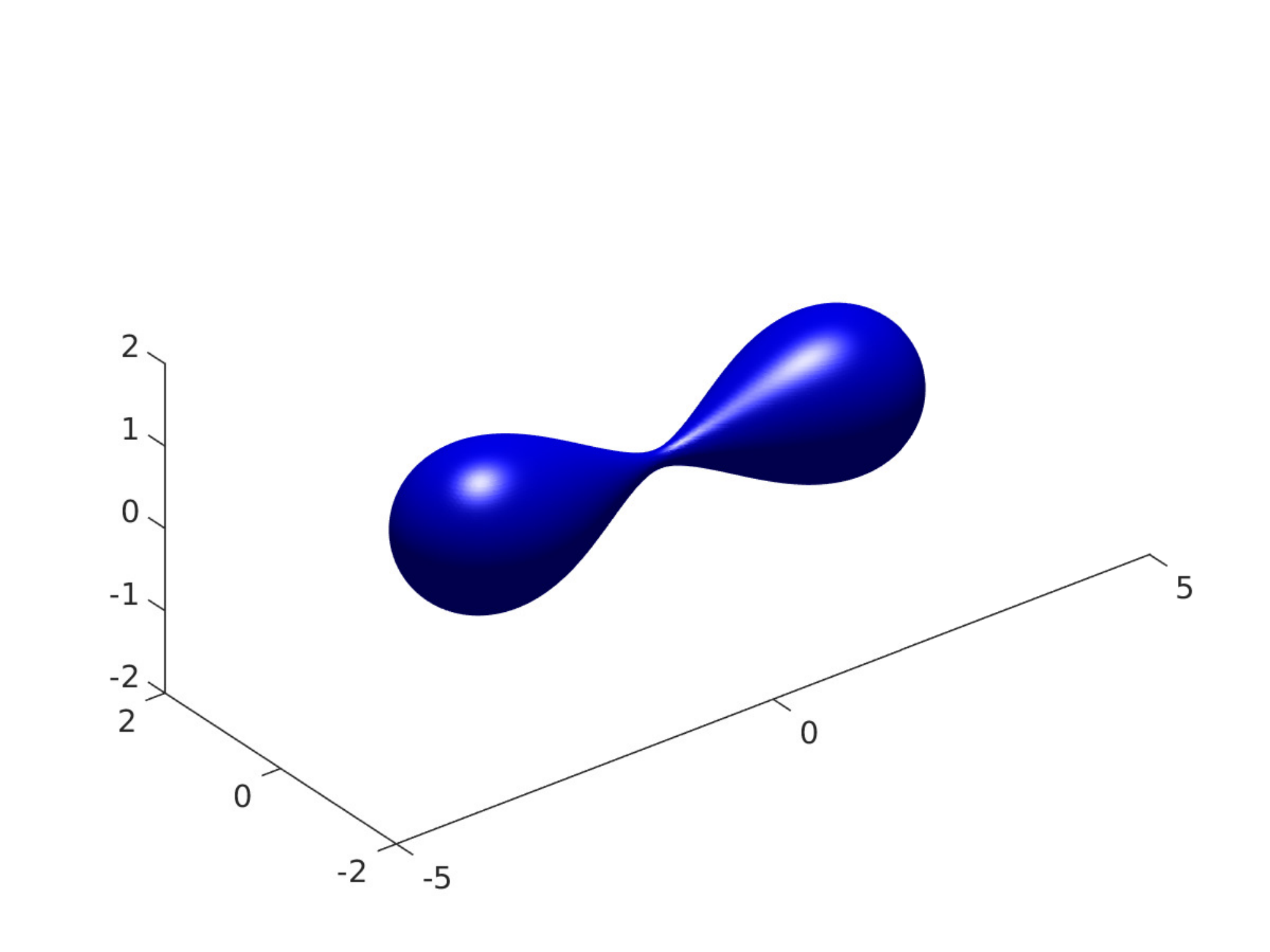}
	\end{minipage}
	\begin{minipage}{0.50\textwidth}
		\includegraphics[width=0.96\textwidth]{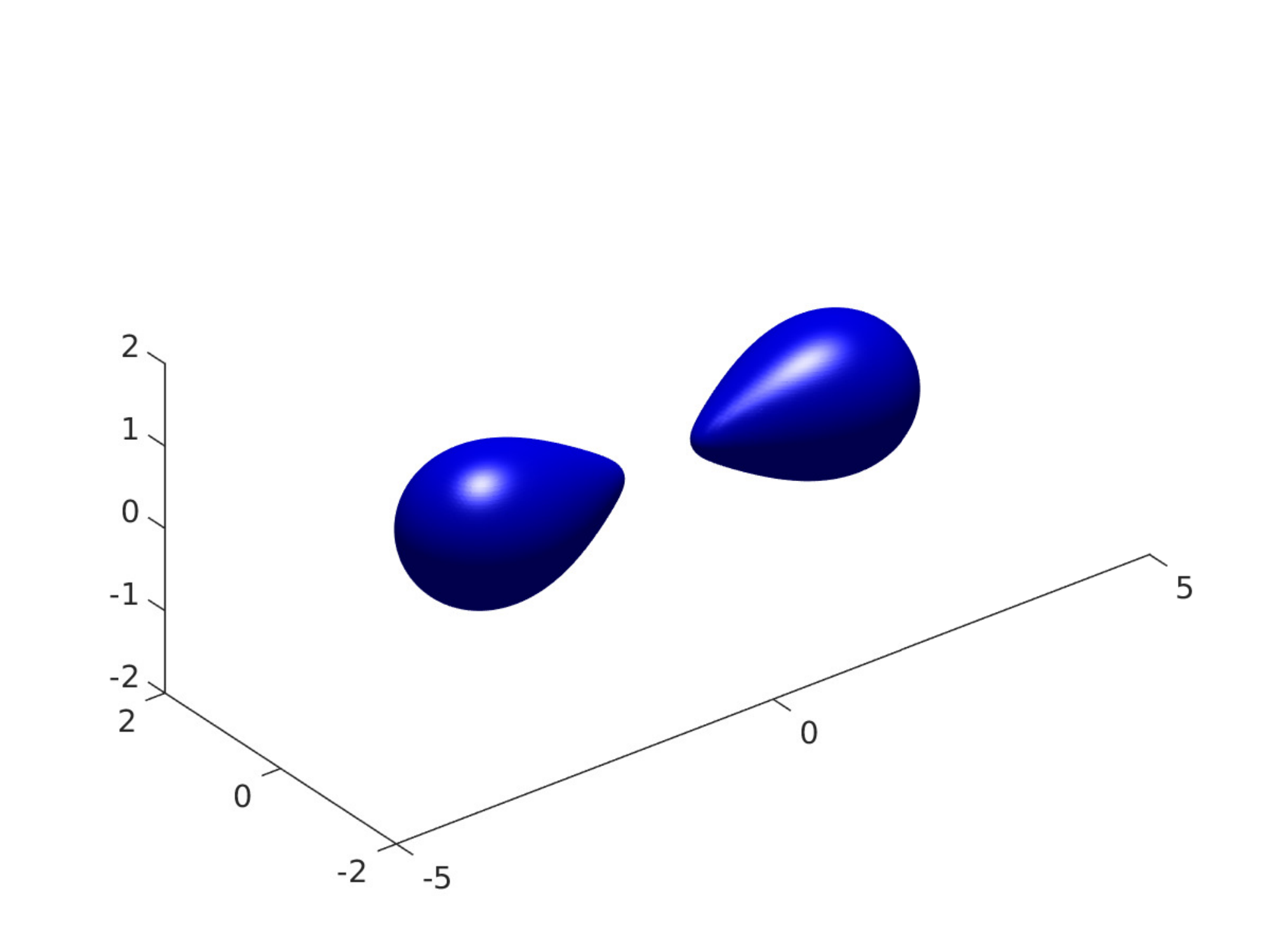}
	\end{minipage}%
	\begin{minipage}{0.50\textwidth}
		\includegraphics[width=0.96\textwidth]{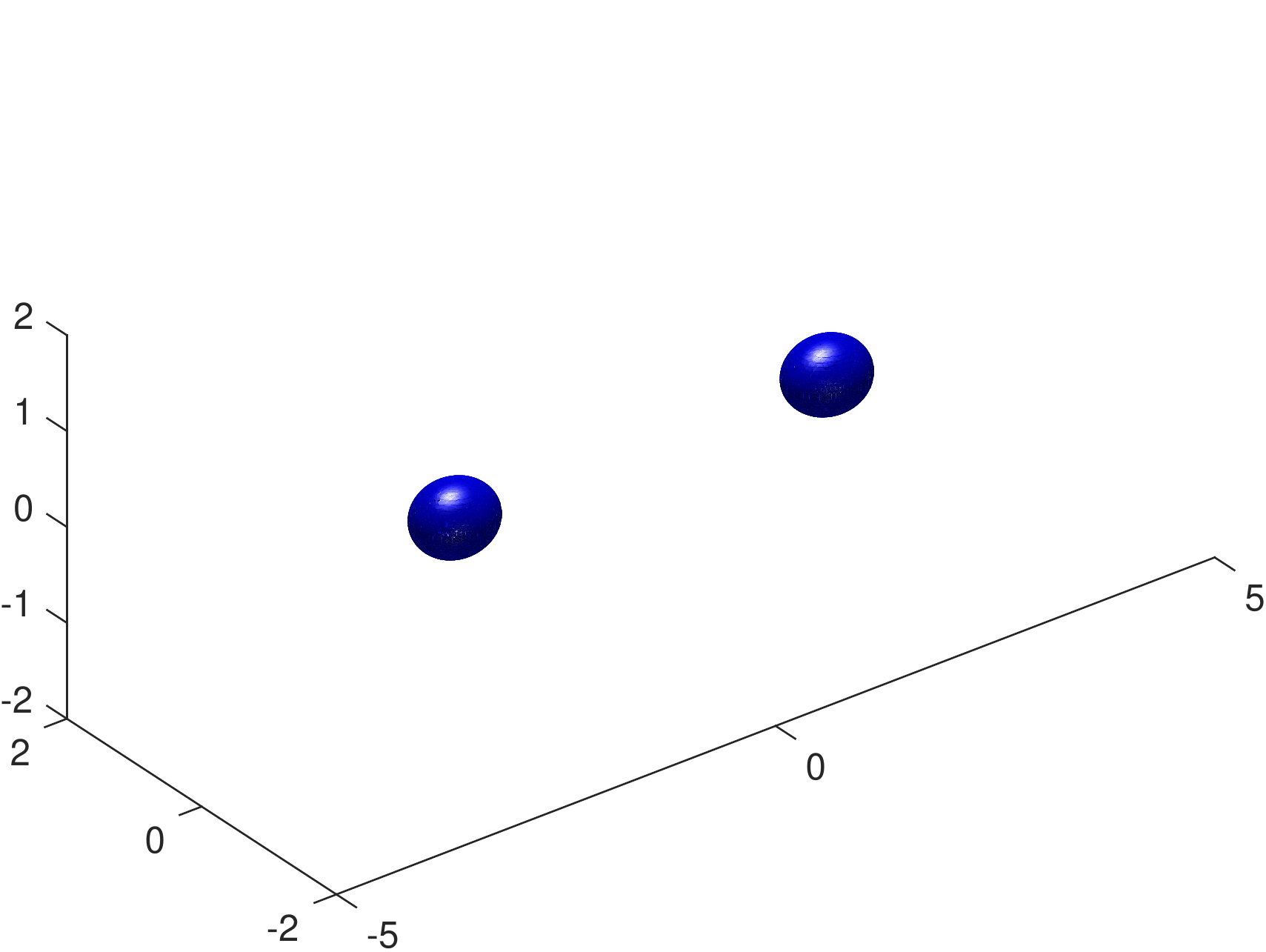}
	\end{minipage}
	\caption{Mean curvature flow of a dumbbell-shaped curve in 3D. From the top left to the bottom right, the plots show the evolution at times $t=0$, $0.10$, $0.30$, $0.525$, $0.55$, and $0.75$.}
	\label{fig:mcm 3d dumbbell}
\end{figure}

Setting the initial condition to be the dumbbell-shaped curve of the top left image in Fig.~\ref{fig:mcm 3d dumbbell}, the curve is then evolved under mean curvature flow. We use a periodic grid of size $256\times 128\times 128$ and solve to time $T=0.75$. With forward Euler, we needed $3000$ time steps for stability leading to a runtime of over 28 minutes in \textsc{Matlab} 2014b on an Intel\textsuperscript{\textregistered}Core\textsuperscript{\texttrademark}i5-4570 CPU@3.20GHz workstation running Linux. With the linearly stabilized ETDRK2, we solved the same problem using 80 time steps in under 2 minutes.

\subsubsection{Anisotropic mean curvature motion}
\label{sssect: aniso mcm}
 We conclude our examples by illustrating the use of a linearly stabilized scheme for an anisotropic motion.
In \cite{oberman2011aniso}, Oberman et al., present a method for anisotropic mean curvature flow: 
\begin{align}
	u_t = (\gamma(\omega) + \gamma''(\omega))\abs{\nabla u} \nabla \cdot \left( \frac{\nabla u}{\abs{\nabla u}} \right),
	\label{aniso mcm}
\end{align}
where $\omega = \arctan(u_y/u_x)$, and 
\begin{align}
	\gamma(\omega) 
	= \gamma_m(\omega) 
	= \frac{1}{m^2+1}(m^2 + 1 - \sin(m\omega)), 
	\quad\text{for } m = 0,2,4,8.
	\label{aniso mcm gamma}
\end{align}
Under isotropic mean curvature motion, a simple closed contour in 2D has a circular limiting shape as it reduces to a point. Under \eqref{aniso mcm} and \eqref{aniso mcm gamma}, the limiting shape will have $m$-fold rotational symmetry.

\begin{figure}
	\centering
	\begin{minipage}{0.50\textwidth}
		\includegraphics[width=0.96\textwidth]{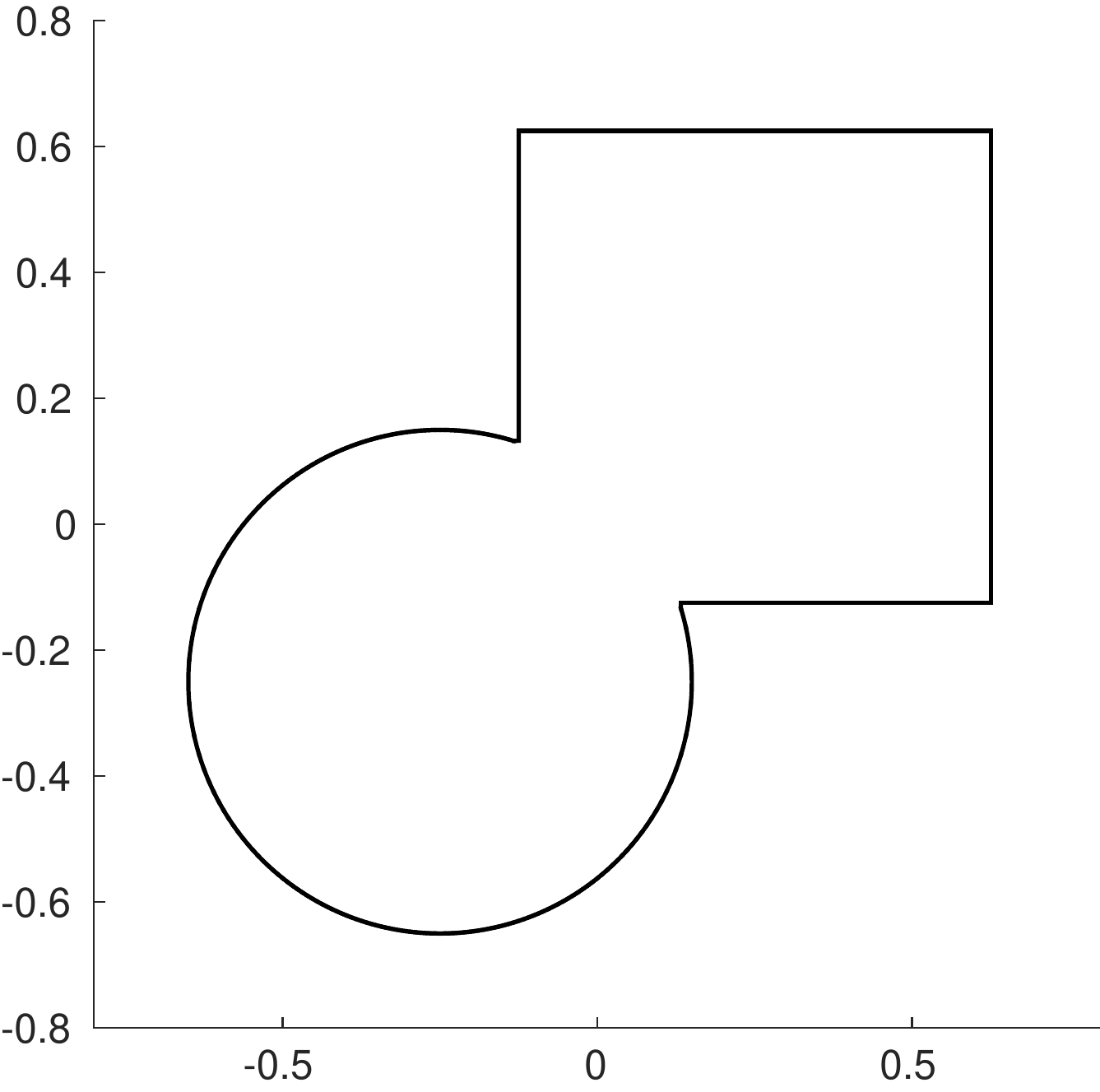}
	\end{minipage}%
	\begin{minipage}{0.50\textwidth}
		\includegraphics[width=0.96\textwidth]{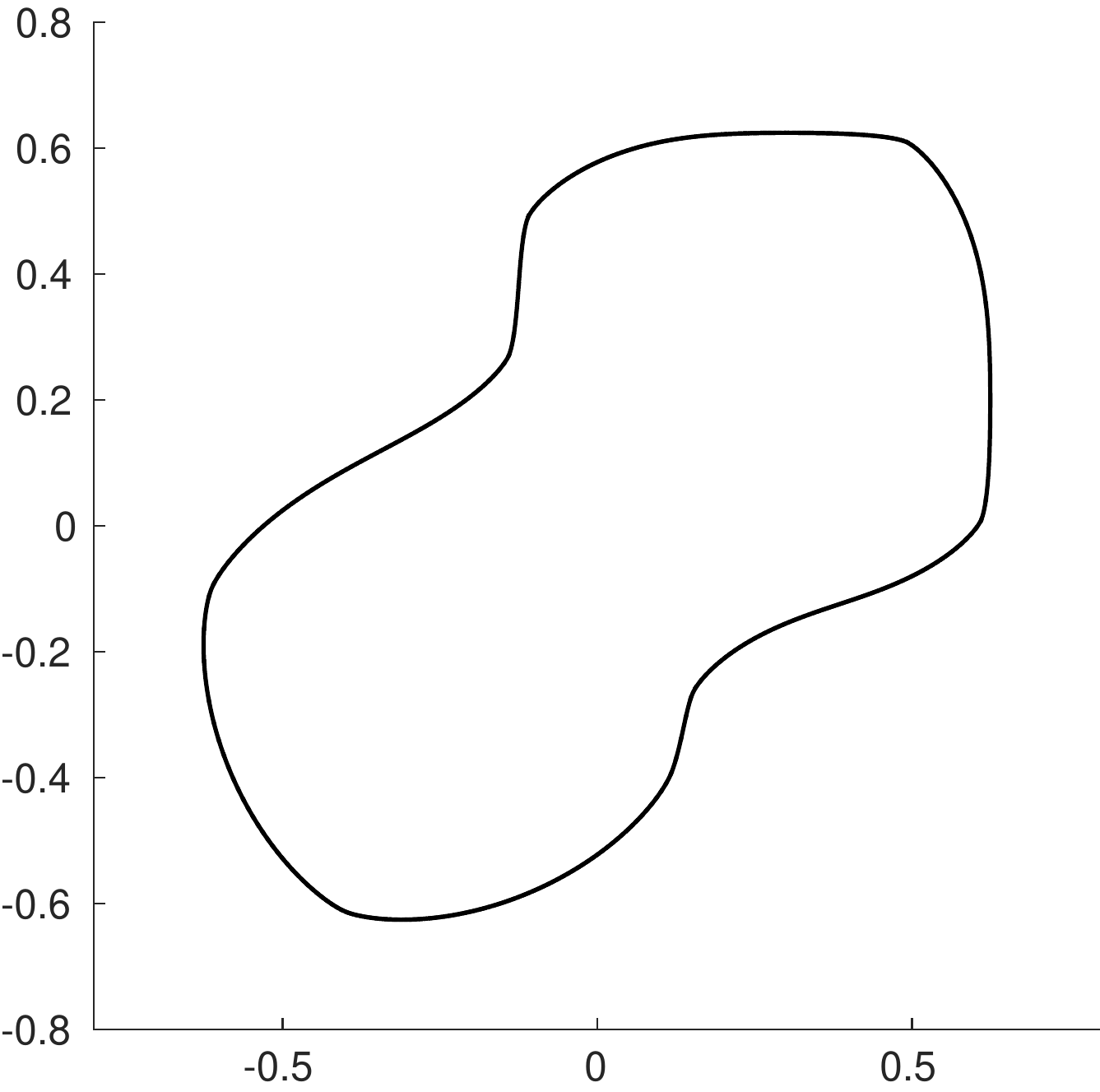}
	\end{minipage}
	\begin{minipage}{0.50\textwidth}
		\includegraphics[width=0.96\textwidth]{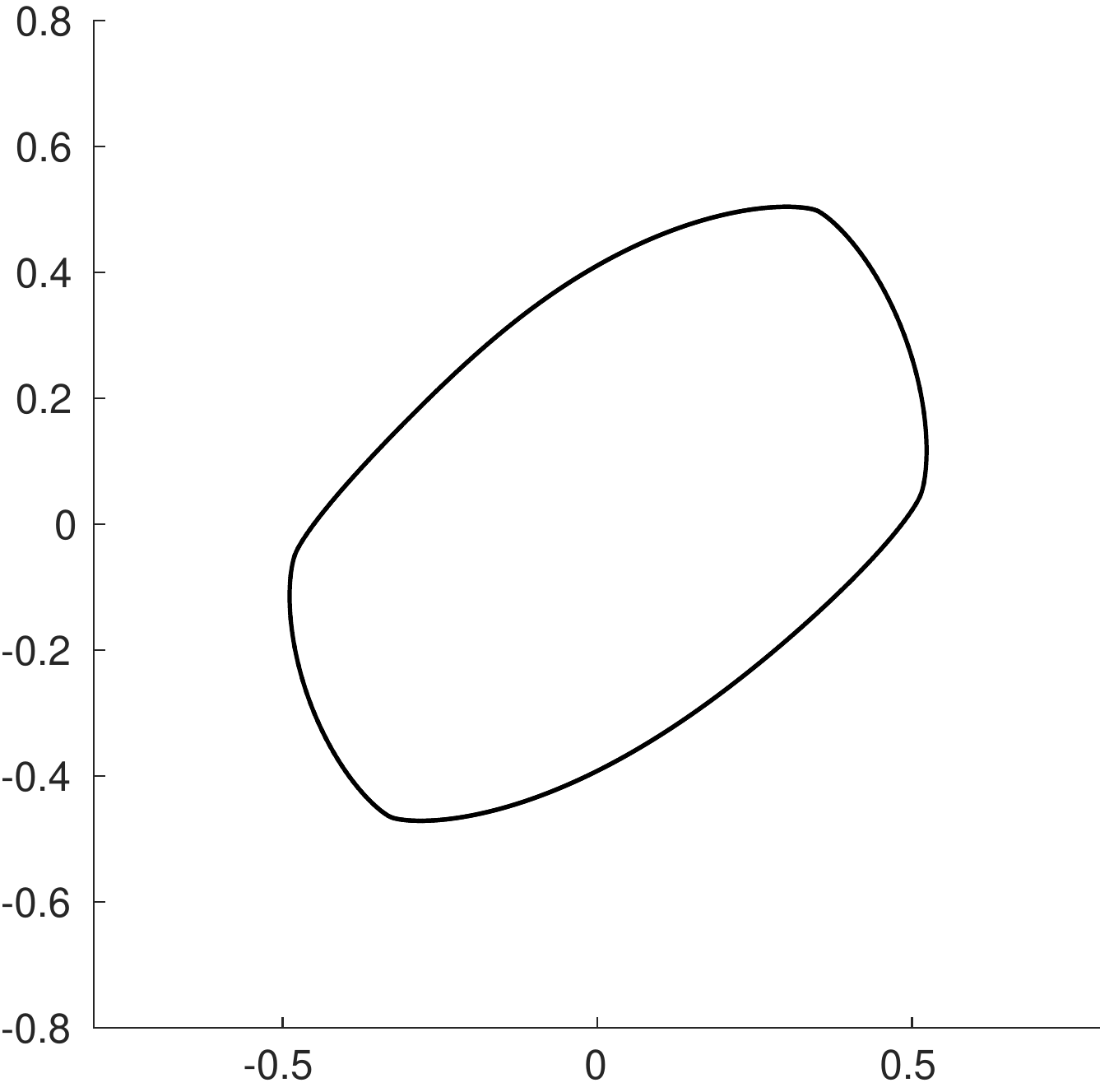}
	\end{minipage}%
	\begin{minipage}{0.50\textwidth}
		\includegraphics[width=0.96\textwidth]{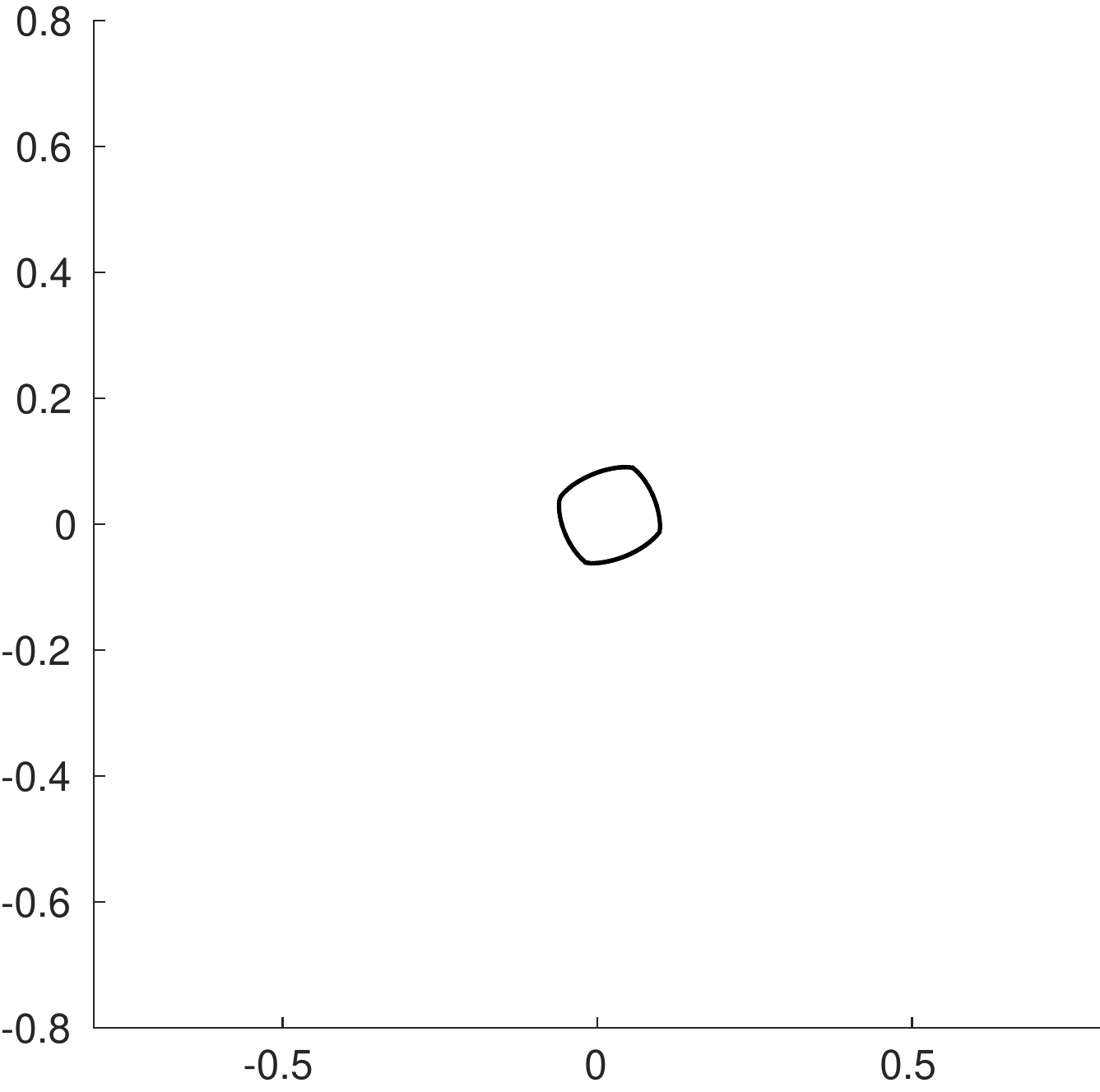}
	\end{minipage}
	\caption{Anisotropic mean curvature flow in 2D. The plots show the evolution of the curve at times $t=0$, $0.01$, $0.06$, $0.16$. The initial curve smooths and shrinks to a curve exhibiting four-fold symmetry as it collapses to a point.}
	\label{fig:aniso mcm}
\end{figure}

Using linearly stabilized schemes, the added factor of $\gamma(\omega) + \gamma''(\omega)$ presents no additional difficulty.   Again, we can stabilize with $p\Delta u$, setting $p=(1+(m^2-1)/(m^2+1))\bar{p}_{\text{min}}$. Shown in Fig.\ \ref{fig:aniso mcm} is an example with $m=4$. The solution is generated on a $256\times 256$ periodic grid using linearly stabilized ETDRK2 with 500 time steps.

\section{Conclusion}
\label{sect: conc}
In this work, we have identified three properties critical for selecting effective linearly stabilized schemes: unconditional stability over an unbounded $p$-parameter range, strong damping, and low polynomial degree in $p$ of the coefficients in the error expansion (relative to the order of the method).

We have proposed a number of new methods based on IMEX multistep methods and exponential Runge-Kutta methods.
Of the second order IMEX methods,  SBDF2 was often effective and possessed superior damping to CNAB.  
On the other hand, CNAB had the advantage of producing small errors in certain problems when applying a small time step  $\Delta t$. 
We found no viable third or higher order IMEX method; all studied possessed only a bounded $p$-parameter range. 
This limitation was removed by considering ETDRK schemes.  We considered ETDRK2 and ETDRK4;
these schemes outperformed the multistep-based methods in certain problems involving small $p$-values.




Of the pre-existing linearly stabilized methods, neither was optimal in our test problems. 
SBDF1 is only first order accurate. The EIN method, although formally second order accurate, exhibited a reduced order of accuracy in many of our numerical experiments due to its error coefficients being high degree polynomials in $p$. These shortcomings were examined in Sect.~\ref{sssect: ein loss of accuracy} and \ref{sssect: shrinking dumbbell} where substantial improvements in accuracy and efficiency were made by using our new methods.

A number of questions have been raised throughout our work that are worthy of further consideration. The derivation of third and higher order  methods excelling in all three of our criteria remains open. 
Moreover, as non-periodic boundary conditions are somewhat complicated for exponential time differencing methods, higher order methods that do not require the matrix exponential would be particularly compelling.
Adaptivity also could be investigated. Both the time step-size and the parameter $p$ are candidates for adaptivity in time, although doing so comes at the cost of carrying out matrix factorization at each time step. 
The analysis of test problem \eqref{ammc} also suggests adaptivity of $p$ in space may lead to interesting results. 


%
%


\begin{acknowledgements}
 \textcolor{black}{We are grateful to the referees for their constructive input.}
 \end{acknowledgements}

\section*{Declarations}
{\small\textbf{Funding}\, 
	The authors gratefully acknowledge the financial support of NSERC Canada
	(RGPIN 2016-04361).}
\\[6pt]
{\small \textbf{Conflict of interest}\,
	The authors have no conflicts of interest to declare that are relevant to the content of this article.} 
\\[6pt]
{\small \textbf{Availability of data and material}\,
	Data sharing not applicable to this article as no datasets were generated or analysed during the current study.}
\\[6pt]
{\small \textbf{Code availability}\, 
	Codes used during the current study are available upon reasonable request.}

\bibliography{master_ref}
\bibliographystyle{spmpsci}      
\end{document}